\documentclass[10pt,a4paper,oneside]{amsart}

\newcounter{commentcounter}

\usepackage{amsmath} 
\usepackage{amssymb} 
\usepackage{amsthm} 
\usepackage{stmaryrd} 
\usepackage[english]{babel} 
\usepackage[font=small,justification=centering]{caption} 
\usepackage[nodayofweek]{datetime}
\usepackage[shortlabels]{enumitem} 
\usepackage[T1]{fontenc} 
\usepackage[utf8]{inputenc} 
\usepackage{ifthen} 
\usepackage{mathabx} 
\usepackage{mathtools} 
\usepackage[dvipsnames]{xcolor} 
\usepackage[pdftex,  colorlinks=true, pagebackref=true]{hyperref} 
    \hypersetup{urlcolor=RoyalBlue, linkcolor=RoyalBlue,  citecolor=black}
\usepackage{setspace} 
\usepackage{tikz-cd} 
\usepackage{xfrac} 
\usepackage[capitalize]{cleveref} 
\usepackage{booktabs}
\usepackage{multirow,multicol}
\renewcommand*{\backref}[1]{}
\renewcommand*{\backrefalt}[4]
{
    \ifcase #1
        No citation in the text.
    \or
        Cited on Page #2.
    \else
        Cited on Pages #2.
    \fi
}

\makeatletter
\@namedef{subjclassname@1991}{Mathematical subject classification 1991}
\@namedef{subjclassname@2000}{Mathematical subject classification 2000}
\@namedef{subjclassname@2010}{Mathematical subject classification 2010}
\@namedef{subjclassname@2020}{Mathematical subject classification 2020}
\makeatother

\newtheorem{thm}{Theorem}[section]
\newtheorem{lemma}[thm]{Lemma}
\newtheorem{corollary}[thm]{Corollary}
\newtheorem{prop}[thm]{Proposition}
\newtheorem{conjecture}[thm]{Conjecture}

\newtheorem*{lemma*}{Lemma}

\newtheorem{thmx}{Theorem}

\theoremstyle{definition}
\newtheorem{defn}[thm]{Definition}
\newtheorem{remark}[thm]{Remark}

\newtheorem{example}[thm]{Example}

\theoremstyle{plain}

    \newtheoremstyle{TheoremNum}
        {8pt}{8pt} 
        {\itshape} 
        {-0.15cm} 
        {\bfseries} 
        {.} 
        { }  
        {\thmname{#1}\thmnote{ \bfseries #3}}
    \theoremstyle{TheoremNum}
    \newtheorem{duplicate}{}

\newcommand*{\claimproofname}{My proof}

\DeclareMathOperator{\Out}{\mathrm{Out}}

\DeclareMathOperator{\Ext}{\mathrm{Ext}}
\DeclareMathOperator{\Tor}{\mathrm{Tor}}

\DeclareMathOperator{\Fix}{\mathrm{Fix}}
\DeclareMathOperator{\pseudorank}{\mathrm{pseudo-rank}}

\DeclareMathOperator{\sym}{Sym}
\DeclareMathOperator{\alt}{Alt}

\DeclareMathOperator{\st}{\mathrm{st}}
\DeclareMathOperator{\lk}{\mathrm{lk}}
\newcommand{\CF}{\mathcal{CF}}
\newcommand{\CFsol}{\CF_{\mathrm{sol}}}
\newcommand{\CFmax}{\CF_{\mathrm{max}}}

\newcommand{\CFp}{\CF_p}

\newcommand{\cala}{{\mathcal{A}}}

\newcommand{\calc}{{\mathcal{C}}}

\newcommand{\calg}{{\mathcal{G}}}
\newcommand{\calh}{{\mathcal{H}}}

\newcommand{\call}{{\mathcal{L}}}

\newcommand{\calw}{{\mathcal{W}}}

\newcommand{\PGL}{\mathrm{PGL}}

\newcommand{\CAT}{\mathrm{CAT}}

\newcommand{\onto}{\twoheadrightarrow}

\newcommand{\FA}{\mathsf{FA}}

\newcommand{\vb}{\mathsf{vb}}

\DeclareMathOperator{\rank}{\mathrm{rank}}

\DeclareMathOperator{\cd}{\mathrm{cd}}
\DeclareMathOperator{\vcd}{\mathrm{vcd}}

\DeclareMathOperator{\Ends}{\mathrm{ends}}
\DeclareMathOperator{\cc}{\mathrm{cc}}

\def\Z{\mathbb{Z}}

\newcommand{\NN}{\mathbb{N}}
\newcommand{\ZZ}{\mathbb{Z}}
\newcommand{\CC}{\mathbb{C}}
\newcommand{\RR}{\mathbb{R}}
\newcommand{\QQ}{\mathbb{Q}}
\newcommand{\FF}{\mathbb{F}}

\usepackage{tikz}
\usetikzlibrary{arrows,quotes}
\usetikzlibrary{shapes.geometric}
\tikzstyle{blackNode}=[fill=black, draw=black, shape=circle]

\title{Profinite properties of Coxeter groups}

\author{Sam Hughes}
\address{Sam Hughes\\
Mathematisches Institut, Rheinische Friedrich-Wilhelms-Universit\"at Bonn, Endenicher Allee 60, 53115 Bonn, Germany}
\email{hughes@math.uni-bonn.de; sam.hughes.maths@gmail.com}

\author{Philip M\"oller}
\address{Philip M\"oller\\ Institute of Mathematics, Heinrich-Heine-University Düsseldorf, Universitätsstra{\upshape{\ss}}e 1, 40225, Düsseldorf, Germany}
\email{philip.moeller.2@hhu.de}

\author{Olga Varghese}
\address{Olga Varghese\\ Institute of Mathematics, University of Münster, Einsteinstra{\upshape{\ss}}e 62, 48149, Münster, Germany}
\email{olga.varghese@uni-muenster.de}

\address{Sam P.~Fisher\\
Mathematical Institute, University of Oxford, Andrew Wiles Building, Woodstock Rd, Oxford OX2 6GG, UK}
\email{fisher@maths.ox.ac.uk}

\makeatletter
\let\@wraptoccontribs\wraptoccontribs
\makeatother
\contrib[\MakeLowercase{with an appendix by}]{Sam P. Fisher \MakeLowercase{and} Sam Hughes}

\date{\today}


\begin{document}
	
\pagenumbering{arabic}
	
\begin{abstract}
	We prove a number of results about profinite completions of Coxeter groups.  For example we prove Coxeter groups are good in the sense of Serre and that various splittings of Coxeter groups arising from actions on trees are detected by the profinite completion.  As an application we prove a number of families of Coxeter groups are profinitely rigid amongst Coxeter groups.  We also prove that Gromov-hyperbolic FC type, extra large type, and odd Coxeter groups are almost profinitely rigid amongst Coxeter groups.  In the appendix, Sam Fisher and Sam Hughes show that the Atiyah Conjecture holds for all Coxeter groups, and that $\ell^2$-Betti numbers and their positive characteristic analogues are profinite invariants of Coxeter groups and of virtually compact special groups.
\end{abstract}

\maketitle

\section{Introduction}
\subsection*{Context}
In 1911, Max Dehn defined the word, conjugacy, and isomorphism problems for finitely presented groups.  Whilst all three problems are unsolvable in full generality, for Coxeter groups, both the word \cite{Tits1969} and conjugacy problems \cite{AlonsoBridson1995,Krammer2009} have been solved.  Despite much effort, the isomorphism problem amongst Coxeter groups remains open, see \cite{Bahls2005,Muhller2006} for a survey and \cite{CharneyDavis2000,Mihalik2007,CapracePrzytycki2010,HowlettMuhlherrNuida2018,SantosRegoSchwer2022} and the references therein for the recent progress.  Whilst Coxeter groups occur very naturally as abstractions of reflection groups, the uninitiated reader may view them as more combinatorial objects.  To this end, let $\Gamma$ be a finite simplicial graph with edge labelling $m\colon E(\Gamma) \to \NN_{\geq 2}$.  The \emph{Coxeter group} $W_\Gamma$ is the group defined by the presentation
\[ \left\langle V(\Gamma)\ \big|\ v^2 \text{ for }v\in V(\Gamma),\ (vw)^{m(\{v,w\})}\text{ for }\{v,w\}\in E(\Gamma)\right\rangle. \]
The key difficulty in resolving the isomorphism problem is that many distinct labelled graphs may correspond to the same isomorphism class of Coxeter groups.

One of the most na\"ive ways to try to understand the nature of an infinite group $G$ is via its actions on finite sets.  This essentially amounts to understanding the finite quotients of a group.  Of course to be able to gain any useful information one must insist that every element acts non-trivially on some finite set, when this happens we say $G$ is \emph{residually finite}.  Pursuing this direction further, one is led to the study of profinite rigidity---the study of $G$ through its profinite completion $\widehat G$, for an excellent introduction the reader is referred to Reid's ICM notes \cite{Reid2018}.  

Let $\calc$ be a class of finitely generated residually finite groups.  The \emph{$\calc$-genus} of a group $G$ is the set of isomorphism classes
\[\calg_\calc(G)\coloneqq \{H\ |\ H\in\calc,\ \widehat H\cong \widehat G \}/\cong.\]
If $|\calg_\calc(G)|=1$, then we say $G$ is \emph{profinitely rigid in $\calc$} or $\calc$\emph{-profinitely rigid}.  If $\calg_\calc(G)$ is finite, then we say $G$ is \emph{almost profinitely rigid in $\calc$}. The relevance to the isomorphism problem is the following well known fact which has appeared in work of Gardam \cite[Proposition~8]{Gardam2019} and Bridson, Conder, and Reid \cite[page~5]{BridsonConderReid2016}.

\begin{lemma*}
    Let $\calc$ be a class of finitely presented residually finite groups.  If $G$ in $\calc$ has $|\calg_\calc(G)|=1$, then the isomorphism problem for $G$ is solvable in $\calc$.
\end{lemma*}

The study of profinite rigidity in the context of Coxeter groups has been an active topic for almost a decade, we briefly survey the literature: \begin{enumerate}
    \item Kropholler and Wilkes showed that right-angled Coxeter groups are profini-\\tely rigid amongst themselves \cite{KrophollerWilkes2016}.  Later, Corson along with the authors of this paper showed that right-angled Coxeter groups are profinitely rigid amongst \emph{all} Coxeter groups \cite{CorsonHughesMollerVarghese2023}.  On the other hand they showed that the genus of some right-angled Coxeter groups can be infinite amongst all finitely generated residually finite groups, building on techniques from \cite{PlatonovTavgen1986,BridsonGrunewald2004}.
    \item Kropholler and Wilkes also showed that the profinite completion of a right-angled Coxeter group splits as a profinite free product if and only if the right-angled Coxeter group itself splits as a free product \cite[Theorem~11]{KrophollerWilkes2016}.  In \cite[Theorem~3.10]{CorsonHughesMollerVarghese2023}, Corson along with the authors showed that $\calw$-profinite rigidity of free products of Coxeter groups reduces to the profinite rigidity of the free factors, where we denote by $\calw$ the class of all Coxeter groups.
    \item Bridson, McReynolds, Spitler, and Reid showed that 14 hyperbolic Coxeter triangle groups are profinitely rigid amongst all finitely generated residually finite groups \cite{BridsonMcReynoldsReidSpitler2021}, building on their breakthrough work on certain $3$-manifold groups \cite{BridsonMcReynoldsReidSpitler2020}.
    \item Santos Rego and Schwer showed that Coxeter groups on at most $3$ vertices are distinguished from each other by their profinite completions \cite[Theorem~4.25]{SantosRegoSchwer2022}, building on work of Bridson--Conder--Reid \cite{BridsonConderReid2016}.
    \item The second and third author showed that irreducible affine Coxeter groups are profinitely rigid amongst Coxeter groups \cite{MollerVarghese2023}, this was generalised to all finitely generated residually finite groups in \cite{CorsonHughesMollerVarghese2024}, see also \cite{PaoliniSklinos2024}.
\end{enumerate}  

\bigskip

\subsection*{Rigidity}
Our first main result gives a number of families of Coxeter groups where we can prove profinite rigidity or almost rigidity amongst Coxeter groups.

\bigskip

\begin{thmx}\label{thmx rigidity}
    The following Coxeter groups are profinitely rigid amongst Coxeter groups:
    \begin{enumerate}[label=(1\alph*)]
        \item reflection groups of regular hyperbolic polygons (\Cref{triangleAndmore}(3));
        \item cocompact hyperbolic simplicial reflection groups (\Cref{thm.Lanner});
        \item virtually abelian Coxeter groups (\Cref{virtually.abelian.Coxetergroups.profinitely.rigid});
        \item odd forest Coxeter groups (\Cref{oddForest_rigid});
        \item Coxeter groups such that all edge-labels are divisible by $4$ (\Cref{thm:even.profinite});
        \item complete Coxeter groups such that all edge-labels are equal to $n\neq 4k+2$, for $k\geq 1$ (\Cref{complete n neq 3 4k+2});
        \item Coxeter groups of rank at most $3$ (\Cref{rank at most 3 rigid});
        \item Coxeter groups of rank $4$ such that all edge-labels are equal to $n\neq 4k+2$, for $k\geq1$ (\Cref{CoxeterGroups4Vertices}).
    \end{enumerate}
    The following Coxeter groups are almost profinitely rigid amongst Coxeter groups:
    \begin{enumerate}[label=(2\alph*)]
        \item Gromov-hyperbolic Coxeter groups of FC type (\Cref{almost hyp FC type});
        \item virtually surface Coxeter groups (\Cref{triangleAndmore}(1));
        \item Gromov-hyperbolic Coxeter groups of odd type (\Cref{OddCoxeterAlmost});
        \item odd Coxeter groups (\Cref{OddCoxeterAlmost}); 
        \item Coxeter groups of extra large type (\Cref{XL almost}).
    \end{enumerate}
\end{thmx}

We now briefly explain the terms in the theorem.  A Coxeter group is a \emph{reflection group of a regular hyperbolic polygon} if its Coxeter graph is an $n$-gon and all labels are equal.  A Coxeter group is a \emph{cocompact hyperbolic simplicial reflection group} if it acts cocompactly on real hyperbolic $n$-space with fundamental domain a simplex, such groups were classified by Lann\'er \cite{Lanner1950} and Vinberg \cite{Vinberg1985} (see also \Cref{rank 3 lanner} and \Cref{tab rank 5 Lanner}). The \emph{rank} of the Coxeter group $W_\Gamma$ is the number of vertices in the defining graph $\Gamma$. 
A Coxeter group $W_\Gamma$ is \emph{complete} if $\Gamma$ is a complete graph.  
A Coxeter group $W_\Gamma$ is of \emph{FC type} if every clique in $\Gamma$ generates a finite subgroup.  
A Coxeter group $W_\Gamma$ is \emph{odd} if every edge label of $\Gamma$ is odd.  
A Coxeter group $W_\Gamma$ is of \emph{extra large type} if every edge label is at least $4$.
A group $G$ is \emph{virtually free/surface/abelian} if it contains a finite index subgroup isomorphic to a free group/surface group/abelian group.

\subsection*{Product decompositions}
A Coxeter group $W_\Gamma$ can be decomposed using the combinatorial structure of $\Gamma$ in the form $W_\Gamma\cong W_{\rm sph}\times W_{\rm aff}\times W_{\rm gen}$, where $W_{\rm sph}$ is trivial or a finite Coxeter group called the \emph{spherical part} of $W_\Gamma$, $W_{\rm aff}$ is trivial or a product of irreducible affine Coxeter groups called the \emph{affine part} of $W_\Gamma$, and $W_{\rm gen}$ is trivial or an infinite non-affine Coxeter group called the \emph{generic part} of $W_\Gamma$.

\smallskip
\begin{duplicate}[\Cref{thmx:products}]
Let $W_\Gamma=W_{\Gamma_{{\rm sph}}}\times W_{\Gamma_{{\rm aff}}}\times W_{\Gamma_{{\rm gen}}}$ and $W_\Omega=W_{\Omega_{{\rm sph}}}\times W_{\Omega_{{\rm aff}}}\times W_{\Omega_{{\rm gen}}}$ be Coxeter groups such that $ \widehat{W_\Gamma}\cong \widehat{W_\Omega}$. Then \[W_{\Gamma_{{\rm sph}}}\cong W_{\Omega_{{\rm sph}}},\quad 
W_{\Gamma_{{\rm aff}}}\cong W_{\Omega_{{\rm aff}}},\quad \text{and} \quad \widehat{W_{\Gamma_{{\rm gen}}}}\cong \widehat{W_{\Omega_{{\rm gen}}}}.\]
In particular,  a Coxeter group $W_\Gamma$ is profinitely rigid amongst Coxeter groups if and only if $W_{\Gamma_{\rm gen}}$ is profinitely rigid amongst Coxeter groups.
\end{duplicate}

We remark that there is a long history of results which show that direct product decompositions of profinite completions do not pass to dense subgroups.  Indeed, the earliest examples of Grothendieck pairs by Platonov--Tavgen \cite{PlatonovTavgen1986} are pairs of groups $\iota\colon P\rightarrowtail F_n^2$, for $n\geq 4$, where $\widehat P\cong\widehat {F_n^2}$  and $P$ is not finitely presented, nor a direct product of two infinite groups.  Other examples have appeared in the work of Bridson--Grunewald \cite{BridsonGrunewald2004}, Bridson \cite{Bridson2016,Bridson2024FbyF, Bridson2024}, and in work of the authors with Corson in the context of Coxeter groups \cite{CorsonHughesMollerVarghese2023}.  Note that our earlier work with Corson does not contradict the above theorem since we provide pairs $P\rightarrowtail W$ where $P$ is not finitely presented, and hence not a Coxeter group.

\subsection*{Invariants}
Our next main result records a large number of properties of Coxeter groups that are profinite invariants amongst Coxeter groups.

\setcounter{thmx}{2}
\begin{thmx}\label{thmx props}
    Let $W_\Gamma$ and $W_\Lambda$ be Coxeter groups and suppose $\widehat {W_\Gamma}\cong \widehat{W_\Lambda}$.  Then,
    \begin{enumerate}
        \item $W_\Gamma$ has $\FA$ if and only if $W_\Lambda$ has $\FA$ (\Cref{{profiniteinvariantFA}});
        \item $|\Ends(W_\Gamma)|=|\Ends(W_\Lambda)|$ (\Cref{profiniteinvariantNumberofends});
        \item $W_\Gamma$ is Gromov-hyperbolic if and only if $W_\Lambda$ is  Gromov-hyperbolic (\Cref{profiniteinvariantHyperbolicity});
        \item $W_\Gamma$ is virtually free if and only if $W_\Lambda$ is virtually free (\Cref{vir surface profinite});  
        \item $W_\Gamma$ is virtually surface if and only if $W_\Lambda$ is virtually surface (\Cref{vir surface profinite});
        \item $W_\Gamma$ is of FC type if and only if $W_\Lambda$ is of FC type (\Cref{profinite.invariant.FC});
        \item $W_\Gamma$ is odd if and only if $W_\Lambda$ is odd (\Cref{props odd}(3));
        \item $\CFp(W_\Gamma)=\CF_p(W_\Lambda)$ for every prime $p$ (\Cref{prop:CF}(1));
        \item $H^\ast(W_\Gamma;\FF_p) \cong H^\ast(W_\Lambda;\FF_p)$ (\Cref{profInv_X_H*}(1));
        \item $\chi(W_\Gamma)=\chi(W_\Lambda)$ (\Cref{profInv_X_H*}(2));
        \item $M(W_\Gamma)\cong M(W_\Lambda)$, where $M(-)$ denotes the Schur multiplier (\Cref{M(W) prof}).
    \end{enumerate}
\end{thmx}

We now explain the terms in the theorem, relevant literature, and some of the methods of proof.  We remark that the question of which invariants of a group are profinite is extremely subtle, see for instance \cite{Bridson2024FbyF}.

A group $G$ has \emph{Serre's fixed point property $\FA$} if every action of $G$ on a simplicial tree without edge inversions has a global fixed point.  Amongst the class of finitely generated residually finite groups, property $\FA$ is not a profinite invariant \cite{CheethamWestLubotzkyReidSpitler2022}, so (1) is very much a special property of Coxeter groups.  We refer the reader to \cite{Aka2012,CottonBarratt2013,Bridson2024} for other results in this vein.  Our proofs of (1) and (6) use profinite Bass-Serre theory developed in \cite{MelnikovZalesskii1989,RibesZalesskii2010}.

The \emph{ends} of a finitely generated group $G$, denoted $\Ends(G)$, is roughly `the number of connected components at infinity' in a finitely generated Cayley graph for $G$.  A classical theorem of Hopf \cite{Hopf44} states that the number of ends of $G$ is either $0$, $1$, $2$, or $\infty$; depending on if $G$ is finite, one-ended, virtually-$\Z$, or splits as a non-trivial graph of groups with finite edge groups. As the previous paragraph suggests, detecting properties of groups acting on trees profinitely is remarkably subtle.  It appears to be completely open whether the number of ends of a group is a profinite invariant.  However, work of Cotton-Barrett suggests that it is likely to be false \cite{CottonBarratt2013}.  In particular, it seems possible that (2) above is a special property of Coxeter groups.

Our proof that Gromov-hyperbolicity is a profinite invariant amongst Coxeter groups depends on a number of results.  Firstly, we use Moussong's characterisation of Gromov-hyperbolic Coxeter groups \cite{Moussong88}, secondly we use deep work of Wilton and Zalesskii on the profinite completions of virtually compact special groups \cite{WiltonZalesskii2017}, and thirdly we use the fact that Gromov-hyperbolic Coxeter groups are virtually compact special which follows from work of Niblo--Reeves \cite{NibloReeves2003}, Williams \cite{Williams99}, and Haglund--Wise \cite{HaglundWise2010}.

The results (4) and (5) fit into a lineage of results proving profinite rigidity of free and surface groups amongst various classes of groups, see for instance: Bridson--Conder--Reid's work on Fuchsian groups \cite{BridsonConderReid2016}, Wilton's work on limit groups \cite{Wilton2021} (see also \cite{FruchterMorales2022} for residually free groups), \cite{HLIPSV2025} for K\"ahler groups, and \cite{Jaikin2023} for other recent progress.  We remark that our methods are most similar in spirit to that of Wilton's \cite{Wilton2021} and to that of Fruchter--Morales' \cite{FruchterMorales2022}, although we also employ work of Gordon, Long, and Reid \cite{GordonLongReid2004} on surface subgroups of Coxeter groups.

For a group $G$, we denote by $\CFp(G)$ the poset of conjugacy classes of finite $p$ subgroups of $G$.  We denote by $\chi(G)$, the Euler characteristic of $G$, see \cite[IX.7]{Brown1982}.  Items (7)-(11) of \Cref{thmx props} are essentially consequences of the claim that Coxeter groups are \emph{good in the sense of Serre}, see \Cref{sec good} for both a definition and proof of this fact.  We mention here that our proof relies on work of Genevois \cite{Genevois2024}.  Item (8), however, is a direct application of work of Boggi--Zalesskii \cite{BoggiZalesskii2024} building on work of Minasyan--Zalesskii \cite{MinasyanZalesskii2016}.  Item (11) depends on work of Howlett \cite{Howlett1988}, which computes the Schur multipliers of Coxeter groups.

\subsection*{Reduction to the one-ended case} 
A very common theme in group theory is the aim to decompose groups into simpler pieces. In geometric group theory, a typical approach is to try to find and understand so-called JSJ-decompositions of a group. In general, it is not clear, whether the profinite completion of a group can detect JSJ-decompositions. However, profinite Bass-Serre theory allows us to ``almost'' detect these decompositions, reducing the question of almost $\mathcal{W}$-profinite rigidity for all Coxeter groups to the one-ended case.

\begin{duplicate}[\Cref{thmx:Reduction1ended}]
    Coxeter groups are almost $\mathcal{W}$-profinitely rigid if and only if $1$-ended Coxeter groups are almost $\mathcal{W}$-profinitely rigid.
\end{duplicate}

\subsection*{Questions}
We leave open a few questions regarding profinite properties of Coxeter groups.  We say that $H\leqslant G$ is a \emph{retract} of $G$ if there exists an epimorphism $r\colon G \twoheadrightarrow H$ which is the identity on $H$. We say that $H$ is a \emph{virtual retract} of $G$ if there exists $G'$ a finite index subgroup of $G$ such that  $H$ is a retract of $G'$.  We say that $H$ is \emph{virtually a virtual retract} of $G$ if a finite index subgroup of $H$ is a virtual retract of $G$.
\emph{
    \begin{enumerate}
        \item If $W_\Gamma$ is a Coxeter group and $W_\Omega\leqslant W_\Gamma$ is a proper parabolic subgroup, does $W_\Omega$ have the full profinite topology?
        \begin{enumerate}
            \item Is $W_\Omega$ virtually a virtual retract of $W_\Gamma$?
            \item What about a virtual retract?
        \end{enumerate}
        \item Are virtually surface Coxeter groups profinitely rigid amongst Coxeter groups?\item  Are non-Euclidean crystallographic groups, that is discrete subgroups of $\PGL_2(\RR)$, distinguished from each other up to isomorphism by their profinite completions?  
    \end{enumerate}
}
Note that the class in (2) includes all hyperbolic Coxeter groups with diagram an $n$-gon. 
 The analogous result for Fuchsian groups is \cite{BridsonConderReid2016}.

\subsection*{Structure of the paper}
In \Cref{sec prelims} we recount the relevant background on Coxeter groups and profinite completions that we will need.  

In \Cref{sec good} we prove that Coxeter groups are good in the sense of Serre (\cref{CoxGpGood}) and then deduce a number of applications.  Most notably, in \Cref{sec good Schur} we prove that the Schur multiplier is a profinite invariant of Coxeter groups and in \Cref{sec good CF} we explain how to apply work of Boggi and Zalesskii \cite{BoggiZalesskii2024} on the poset of finite subgroups up to conjugacy to Coxeter groups.

In \Cref{sec product decomps} we prove the product decomposition theorem (\Cref{thmx:products}).

In \Cref{sec low vcd} we investigate Coxeter groups with virtual cohomological dimension equal to one (\Cref{sec low vcd vir free}) and two (\Cref{sec low vcd vir surface}).

In \Cref{sec trees} we study the actions of Coxeter groups on profinite trees.  In \Cref{sec trees FA FC} we prove that property $\FA$ is a profinite invariant of Coxeter groups (\Cref{profiniteinvariantFA}) and that being of FC type is a profinite invariant of Coxeter groups (\Cref{profinite.invariant.FC}).  In \Cref{sec trees ends class A} we prove that Coxeter groups are in de Bessa, Porto, and Zalesskii's class $\cala$ \cite{AndersonBessaZalesskii2023} (\Cref{ClassA}) and deduce a number of consequences. 
 Most notably, we prove that among Coxeter groups the number of ends is a profinite invariant (\Cref{profiniteinvariantNumberofends}), and we also prove \Cref{thmx:Reduction1ended}.

In \Cref{sec lanner} we prove that cocompact hyperbolic simplicial reflection groups are profinitely rigid amongst Coxeter groups.

In \Cref{sec Special types} we prove rigidity and almost rigidity for a number of classes of Coxeter groups defined via combinatorics of the defining graphs.  This establishes the remaining cases of \Cref{thmx rigidity}.

In \Cref{sec appendix} we prove the Strong Atiyah Conjecture for Coxeter groups.  We also prove that $\ell^2$-Betti numbers (and certain modulo $p$ analogues) are profinite invariants amongst good residually (locally indicable and amenable) groups with sufficient finiteness properties.  As applications, using work of Fisher \cite{Fisher_Improved,Fisher_freebyZ} we deduce profinite invariance of virtual homological fibring for (virtually) RFRS groups and consider other consequences relating to free-by-cyclic groups in cohomological dimension $2$.  We also highlight that these results apply to Coxeter groups, see \Cref{app:cor.Cox}.

\subsection*{Acknowledgements}
Sam Hughes was supported by a Humboldt Research Fellowship at Universit\"at Bonn and by the Deutsche Forschungsgemeinschaft (DFG, German Research Foundation) under Germany's Excellence Strategy - EXC-2047/1 - 390685813. Olga Varghese was funded by the Deutsche Forschungsgemeinschaft (DFG, German Research Foundation) under Germany's Excellence Strategy EXC 2044 –390685587, Mathematics Münster: Dynamics–Geometry–Structure.  

The authors are grateful to Kevin Schreve for pointing out a mistake in an earlier version of this paper and in \cite{Genevois2024}.  The authors thank Anthony Genevois for helpful conversations.

\section{Preliminaries}\label{sec prelims}

\subsection{Background on Coxeter groups}
We begin by recalling definitions of a Coxeter graph and the associated Coxeter group.
\begin{defn}
A \emph{Coxeter graph} $\Gamma$ is a finite simplicial graph with the vertex set $V(\Gamma)$, the edge set $E(\Gamma)$ and with an edge-labeling $m\colon E(\Gamma)\rightarrow\mathbb{N}_{\geq 2}$.  The \emph{Coxeter group} $W_\Gamma$ associated to $\Gamma$ is the group with the presentation 
$$W_\Gamma=\left\langle V(\Gamma)\ \big| \ v^2\text{ for }v\in V(\Gamma), (vw)^{m(\left\{v,w\right\})} \text{ for }\left\{v,w\right\}\in E(\Gamma)\right\rangle.$$ 
A \emph{Coxeter system} $(W,S)$ is a pair consisting of a Coxeter group $W$ and a generating set $S$ for $W$ which is precisely the vertex set of a Coxeter graph $\Gamma$ such that $W\cong W_\Gamma$.
\end{defn}
We denote by $\calw$ the class consisting of all Coxeter groups.
\bigskip

\begin{figure}[htb]
 	\begin{center}
	\captionsetup{justification=centering}
 		\begin{tikzpicture}
 			\draw[fill=black]  (0,0) circle (2pt);
             \draw[fill=black]  (1,0) circle (2pt);
             \draw[fill=black]  (0.5,0.7) circle (2pt);
            \draw (0,0)--(1,0);
            \draw (1,0)--(0.5, 0.7);
             \draw (0,0)--(0.5, 0.7);
             \node at (0.5, -0.25) {$p$};
             \node at (0.1, 0.45) {$r$};
             \node at (0.9, 0.45) {$q$};
         \end{tikzpicture}
         \caption{$\Delta(p,q,r)$.}
         \label{fig:triangle}
     \end{center}
 \end{figure}
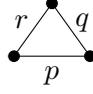

Let $\Gamma$ be a Coxeter graph. A Coxeter group $W$ is said to be of \emph{type} $\Gamma$ if $W\cong W_\Gamma$. For $p,q,r\in\NN_{\geq 2}$ we denote by $\Delta(p,q,r)$ the Coxeter graph that is isomorphic to the graph in \Cref{fig:triangle}.
 
It is straightforward to see that $W_{\Delta(2,2,3)}$ is isomorphic to $W_\Omega$ where $\Omega$ is of type $(\left\{v,w\right\}, \left\{\left\{v,w\right\}\right\})$ and $m(\left\{v,w\right\})=6$. In \cite{MihalikRatcliffeTschantz2007}, Mihalik, Ratcliffe and Tschantz describe an algorithm for constructing from any set of Coxeter generators of $W$ a set of Coxeter generators $R$ of maximal cardinality. This observation leads to the following definition.

\begin{defn}
Let $W_\Gamma$ be a Coxeter group. 
\begin{enumerate}
\item The \emph{rank} of $W_\Gamma$ is defined as: $\rank(W_\Gamma):=|V(\Gamma)|$.
\item The \emph{pseudo-rank} of $W_\Gamma$ is defined as:
\[\pseudorank(W_\Gamma):=\max\left\{\rank(W_\Omega)\mid \Omega \text{ a Coxeter graph with } W_\Omega\cong W_\Gamma\right\}. \]
\end{enumerate}
\end{defn}

Note that the pseudo-rank is always finite by \cite{MihalikRatcliffeTschantz2007}.

\begin{defn}
A Coxeter group $W_\Gamma$ is called \emph{graph rigid} if it cannot be defined by two non-isomorphic graphs.
\end{defn}

\begin{lemma}~\cite[Theorem~3.3]{Nuida2006}\label{DihedralGroups.Graphrigid}
Let $\Gamma_n=(\left\{v,w\right\}, \left\{\left\{v,w\right\}\right\})$ be a Coxeter graph with the edge-labeling $m(\left\{v,w\right\})=n$, $n\in\NN_{\geq 2}$. The Coxeter group $W_{\Gamma_n}$ is isomorphic to the Dihedral group $D_n$ of order $2n$ and is graph rigid if and only if $n\neq 4k+2$ for $k\geq 1$. If $n=4k+2$, then $W_{\Gamma_n}\cong W_\Omega$ where $\Omega\cong\Gamma_n$ or $\Omega\cong\Delta(2,2,2k+1)$.
\end{lemma}

Some examples of graph rigid Coxeter groups are collected in the next proposition.
\begin{prop}\label{graph.rigid.Coxetergroups}
Let $W_\Gamma$ be a 
group. If
\begin{enumerate}
\item $\Gamma$ is an $n$-gon where $n\geq 3$ and $W_\Gamma$ is infinite or
\item $m(E(\Gamma))\subseteq\left\{2\right\}\cup 4\NN$ or
\item $W_\Gamma$ is an irreducible affine Coxeter group,
\end{enumerate}
then $W_\Gamma$ is graph rigid.
\end{prop}
\begin{proof}
Items (1) and (3) follow from \cite[Main Theorem]{CharneyDavis2000} coupled with the observation that a Coxeter group with defining graph an $n$-gon is capable of acting effectively, properly and cocompactly on $\mathbb{E}^2$ or $\mathbb{H}^2$, see the proof of \cite[Theorem 2.1]{GordonLongReid2004}. 
Item (3) is precisely \cite[Proposition 17.7]{HowlettMuhlherrNuida2018} and Item (2) is precisely \cite[Theorem 4.11]{Radcliffe2001}.
\end{proof}

\subsection{Coxeter--Dynkin diagrams}
A \emph{Coxeter--Dynkin diagram} is a finite simplicial graph $\Gamma$ with vertices $V(\Gamma)$, edges $E(\Gamma)$ and the edge-labelling $m\colon E(\Gamma)\to\mathbb{N}_{\geq 3}\cup\left\{\infty\right\}$. The Coxeter group $W_\Gamma$ is defined as follows:
 $$W_{\Gamma}=\left\langle V(\Gamma)\middle\vert \begin{array}{l} 
v^2\text{ for all }v\in V(\Gamma), (vw)^2 \text{ if }\left\{v,w\right\}\notin E(\Gamma),\\ (vw)^{m(\left\{v,w\right\})}\text{ if } \left\{v,w\right\}\in E(\Gamma) \text{ and } m(\left\{v,w\right\})<\infty\end{array} \right\rangle.$$

Note that no edge here corresponds to a commutation relation and an unlabelled edge is assumed to have the label $3$.
 
Let $\Gamma$ be a Coxeter--Dynkin diagram and $W_\Gamma$ the corresponding Coxeter group. By definition, $W_\Gamma$ is called \emph{irreducible} if $\Gamma$ is connected. 

\begin{figure}[htb]
	\begin{center}
	\captionsetup{justification=centering}
		\begin{tikzpicture}[scale=0.8]
			\draw[fill=black]  (0,0) circle (2pt);
			\draw[fill=black]  (1,0) circle (2pt);
            \draw[fill=black]  (2,0) circle (2pt);
			\draw (0,0)--(1,0);
            \draw[dashed] (1,0)--(2,0);
            \draw[fill=black] (3,0) circle (2pt);
            \draw (2,0)--(3,0);
			\node at (-0.9,0) {$\mathtt{A}_{n(n\geq 1)}$};

            \draw[fill=black]  (0,-1.5) circle (2pt);
			\draw[fill=black]  (1,-1.5) circle (2pt);
            \draw[fill=black]  (2,-1.5) circle (2pt);
			\draw (0,-1.5)--(1,-1.5);
            \draw[dashed] (1,-1.5)--(2,-1.5);
            \draw[fill=black] (3,-1.5) circle (2pt);
            \draw (2,-1.5)--(3,-1.5);
            \node at (2.5, -1.3) {$4$};
            \node at (-0.9,-1.5) {$\mathtt{B}_{n(n\geq 2)}$};

            \draw[fill=black]  (0,-3) circle (2pt);
			\draw[fill=black]  (1,-3) circle (2pt);
            \draw[fill=black]  (2,-3) circle (2pt);
			\draw (0,-3)--(1,-3);
            \draw[dashed] (1,-3)--(2,-3);
            \draw[fill=black] (3,-2.3) circle (2pt);
            \draw[fill=black] (3,-3.7) circle (2pt);
            \draw (2,-3)--(3, -2.3);
            \draw (2,-3)--(3, -3.7);
            \node at (-0.9,-3) {$\mathtt{D}_{n(n\geq 4)}$};

            \draw (0,-4.5)--(4,-4.5);
            \draw (2,-3.5)--(2,-4.5);
            \draw[fill=black] (0,-4.5) circle (2pt); 
            \draw[fill=black] (1,-4.5) circle (2pt);
            \draw[fill=black] (2,-4.5) circle (2pt);
            \draw[fill=black] (3,-4.5) circle (2pt);
            \draw[fill=black] (4,-4.5) circle (2pt);
            \draw[fill=black] (2,-3.5) circle (2pt);
            \node at (-0.5,-4.5) {$\mathtt{E}_6$};

            \draw (0,-6)--(5,-6);
            \draw (2,-5)--(2,-6);
            \draw[fill=black] (0,-6) circle (2pt); 
            \draw[fill=black] (1,-6) circle (2pt);
            \draw[fill=black] (2,-6) circle (2pt);
            \draw[fill=black] (3,-6) circle (2pt);
            \draw[fill=black] (4,-6) circle (2pt);
            \draw[fill=black] (5,-6) circle (2pt);
            \draw[fill=black] (2,-5) circle (2pt);
            \node at (-0.5,-6) {$\mathtt{E}_7$};

            \draw (0,-7.5)--(6,-7.5);
            \draw (2,-6.5)--(2,-7.5);
            \draw[fill=black] (0,-7.5) circle (2pt); 
            \draw[fill=black] (1,-7.5) circle (2pt);
            \draw[fill=black] (2,-7.5) circle (2pt);
            \draw[fill=black] (3,-7.5) circle (2pt);
            \draw[fill=black] (4,-7.5) circle (2pt);
            \draw[fill=black] (5,-7.5) circle (2pt);
            \draw[fill=black] (6,-7.5) circle (2pt);
            \draw[fill=black] (2,-6.5) circle (2pt);
            \node at (-0.5,-7.5) {$\mathtt{E}_8$};

            \draw[fill=black] (8,0) circle (2pt);
            \draw[fill=black] (9,0) circle (2pt);
            \draw[fill=black] (10,0) circle (2pt);
            \draw[fill=black] (11,0) circle (2pt);
            \draw (8,0)--(11,0);
            \node at (9.5, 0.2) {$4$};
            \node at (7.5,0) {$\mathtt{F}_4$};

            \draw[fill=black] (8,-1.5) circle (2pt);
            \draw[fill=black] (9,-1.5) circle (2pt);
            \draw (8,-1.5)--(9,-1.5);
            \node at (8.5,-1.3)  {$6$};
            \node at (7.5,-1.5) {$\mathtt{G}_2$};

            \draw[fill=black] (8,-3) circle (2pt);
            \draw[fill=black] (9,-3) circle (2pt);
            \draw[fill=black] (10,-3) circle (2pt);
            \draw (8,-3)--(10,-3);
            \node at (9.5, -2.8){$5$};
            \node at (7.5,-3) {$\mathtt{H}_3$};

            \draw[fill=black] (8,-4.5) circle (2pt);
            \draw[fill=black] (9,-4.5) circle (2pt);
            \draw[fill=black] (10,-4.5) circle (2pt);
            \draw[fill=black] (11,-4.5) circle (2pt);
            \draw (8,-4.5)--(11, -4.5);
            \node at (10.5, -4.3){$5$};
            \node at (7.5,-4.5) {$\mathtt{H}_4$};

            \draw[fill=black] (8,-6) circle (2pt);
            \draw[fill=black] (9,-6) circle (2pt);
            \draw (8,-6)--(9,-6);
            \node at (8.5, -5.8){$m$};
            \node at (7.3,-6) {$\mathtt{I}_2(m)$};
				
		\end{tikzpicture}
	\caption{Coxeter--Dynkin diagram of type $\mathtt{X}_n$ where $\mathtt{X}_n$ has $n$ vertices.  Note that $m\geq3$ and $\mathtt{A}_2=\mathtt{I}_2(3),$ $\mathtt{B}_2=\mathtt{I}_2(4)$, 
    and $\mathtt{G}_2=\mathtt{I}_2(6)$.}
    \label{fig:FiniteCoxeterGroups}
	\end{center}
\end{figure}
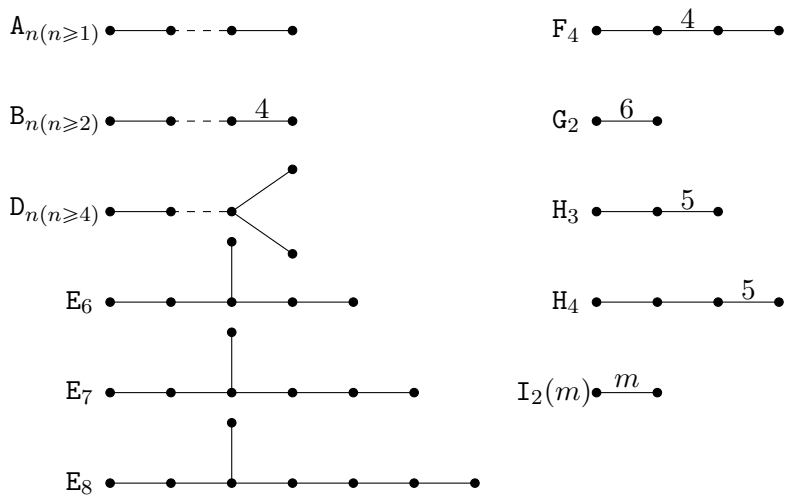

A classical result of Coxeter \cite{Coxeter1935} states that if $W_\Gamma$ is finite, then every connected component of the Coxeter--Dynkin diagram $\Gamma$ is isomorphic to one of the graphs in \Cref{fig:FiniteCoxeterGroups}.

Note that given a Coxeter graph $\Gamma$ it is straightforward to obtain a Coxeter-Dynkin diagram $\mathtt{\Omega}$ using $\Gamma$, in particular $|V(\Gamma)|=|V(\mathtt{\Omega})|$. Hence there is a canonical bijection between Coxeter graphs and Coxeter--Dynkin diagrams. 

\begin{remark}
    For the arguments in our proofs we usually work with Coxeter graphs. However, in some places classification results stated for Coxeter--Dynkin diagrams enter our arguments.  In this case we use letters of the form $\mathtt{A}$ to differentiate between the two conventions.
\end{remark}

\subsection{Parabolic subgroups}
Given a Coxeter graph $\Gamma$, for each $X\subseteq V(\Gamma)$, the subgroup which is generated by $X$ is canonically isomorphic to the Coxeter group $W_\Lambda$ where $\Lambda$ is the subgraph of $\Gamma$ induced by $X$, see \cite[Theorem 4.1.6]{Davis2008}. The subgroup $W_\Lambda$ is called \emph{special parabolic} and any conjugate of a special parabolic subgroup is called \emph{parabolic}. Note that instead of $W_{\Lambda}$, we usually just write $W_{V(\Lambda)}$. 

The following lemma follows readily from an analysis of the standard presentation of a Coxeter group. A standard reference for the right-angled case and more general graph products is \cite[Lemma~3.20]{Green1990}.

\begin{lemma}
\label{amalgam.decompostion.Coxeter}
Let $W_\Gamma$ be a Coxeter group. If there exist two vertices $v,w\in V(\Gamma)$ such that $\left\{v,w\right\}\notin E(\Gamma)$, then
    $$W_\Gamma\cong W_{\st(v)}*_{W_{\lk(v)}}W_{V(\Gamma)-\left\{v\right\}}$$ 
    where $\lk(v)=\left\{w\in V(\Gamma)\mid\left\{v,w\right\}\in E(\Gamma)\right\}$ and $\st(v)=\lk(v)\cup\left\{v\right\}$.
 \end{lemma}

\bigskip

\subsection{The poset of conjugacy classes of finite subgroups}
For a given group $G$ we denote by $\CF(G)$ the set of conjugacy classes of all finite subgroups; by $\CFsol(G)$ the set of conjugacy classes of all finite soluble subgroups; and by $\CFp(G)$ the set of conjugacy classes of all subgroups of $p$-power order in $G$. 
We define a partial order on $\mathcal{CF}(G)$ as follows: $[A]\leq [B]$ if there exists a $g\in G$ such that $A\subseteq gBg^{-1}$. We denote by $\CFmax(G)$ the set of maximal elements of $\CF(G)$.

Let $W_\Gamma$ be a Coxeter group. A subset $J\subseteq V(\Gamma)$ is called a \emph{spherical subset} of $V(\Gamma)$ if the special parabolic subgroup $W_{J}$ is finite.
\begin{prop}~\cite[Theorem 1.9]{BradyMcCammondMühlherrNeumann2002}
\label{FactsCoxeter}
Let $W_\Gamma$ be a Coxeter group. The conjugacy classes of the maximal finite subgroups of $W_\Gamma$ are in one-to-one correspondence with the maximal spherical subsets of $V(\Gamma)$.
\end{prop}
More precisely, we denote by $\CFmax(W_\Gamma)$ the set of conjugacy classes of maximal finite subgroups of $W_\Gamma$ and by $\mathcal{S}_{{\rm max}}(\Gamma)$ the set of maximal spherical subsets of $V(\Gamma)$. Then the map $\Phi\colon \mathcal{S}_{{\rm max}}(\Gamma)\to\CFmax(W_\Gamma)$ defined via $\Phi(I)=[W_I]$ is bijective.

In particular, let $J_1, \ldots, J_n$ be the maximal spherical subsets of $V(\Gamma)$. Then $\Gamma$ is a union of induced subgraphs $\Gamma_1, \ldots, \Gamma_n$ where $V(\Gamma_i)=J_i$ for $i=1,\ldots, n$. 

As an immediate consequence of \Cref{FactsCoxeter} we have the following result.
\begin{corollary}\label{graph.coveringbymaximal}
Let $W$ be a Coxeter group and suppose that $\CFmax(W)=\left\{[A_1], \ldots, [A_n]\right\}$. Then
\begin{enumerate}
\item The representatives $A_1,\ldots, A_n$ are isomorphic to parabolic subgroups of $W$.
\item Let $\Gamma$ be a Coxeter graph such that $W\cong W_\Gamma$. Then there exist defining graphs $\Gamma_i$ of $A_i$ for $i=1,\ldots, n$, such that $\Gamma=\Gamma_1\cup\ldots\cup \Gamma_n$. 
\item Let $\Omega$ be a Coxeter graph such that $W\cong W_\Omega$. Then $$\rank(W_\Omega)\leq\pseudorank(A_1)+\ldots+\pseudorank(A_n).$$
\end{enumerate}
\end{corollary}

\begin{remark}
It was shown in \cite{CorsonHughesMollerVarghese2023} that, given two right-angled Coxeter groups $W_\Gamma$ and $W_\Omega$, if $\CF(W_\Gamma)=\CF(W_\Omega)$, then $W_\Gamma\cong W_\Omega$. One can ask if this holds for all Coxeter groups. Unfortunately, this is not the case. Consider for example the Coxeter groups defined by the graphs in \Cref{fig:Delta}.  In this case $\CF(W_\Gamma)=\CF(W_\Omega)$, but $W_\Gamma$ is not isomorphic to $W_\Omega$ since $W_\Gamma$ has Serre's property $\FA$ but $W_\Omega$ does not satisfy $\FA$ (see \Cref{CoxetergroupFA}).  

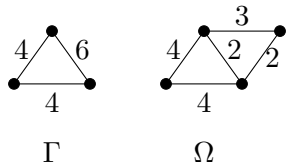
\begin{figure}[htb]
	\begin{center}
	\captionsetup{justification=centering}
		\begin{tikzpicture}
			\draw[fill=black]  (0,0) circle (2pt);
            \draw[fill=black]  (1,0) circle (2pt);
            \draw[fill=black]  (0.5,0.7) circle (2pt);
            \draw (0,0)--(1,0);
            \draw (1,0)--(0.5, 0.7);
            \draw (0,0)--(0.5, 0.7);
            \node at (0.5, -0.25) {$4$};
            \node at (0.1, 0.45) {$4$};
            \node at (0.9, 0.45) {$6$};

            \draw[fill=black]  (2,0) circle (2pt);
            \draw[fill=black]  (3,0) circle (2pt);
            \draw[fill=black]  (2.5,0.7) circle (2pt);
            \draw[fill=black]  (3.5,0.7) circle (2pt);
            \draw (2,0)--(3,0);
            \draw (3,0)--(2.5, 0.7);
            \draw (2,0)--(2.5, 0.7);
             \draw (3.5,0.7)--(3,0);
            \draw (3.5,0.7)--(2.5, 0.7);
            \node at (2.5, -0.25) {$4$};
            \node at (2.1, 0.45) {$4$};
            \node at (2.9, 0.45) {$2$};
            \node at (3, 0.9) {$3$};
            \node at (3.4, 0.35) {$2$};

            \node at (0.5, -0.9) {$\Gamma$};
            \node at (2.5, -0.9) {$\Omega$};
            
        \end{tikzpicture}
        \caption{Two Coxeter diagrams giving non-isomorphic Coxeter groups with isomorphic posets of finite subgroups up to conjugacy.}
	      \label{fig:Delta}
    \end{center}
\end{figure}
\end{remark}

\subsection{Profinite completions and invariants}
Let $G$ denote a group and $\mathcal{N}$ the set of all finite index normal subgroups of $G$. We equip every quotient $G/N$, $N\in \mathcal{N}$ with the discrete topology and endow the product $\prod_{N\in \mathcal{N}}G/N$ with the product topology. Moreover, we define a homomorphism $i\colon G\to\prod_{N\in \mathcal{N}}G/N$ by $g\mapsto (gN)_{N\in \mathcal{N}}$. This homomorphism is injective if and only if $G$ is residually finite. The \textit{profinite completion} of $G$, denoted by $\widehat{G}$, is defined as $\widehat{G}\colon=\overline{i(G)}$.

The next theorem shows that the set $\mathcal{F}(G)$ of isomorphism classes of finite quotients of a finitely generated residually finite group $G$ encodes the same information as $\widehat{G}$.

\begin{thm} \emph{\cite{DixonFormanekPolandRibes1982}}
\label{Dixon}
Let $G$ and $H$ be finitely generated residually finite groups. Then $\mathcal{F}(G)=\mathcal{F}(H)$ if and only if $\widehat{G}\cong\widehat{H}$.
\end{thm}

Note that by the work of Nikolov--Segal \cite{NikolovSegal2007I,NikolovSegal2007II} we have that $\widehat{G}$ is isomorphic to $\widehat{H}$ as a topological group if and only if $\widehat{G}$ is isomorphic to $\widehat{H}$ as an abstract group.

\begin{defn}
Let $\mathcal{C}$ be a class of finitely generated residually finite groups. A property $\mathcal{P}$ is a $\mathcal{C}$-profinite invariant, if for $G,H\in\mathcal{C}$, whenever $G$ has property $\mathcal{P}$ and $\widehat{G}\cong\widehat{H}$ , then $H$ has $\mathcal{P}$.
\end{defn}

\subsection{Separable subgroups}
In this section we record some examples of separable subgroups in Coxeter groups that will be used later in the paper.

By definition, a subgroup $H\subseteq G$  is \emph{separable} if for every $g\notin H$ there exists a homomorphism $\varphi\colon G\to F$ to a finite group such that $\varphi(g)\notin\varphi(H)$. We note that $G$ is residually finite if and only if the trivial subgroup is separable.
It is easy to see that a subgroup $H\subseteq G$ is closed in the profinite topology on $G$ if and only if it is separable.

\begin{lemma}~\cite[Theorem 1.4]{CooperLongReid1998}
    \label{parabolic is separable}
    Let $W_\Omega$ be a special parabolic subgroup of a Coxeter group $W_\Gamma$.  Then $W_\Omega$ is separable. 
\end{lemma}

Recall that a subgroup $H\subseteq G$ is called a \emph{virtual retract},  
if there exists a subgroup $K\subseteq G$ such that $[G:K]<\infty$, $H\subseteq K$ and there is a homomorphism $\varphi\colon K\to H$ which restricts to the identity map on $H$. It is known that virtual retracts have the full profinite topology. Since every special parabolic subgroup in a right-angled Coxeter group is a canonical retract, we obtain the following result.

\begin{lemma}
\label{FullProfTopParaSubgpRACG}
Let $W_\Gamma$ be a right-angled Coxeter group and $W_\Omega$ be a special parabolic subgroup. Then $\overline{W_\Omega}\cong\widehat{W_\Omega}$.
\end{lemma}

\begin{defn}
Let $\mathtt{\Gamma}$ be a Coxeter-Dynkin diagram. The associated Coxeter group $W_\mathtt{\Gamma}$ is called \emph{aﬃne-free} if $\mathtt{\Gamma}$ has no irreducible affine subdiagram of rank $n\geq 3$.
\end{defn}
Examples of affine-free Coxeter groups are Gromov-hyperbolic Coxeter groups and Coxeter groups of FC type.

\begin{lemma}\label{virtuallyVRisenough}
    Let $G$ be a residually finite group.  Suppose $H\leqslant G$ is virtually a virtual retract.  If $H\leqslant K\leqslant G$ with $[K:H]<\infty$ and $K$ separable in $G$, then $K$ has the full profinite topology.
\end{lemma}
\begin{proof}
   By passing to a finite index subgroup of $H$ we may suppose that $H$ is characteristic in $K$.  Now, since $H$ is a virtual retract we have that $H$ has the full profinite topology in $G$.  Let $L$ be a finite index subgroup of $K$.  By \cite[Discussion before Lemma~4.6]{Reid2015}, we need to show that there exists a finite index subgroup $G_L\leqslant G$ such that $G_L\cap L\leqslant L$.  To do this let $H_L=H\cap L$ and note that $H_L$ has finite index in $H$ and so is closed in the profinite topology on $G$.  Thus, there exists $G_L\leqslant G$ such that $G_L\cap H \leqslant H_L$ as required.
\end{proof}

\begin{lemma}\label{lem:affine-free_parabs}
    If $W_\Gamma$ is affine-free and $W_\Omega$ is a special parabolic subgroup, then $W_\Omega$ is virtually a virtual retract. Moreover, $\overline W_\Omega \cong \widehat W_\Omega$
\end{lemma}
\begin{proof}
Let $W_\Gamma$ be an affine-free Coxeter group and $W_\Omega$ be a special parabolic subgroup. By \cite{NibloReeves2003}, $W_\Gamma$ acts cocompactly on the associated Niblo--Reeves $\CAT(0)$ cube complex. Hence, $W_\Gamma$ is a virtually compact special group. Furthermore, $W_\Omega$ is a quasiconvex subgroup of $W_\Gamma$, see \cite{NibloReeves2003}. It was proven in \cite{HaglundWise2008} that any quasiconvex subgroup $H$ in a virtually compact special group $G$ is virtually a virtual retract. Since virtually virtual retracts have the full profinite topology by \Cref{virtuallyVRisenough}, we know that $W_\Omega$ has the full profinite topology. 
\end{proof}

It is not known if all special parabolic subgroups in an arbitrary Coxeter group are virtual retracts, but if the special parabolic subgroup is virtually abelian, then it follows from the following lemma that this group is a virtual retract.

\begin{lemma}\cite[Corollary 1.6, Proposition 1.5]{Minasyan2021}
\label{virtullyAbelianVirtualRetracts}\label{lemma:AmenableSubgroupsFullProfiniteTopology}
Let $W_\Gamma$ be a Coxeter group and $H$ be a subgroup. If $H$ is a finitely generated virtually abelian subgroup, then $H$ is a virtual retract of $W_\Gamma$. Moreover, $\overline{H}\cong\widehat{H}$. 
\end{lemma}

\subsection{Profinite pushouts}
For profinite groups $\mathbf G_1$ and $\mathbf G_2$ with a common closed subgroup $\mathbf H$, we denote the pushout  $\mathbf G_1$ and $\mathbf G_2$ over $\mathbf H$ by $\mathbf G=\mathbf G_1\coprod_{\mathbf H} \mathbf G_2$.  If the natural maps from $\mathbf G_1$ and $\mathbf G_2$ to $\mathbf G$ are embeddings then we call $\mathbf G$ the \emph{profinite amalgamated product} of $\mathbf G_1$ and $\mathbf G_2$ along $\mathbf H$.

Following the proof of \cite[Lemma~3.3]{CorsonHughesMollerVarghese2023} we obtain the following result.
\begin{lemma}
\label{AmalgamProfCompletion}
Let $G\cong A\ast_C B$ be a finitely generated residually finite group. If
$A$, $B$ and $C$ are virtual retracts of $G$,
then $\widehat{G}\cong\widehat{A}\coprod_{\widehat{C}}\widehat{B}$.
\end{lemma}

In general if $G\cong A\ast_C B$ is a finitely generated residually finite group, then $\widehat G\cong \overline A\coprod_{\overline C}\overline B$ is a profinite amalgamated product.  The profinite amalgamated product is \emph{non-trivial} if $\overline A\neq \overline C \neq \overline B$.  In this case $\widehat G$ acts on an infinite profinite tree $T_{\overline A\coprod_{\overline C}\overline B}$ called the \emph{profinite Bass--Serre tree} of the splitting. We refer the reader to \cite{Ribes2017} for more information.

\subsection{Gromov-hyperbolic Coxeter groups}
We need the following characterisation of Gromov-hyperbolicity for Coxeter groups provided by Moussong.

\begin{thm}\label{NoZsquare}\emph{\cite[Theorem B]{Moussong88}}
    Let $W_\Gamma$ denote a Coxeter group. Then the following are equivalent:
    \begin{enumerate}
        \item $W_\Gamma$ is Gromov-hyperbolic.
        \item $W_\Gamma$  contains no subgroup isomorphic to $\Z\times \Z$.
        \item $W_\Gamma$ does not have irreducible affine special parabolic subgroups of rank $\geq 3$ and $W_\Gamma$ does not have infinite special parabolic subgroups $W_\Omega$, $W_\Delta$ such that $W_\Omega$ commutes with $W_\Delta$. 
    \end{enumerate}
\end{thm}

We finish this section by proving that Gromov-hyperbolicity is a profinite invariant amongst Coxeter groups.

\begin{thm}
\label{profiniteinvariantHyperbolicity}
Let $W,W'$ be Coxeter groups such that $\widehat{W}\cong \widehat{W'}$.  Then, $W$ is Gromov-hyperbolic if and only if $W'$ is Gromov-hyperbolic.  In this case $\widehat{W}$ does not contain a $\widehat {\Z}^2$ subgroup.
\end{thm}
\begin{proof}
If $W$ is not Gromov-hyperbolic, then $W$ has a subgroup isomorphic to $\Z\times\Z$ by \Cref{NoZsquare}. By \Cref{lemma:AmenableSubgroupsFullProfiniteTopology} we obtain that $\overline{\Z\times\Z}\cong \widehat{\Z}\times\widehat{\Z}$.

By \cite[Corollary~1.5]{CapraceMuhlherr2005} (building on \cite{NibloReeves2003}) and \cite[Theorem 8.1]{HaglundWise2010} Gromov-hyperbolic Coxeter groups are virtually compact special. Therefore, by \cite[Theorem~D]{WiltonZalesskii2017} the profinite completion of a Gromov-hyperbolic Coxeter group does not contain a subgroup isomorphic to $\widehat{\Z}\times\widehat{\Z}$.
\end{proof}

\section{Goodness and applications}\label{sec good}
\subsection{Coxeter groups are good} 

In the context of profinite completions, \textit{goodness} in the sense of Serre, also known as cohomological goodness, is a cohomological property with many implications. Hence, we recall its definition here, see \cite{Serre1997} for more information.
\begin{defn}
    A finitely generated residually finite group $G$ is called \textit{good} if for every finite $G$-module $M$ and every $q\geq 0$, the induced map on cohomology $H^q(\widehat{G}; M)\to H^q(G;M)$ is an isomorphism.  If the map is an isomorphism for $q\leq n$, then we say $G$ is \emph{$n$-good}.
\end{defn}

In this section we prove that all Coxeter groups are \emph{good} in the sense of Serre.
To begin we recount a construction and a result of Genevois \cite{Genevois2024};
although the reader can safely skip the next definition, but we include it for the more detail-oriented reader.

\begin{defn}
Let $\Gamma=(V(\Gamma),E(\Gamma))$ be a simplicial graph, let $m\colon E(\Gamma)\to \NN_{\geq 2}$ be a labelling of its edges, and let $\calg=\{G_v\mid v\in V(\Gamma)\}$ be a collection of groups. We refer to the groups in $\calg$ as the \emph{vertex-groups}. Assume that, for every edge $\left\{u, v\right\}\in E(\Gamma)$, if $m(\{u, v\}) > 2$ then $G_u$, $G_v$ have order two. The \emph{periagroup} $\Pi(\Gamma, \calg, m)$ is the group with relative presentation
\[\left\langle G_v,\ v\in V(\Gamma)\ |\ \langle G_u,G_v\rangle^{m(\{u,v\})}=\langle G_v,G_u\rangle^{m(\{u,v\})}, \{u,v\}\in E(\Gamma) \right\rangle,\]
where $\langle a,b\rangle^n$ refers to the word obtained from $ababab\dots$ by keeping only the first $n$ letters; and $\langle G_u, G_v\rangle^n=\langle G_v,G_u\rangle^n$ is a shorthand for $\langle a,b\rangle^n=\langle b,a \rangle^n$ for all non-trivial $a\in G_u$ and $b\in G_v$.
\end{defn}

\begin{remark}
   If every group $G\in\calg$ is isomorphic to the cyclic group of order $2$ and the graph $\Gamma$ is finite, then the periagroup $\Pi(\Gamma,\calg,m)$ is exactly the Coxeter group with graph $\Gamma$ and edge-labelling $m$.
\end{remark}

\begin{thm}~\cite[Theorem~1.2]{Genevois2024}
Let $\Pi(\Gamma,\calg,m)$ be a periagroup with $\Gamma$ finite. Then, there exists a finite graph $\Phi$ and a collection $\calh$ indexed by $V(\Phi)$ of groups isomorphic to groups from $\calg$ such that the periagroup $\Pi(\Gamma,\calg,m)$  virtually embeds into the graph product $\Phi \calh$ as a virtual retract.
\end{thm}

\begin{corollary}\label{cor:Cox_vir_retract_RACG}
   Let $W$ be a Coxeter group.  Then, $W$ is virtually a virtual retract of a right-angled Coxeter group.
\end{corollary}

\begin{prop}\label{CoxGpGood}
    If $W$ is a Coxeter group, then $W$ is good.
\end{prop}
\begin{proof}
    We first suppose that $W_\Gamma$ is right-angled.  We will show $W_\Gamma$ is good by induction on the number of vertices $v$ in the defining graph for $\Gamma$.  Note that the argument in this case is essentially identical to the argument for right-angled Artin groups in \cite[Proof of Proposition~3.8]{MinasyanZalesskii2016}.  
    
    The base case, $v=0$ is trivial, since then $W_\Gamma=\{1\}$.  Now, suppose $v>0$ and suppose $\Gamma$ is not complete, as otherwise $W_\Gamma$ is finite and thus good. Therefore we find a vertex $x \in V(\Gamma)$ such that $\st(x)\neq V(\Gamma)$. By \Cref{amalgam.decompostion.Coxeter} we have 
    $$W_\Gamma\cong W_{\st(x)} \ast_{W_{\lk(x)}} W_{V(\Gamma)-\{x\}}.$$ 
    Note that $\st(x)$, $\lk(x)$ and $V(\Gamma)-\{x\}$ have less vertices than $V(\Gamma)$.
    
    All three subgroups 
    $W_{\st(x)}, W_{\lk(x)}$ and $W_{V(\Gamma)-\{x\}}$ 
    have the full profinite topology induced on them 
    by \Cref{FullProfTopParaSubgpRACG} and are good by induction hypothesis.  Hence, $W_\Gamma$ is good by \cite[Proposition~3.6]{GrunewaldJaikinZalesskii2008} (see also \cite[Corollary~3.11]{Lorensen2008}).

    Now, suppose $W_\Gamma$ is an arbitrary Coxeter group.  By \Cref{cor:Cox_vir_retract_RACG} we see that $W_\Gamma$ has a finite index subgroup $H$ that is a virtual retract of some right-angled Coxeter group $W_\Lambda$.  Hence, by \cite[Lemma~3.1]{MinasyanZalesskii2016} $H$ is good.  But goodness is a commensurability invariant \cite[Lemma~3.2]{GrunewaldJaikinZalesskii2008} so we conclude $W_\Gamma$ is good as well.
\end{proof}

\begin{remark}
    Note that essentially the same argument implies that a graph product or a periagroup, where the defining graph is finite and the vertex groups are good, is good. In particular, Dyer groups are good.
\end{remark}

\subsection{Some easy applications}
Recall that via the Tits representation Coxeter groups are linear over $\CC$.  Since Coxeter groups are finitely generated, Selberg's Theorem applies \cite{Selberg1960} and we see that Coxeter groups are virtually torsion-free.
The next lemma shows that the minimal index of a normal torsion-free subgroup in a Coxeter group is a profinite invariant amongst Coxeter groups. This result will be useful in proving that some classes of Coxeter groups are almost $\calw$-profinite rigid.

\begin{lemma}\label{bound.on.edgelabels.profinite}
Let $W_\Gamma$ and $W_\Omega$ be Coxeter groups. Let $H\subseteq W_\Gamma$ be a normal torsion-free subgroup of minimal index $d=[W_\Gamma:H]$. If $\widehat{W_\Gamma}\cong\widehat{W_\Omega}$, then $W_\Omega$ has a normal torsion-free subgroup $L$ of minimal index $d$.
In particular, $m(E(\Omega))\subseteq\left\{2,\ldots, d\right\}$.
\end{lemma}
\begin{proof}
Since $W_\Gamma$ is good by \Cref{CoxGpGood}, the profinite completion of a torsion-free subgroup of finite index in $W_\Gamma$ is torsion-free, see \cite[Lemma 3.3]{MinasyanZalesskii2016}. Further, $W_\Gamma$ and $\widehat{W_\Gamma}$ have the same finite quotients, in particular $W_\Gamma$ and $\widehat{W_\Gamma}$ have the same finite quotients where the corresponding normal subgroups of finite index are torsion-free. By assumption, $\widehat{W_\Gamma}\cong\widehat{W_\Omega}$, so $W_\Omega$ has a normal torsion-free subgroup $L$ of minimal index $d$.  

Let $e$ be an edge in $\Omega$ with edge-label $n$. Then the associated special parabolic subgroup $W_\Delta$ is of type $\mathtt{I}_2(n)$. The canonical map $\psi\colon W_\Omega\twoheadrightarrow W_\Omega/L$ restricted to $W_\Delta$ is injective, since $L$ is torsion-free. Thus, we obtain
\[n\leq|W_\Delta|\leq |W_\Omega/L|=d.\qedhere\]
\end{proof}

As a corollary we obtain:
\begin{corollary}\label{cor:almostRigidity}
   Let $W_\Gamma$ and $W_\Omega$ be Coxeter groups. Assume that $\widehat{W_\Gamma}\cong\widehat{W_\Omega}$. If there is a bound on $|V(\Omega)|$ in terms of the combinatorics of $\Gamma$, then $W_\Gamma$ is almost $\mathcal{W}$-profinitely rigid.  
\end{corollary}
\begin{proof}
By \Cref{bound.on.edgelabels.profinite} there exists $d\in\mathbb{N}$ such that each edge-label in $\Gamma$ and $\Omega$ is at most $d$.
Since there are only finitely many Coxeter graphs with at most a given number of vertices and a uniform bound on the edge labels we obtain almost $\calw$-profinite rigidity of $W_\Gamma$. 
\end{proof}

A group $G$ is of \emph{type $\mathsf{F}$} if there exists a finite model for a $\mathsf{K}(G,1)$, that is a finite connected CW complex $X$ with $\pi_1(X)\cong G$ and $\pi_n(X)=0$ for $n\geq 2$.  We say $G$ is of \emph{type} $\mathsf{VF}$ if $G$ has a finite index subgroup of type $\mathsf{F}$.

Let $G$ be a group of type $\mathsf{VF}$ and let $H$ be a finite index subgroup of type $\mathsf{F}$.  Following Brown \cite[IX.7]{Brown1982}, the \emph{Euler characteristic} of $G$ is defined to be
\[\chi(G)\coloneqq \frac{\chi(H)}{|G:H|},\]
where $\chi(H)$ is the Euler characteristic of any finite $\mathsf K(H,1)$ complex.  Note that this definition is independent of the choice of $H$.

We record two important consequences of goodness for us.  The argument to deduce profinite invariance of these properties is completely standard, but we include it for completeness.

\begin{corollary}\label{profInv_X_H*}
    Let $W_\Gamma$ and $W_\Lambda$ be Coxeter groups (or more generally good groups of type $\mathsf{VF}$).  If $\widehat{W_\Gamma}\cong\widehat {W_\Lambda}$, then 
    \begin{enumerate}
        \item $H^\ast(W_\Gamma;\FF_p)\cong H^\ast(W_\Lambda;\FF_p)$;
        \item $\chi(W_\Gamma)=\chi(W_\Lambda)$.
    \end{enumerate}
\end{corollary}
\begin{proof}
    We first prove (1).  Since Coxeter groups are good (\Cref{CoxGpGood}) we have, for every $n$, natural isomorphisms
    \[H^n(W_\Gamma;\FF_p) \xleftarrow{} H^n(\widehat{W_\Gamma};\FF_p) \to H^n(\widehat{W_\Lambda};\FF_p) \xrightarrow{} H^n( W_\Lambda;\FF_p).\]
     Moreover, by naturality of the cup product these give rise to ring isomorphisms
    \[H^\ast(W_\Gamma;\FF_p) \xleftarrow{} H^\ast(\widehat{W_\Gamma};\FF_p) \to H^\ast(\widehat{W_\Lambda};\FF_p) \xrightarrow{} H^\ast(W_\Lambda;\FF_p).\]

    We now prove (2).  Let $H_\Gamma$ be a finite index torsion-free subgroup of $W_\Gamma$ and let $H_\Omega$ denote the corresponding finite index subgroup of $W_\Omega$ under the isomorphism $\widehat {W_\Gamma}\cong \widehat{W_\Omega}$.  Since Coxeter groups are good, we see that $H_\Gamma$ and $H_\Omega$ are good.  Hence, $\widehat {H_\Gamma}\cong \widehat{H_\Lambda}$ is torsion-free by \cite[Lemma~3.3]{MinasyanZalesskii2016}.  In particular, $H_\Omega$ is torsion-free.  Now, by the argument in (1), we have $H^\ast(H_\Gamma;\FF_p)\cong H^\ast(H_\Omega;\FF_p)$, and so $\chi(H_\Gamma)=\chi(H_\Omega)$.  Since $|W_\Gamma:H_\Gamma|=|W_\Lambda:H_\Lambda|$, the result follows.
\end{proof}

Note that the Euler characteristic is not a profinite invariant in general, see  for example \cite{KammeyerKionkeRaimbaultSauer2020}.

\subsection{Schur multipliers} \label{sec good Schur}
In this section we show that the Schur multiplier is a profinite invariant of Coxeter groups.  We will later use this result in \Cref{Section:OddCoxeterGroups}.

Let $G$ be a group and $\CC^\times$ be the multiplicative group of complex numbers with trivial $G$-action. The \emph{Schur multiplier},
$M(G)$ is defined to be the cohomology group $H^2(G;\CC^\times)$.  Note that when $M(G)$ is finite, we have 
\[H^2(G;\CC^\times)\cong \hom(H_2(G;\Z),\CC^\times) \cong H_2(G;\Z).\]

\begin{remark}
A result of Howlett \cite[Theorem~A]{Howlett1988}, states that for a Coxeter group $W$, the Schur multiplier $M(W)$ is a finite abelian $2$-group.  Thus, when discussing the Schur multiplier $M(W)$ we can and will use the isomorphism $M(W)=H_2(W;\Z)$.
\end{remark}

\begin{lemma}\label{lem M(W) compute}
    Let $W_\Gamma$ be a Coxeter group.  Then, the Schur multiplier of $W_\Gamma$ is isomorphic to the direct sum of
    \[\dim_{\FF_2} H^2(W_\Gamma;\FF_2) - 2\dim_{\FF_2}W_\Gamma^{\mathrm{ab}} \]
    copies of $\FF_2$, noting that $W_\Gamma^{\mathrm{ab}}$ is isomorphic to $\FF_2^n$ for some $n$.
\end{lemma}
\begin{proof}
    By the Universal Coefficient Theorem for cohomology we have
    \begin{equation}\label{eqn SchurMult UCT coho}
        H^2(W_\Gamma;\FF_2) \cong H_2(W_\Gamma;\FF_2) \oplus \Ext^1_\Z(W_\Gamma^{\mathrm{ab}};\FF_2)
    \end{equation}
    and by the Universal Coefficient Theorem for homology we have
    \begin{equation}\label{eqn SchurMult UCT ho}
        H_2(W_\Gamma;\FF_2) \cong (H_2(W_\Gamma;\ZZ)\otimes_\Z \FF_2) \oplus \Tor_1^\Z (W_\Gamma^{\mathrm{ab}},\FF_2).
    \end{equation}
    Now, $W_\Gamma^{\mathrm{ab}}$ is isomorphic to $\FF_2^n$ for some $n$, so it follows that $\Ext^1_\Z(W_\Gamma^{\mathrm{ab}},\FF_2) \cong 
 W_\Gamma^{\mathrm{ab}}$ and that $\Tor_1^\Z (W_\Gamma^{\mathrm{ab}},\FF_2)\cong W_\Gamma^{\mathrm{ab}}$.  Substituting these isomorphisms and \eqref{eqn SchurMult UCT ho} into \eqref{eqn SchurMult UCT coho}, we obtain
 \[H^2(W_\Gamma;\FF_2) \cong (H_2(W_\Gamma;\Z)\otimes_\Z \FF_2) \oplus W_\Gamma^{\mathrm{ab}} \oplus W_\Gamma^{\mathrm{ab}}. \]
 But, $H_2(W_\Gamma;\Z)$ is isomorphic to $\FF_2^m$ for some $m$ by \cite[Theorem~A]{Howlett1988}.  In particular, $H_2(W_\Gamma;\Z)\otimes_\Z \FF_2 \cong H_2(W_\Gamma;\FF_2)$. Making this substitution we obtain
 \[H_2(W_\Gamma;\FF_2) \cong H_2(W_\Gamma;\ZZ) \oplus W_\Gamma^{\mathrm{ab}} \oplus W_\Gamma^{\mathrm{ab}}.\]
 Whence, the lemma.
\end{proof}

\begin{prop}\label{M(W) prof}
    Let $W_\Gamma$ and $W_\Lambda$ be Coxeter groups.  If $\widehat {W_\Gamma}\cong \widehat{W_\Lambda}$, then $M(W_\Gamma)\cong M(W_\Lambda)$.
\end{prop}
\begin{proof}
    By \Cref{lem M(W) compute}, it suffices to show that quantities \[\dim_{\FF_2} H^2(W_\Gamma;\FF_2) - 2\dim_{\FF_2}W_\Gamma^{\mathrm{ab}} \quad \text{and} \quad \dim_{\FF_2} H^2(W_\Lambda;\FF_2) - 2\dim_{\FF_2}W_\Lambda^{\mathrm{ab}}\]
    are equal.  That the first term in each expression is equal follows from \Cref{CoxGpGood}, since we have isomorphisms
    \[ 
    H^2(W_\Gamma;\FF_2) \xleftarrow{} H^2(\widehat{W_\Gamma};\FF_2) \to H^2(\widehat{W_\Lambda};\FF_2) \xrightarrow{} H^2(W_\Lambda;\FF_2).  \]
    That the second term in each expression above is equal follows from the fact that the abelianisation is a profinite invariant amongst all finitely generated residually finite groups.
\end{proof}

\subsection{Conjugacy classes of finite subgroups} \label{sec good CF}
Let $G$ be a group and $H, L\leqslant G$ be two non-conjugate subgroups. By definition $H$ is \emph{conjugacy separable} from $L$ if there exists a homomorphism $\varphi\colon G\to F$, with $F$ finite, such that $\varphi(H)$ is not conjugate to $\varphi(L)$.  

The normaliser of a subgroup $H\subseteq G$ is denoted by $N_G(H)$ and is defined as $N_G(H):=\left\{g\in G\mid gHg^{-1}=H\right\}$. The centraliser of a subgroup $H\subseteq G$ is denoted by 
$C_G(H):=\left\{g\in G\mid ghg^{-1}=h\text{ for all }h\in H\right\}$.

\begin{prop}\label{prop:CF}
    Let $W$ be a Coxeter group, let $\iota\colon W\to\widehat W$ denote the canonical map, and let $p$ be a prime.  The following conclusions hold:
    \begin{enumerate}
        \item $\iota$ induces an order isomorphism $\CFp(W)\to \CFp(\widehat{W})$.\label{CF:p}
        \item $p$-torsion elements of $W$ are conjugacy distinguished.
        \item For every finite $p$-subgroup $H$ of $W$ we have $C_{\widehat W}(H)=\overline{C_W(H)}$ and $N_{\widehat W}(H)=\overline{N_W(H)}$.\label{CF:p_norms}
    \end{enumerate}
    Suppose additionally that $W$ is virtually compact special and virtually toral relatively hyperbolic.  
    The following conclusions holds:
    \begin{enumerate}[resume]
        \item $\iota$ induces an order isomorphism $\CFsol(W)\to \CFsol(\widehat{W})$ and an order monomorphism $\CF(W)\hookrightarrow\CF(\widehat{W})$.\label{CF:sol}
        \item For every finite subgroup $H$ of $G$ we have $N_{\widehat W}(H)=\widehat{N_W(H)}$.  For every finitely generated subgroup $K$ of $G$ we have $C_{\widehat W}(K)=\overline{C_W(K)}$.\label{CF:sol_norm}
    \end{enumerate}
\end{prop}
\begin{proof}
Items (4) and (5) are an application of a result of Boggi and Zalesskii \cite[Theorem~6.4]{BoggiZalesskii2024}.  Items (1), (2), and (3) follow from a result of Boggi and Zalesskii \cite[Theorem~A]{BoggiZalesskii2024} (note that it is a generalisation of \cite[Corollary 3.5]{MinasyanZalesskii2016}).  We briefly explain how to verify the hypothesis of Boggi and Zalesskii's theorem:
\begin{enumerate}[label=(\alph*)]
    \item \emph{Goodness.}  This follows from \Cref{CoxGpGood}.
    
    \item \emph{Finite virtual $p$-cohomological type.}  We have $\vcd(W)<\infty$ so in particular, $\vcd_p(W)<\infty$.  Moreover, $W$ is type $\mathsf{F}_\infty$ so if $M$ is a finite $\FF_pW$-module, then $H^n(W;M)$ is finite for all $n\geq0$.
    
    \item \emph{The natural map $H^n(\widehat W;M)\to H^n(W;M)$ is an isomorphism for every discrete $\FF_p\llbracket \widehat W\rrbracket$-module $M$.}  This essentially follows from \Cref{CoxGpGood}, which states that $W$ is good in the sense of Serre.  Note that goodness implies that the natural map $H^n(\widehat W;A)\to H^n(W;A)$ is an isomorphism for every finite $\widehat W$-module.  To reconcile this with the desired conclusion, note that a discrete $\FF_p\llbracket \widehat W\rrbracket$-module $M$ is a direct limit of its finite submodules, and that the cohomology of profinite groups with coefficients in discrete modules commutes with direct limits (c.f \cite[page 11]{Serre1997}). \qedhere
\end{enumerate}
\end{proof}

\section{Product decompositions}\label{sec product decomps}
A Coxeter group $W_\Gamma$ can be decomposed in the form 
$W_\Gamma\cong W_{\Gamma_{{\rm sph}}}\times W_{\Gamma_{{\rm aff}}}\times W_{\Gamma_{\rm {gen}}}$, 
where $W_{\Gamma_{{\rm sph}}}$ is trivial or a finite Coxeter group called the \emph{spherical part} of $W_\Gamma$, $W_{\Gamma_{{\rm aff}}}$ is trivial or a product of irreducible affine Coxeter groups called the \emph{affine part} of $W_\Gamma$, and $W_{\Gamma_{\rm {gen}}}$ is trivial or an infinite non-affine Coxeter group called the \emph{generic part} of $W_\Gamma$. More precisely, if $W_{\Gamma_{\rm {gen}}}$ is non-trivial, then every irreducible direct factor in $W_{\Gamma_{\rm {gen}}}$ is an infinite non-affine Coxeter group, see \cite{ParisVarghese2024}.

The profinite completion of a group splits as a direct product if the base group splits, but an isomorphism of such a direct product does not necessarily respect this decomposition. However, the natural splitting of a Coxeter group $W_\Gamma\cong W_{\Gamma_{{\rm sph}}}\times W_{\Gamma_{{\rm aff}}}\times W_{\Gamma_{\rm {gen}}}$  is respected by profinite isomorphisms. More precisely:

\setcounter{thmx}{1}
\begin{thmx}\label{thmx:products}
Let $W_\Gamma=W_{\Gamma_{{\rm sph}}}\times W_{\Gamma_{{\rm aff}}}\times W_{\Gamma_{{\rm gen}}}$ and $W_\Omega=W_{\Omega_{{\rm sph}}}\times W_{\Omega_{{\rm aff}}}\times W_{\Omega_{{\rm gen}}}$ be  Coxeter groups such that $ \widehat{W_\Gamma}\cong \widehat{W_\Omega}$. Then \[W_{\Gamma_{{\rm sph}}}\cong W_{\Omega_{{\rm sph}}},\quad 
W_{\Gamma_{{\rm aff}}}\cong W_{\Omega_{{\rm aff}}},\quad \text{and} \quad \widehat{W_{\Gamma_{{\rm gen}}}}\cong \widehat{W_{\Omega_{{\rm gen}}}}.\]
In particular, a  Coxeter group $W_\Gamma$ is profinitely rigid amongst  Coxeter groups if and only if $W_{\Gamma_{\rm gen}}$ is profinitely rigid amongst  Coxeter groups.
\end{thmx}

We first prove two simple lemmata.

\begin{lemma}\label{NoFiniteNormalSugps}
    Suppose $W_\Gamma$ is a  Coxeter group.  Let $W_\Gamma\cong \{1\}\times W_{\Gamma_{{\rm aff}}}\times W_{\Gamma_{{\rm gen}}}$ and $N\trianglelefteq \widehat{W_\Gamma}$ be a normal subgroup. If $N$ is finite, then $N\cong\{1\}$.
\end{lemma}
\begin{proof}
    Suppose $N$ is not trivial. Then there exists a prime $p$ and an element $n\in N$ of order $p$. By \Cref{prop:CF} we have $\mathcal{CF}_p(W_\Gamma)\cong \mathcal{CF}_p(\widehat{W_\Gamma})$. Therefore we find an $m\in W_\Gamma$ and a $w\in \widehat{W_\Gamma}$ such that $n=wmw^{-1}$, since $N$ is normal, $m\in N\cap W_\Gamma$. Therefore $\{1\}\neq N\cap W_\Gamma \trianglelefteq W_\Gamma$ is a non-trivial finite normal subgroup, a contradiction to \cite[Theorem~1.1]{ParisVarghese2024}, since direct products of affine and generic Coxeter groups do not admit non-trivial finite normal subgroups.
\end{proof}

\begin{lemma}\label{NoVirtAbelianNormalSUbgps}
    Suppose $W_\Gamma$ is a  Coxeter group.    Let $W_\Gamma\cong \{1\}\times \{1\}\times W_{\Gamma_{{\rm gen}}}$ and $N\trianglelefteq \widehat{W_{\Gamma}}$ be a non-trivial normal subgroup. If $N$ is virtually abelian, then $N$ is torsion-free. 
\end{lemma}
\begin{proof}
Let $N\trianglelefteq \widehat{W_{\Gamma}}$ be a non-trivial virtually abelian normal subgroup.
    Suppose $N$ contains a torsion element. Then there exists a prime $p$ and an element $n\in N$ of order $p$. By \Cref{prop:CF} we have $\mathcal{CF}_p(W_\Gamma)\cong \mathcal{CF}_p(\widehat{W_\Gamma})$. Therefore we find an $m\in W_\Gamma$ and a $w\in \widehat{W_\Gamma}$ such that $n=wmw^{-1}$, since $N$ is normal, $m\in N\cap W_\Gamma$. Therefore $\{1\}\neq N\cap W_\Gamma \trianglelefteq W_\Gamma$ is a non-trivial virtually abelian normal subgroup, a contradiction to \cite[Theorem~1.1]{ParisVarghese2024}, since generic type Coxeter groups do not admit non-trivial virtually abelian normal subgroups.
\end{proof}

Now we are able to prove \Cref{thmx:products}.
\begin{proof}[Proof of \Cref{thmx:products}]
 Let
 $$f\colon \widehat{W_{\Gamma_{{\rm sph}}}}\times\widehat{W_{\Gamma_{{\rm aff}}}}\times \widehat{W_{\Gamma_{{\rm gen}}}}\to\widehat{W_{\Omega_{{\rm sph}}}}\times\widehat{W_{\Omega_{{\rm aff}}}}\times\widehat{W_{\Omega_{{\rm gen}}}}$$

denote an isomorphism. Let $\pi^{\Gamma}_i$ denote the projection to the $i$-th coordinate in $\widehat{W_\Gamma}$ and $\pi^{\Omega}_i$ in $\widehat{W_\Omega}$.

    Since $\widehat{W_{\Gamma_{{\rm sph}}}}$ is normal in $\widehat{W_\Gamma}$, $\pi^{\Omega}_j\left(f\left(\widehat{W_{\Gamma_{{\rm sph}}}}\right)\right)$ is a finite normal subgroup in the $j$-th component of $\widehat{W_\Omega}$. By \Cref{NoFiniteNormalSugps}, the projections to the second and third coordinate are trivial. Therefore we conclude $f\left(\widehat{W_{\Gamma_{{\rm sph}}}}\right)\subseteq \widehat{W_{\Omega_{{\rm sph}}}}$. Repeating this argument with $f^{-1}$ yields $W_{\Gamma_{{\rm sph}}}\cong W_{\Omega_{{\rm sph}}}$.

Hence $f$ induces an isomorphism $$g\colon \widehat{W_{\Gamma_{{\rm aff}}}}\times \widehat{W_{\Gamma_{{\rm gen}}}}\to \widehat{W_{\Omega_{{\rm aff}}}}\times \widehat{W_{\Omega_{{\rm gen}}}}.$$

    Now we repeat the same argument with $g$, and as normal subgroup we use the virtually abelian normal subgroup $\widehat{W_{\Gamma_{{\rm aff}}}}$. Note that this contains torsion as Coxeter groups are generated by involutions. This time we invoke \Cref{NoVirtAbelianNormalSUbgps} to obtain that $g\left(\widehat{W_{\Gamma_{{\rm aff}}}}\right)\subseteq \widehat{W_{\Omega_{{\rm aff}}}}$. Similarly $g^{-1}\left(\widehat{W_{\Omega_{{\rm aff}}}}\right)\subseteq \widehat{W_{\Gamma_{{\rm aff}}}}$. Since $g\circ g^{-1} \left(\widehat{W_{\Omega_{{\rm aff}}}}\right) = \widehat{W_{\Omega_{{\rm aff}}}}$, the restricted map 
    $$g_{\mid \widehat{W_{\Gamma_{{\rm aff}}}}}\colon \widehat{W_{\Gamma_{{\rm aff}}}}\to \widehat{W_{\Omega_{{\rm aff}}}}$$ is surjective and since it is also injective we have an  isomorphism between $\widehat{W_{\Gamma_{{\rm aff}}}}$ and $\widehat{W_{\Omega_{{\rm aff}}}}$.  Thus, by \cite[Theorem~1.1]{CorsonHughesMollerVarghese2024} we obtain $W_{\Gamma_{{\rm aff}}}\cong W_{\Omega_{{\rm aff}}}$.

    For the generic part consider 
    $$\pi_{\rm {gen}}^{\Omega}\circ g\colon\widehat{W_{\Gamma_{{\rm aff}}}}\times \widehat{W_{\Gamma_{{\rm gen}}}}\to \widehat{W_{\Omega_{{\rm aff}}}}\times \widehat{W_{\Omega_{{\rm gen}}}}\twoheadrightarrow \widehat{W_{\Omega_{{\rm gen}}}} .$$ 
    This has kernel $\widehat{W_{\Gamma_{{\rm aff}}}}$, therefore we conclude via the homomorphism theorem that $\widehat{W_{\Gamma_{{\rm gen}}}}\cong \widehat{W_{\Omega_{{\rm gen}}}}$ as desired.
\end{proof}

\begin{corollary}\label{virtually.abelian.Coxetergroups.profinitely.rigid}
Virtually abelian Coxeter groups are $\mathcal{W}$-profinitely rigid.
\end{corollary}
\begin{proof}
Let $W_\Gamma$ be a virtually abelian Coxeter group and $W_\Omega$ be a Coxeter group such that $\widehat{W_\Gamma}\cong\widehat{W_\Omega}$.  Note that $W_\Omega$ is also virtually abelian because $\widehat W_\Gamma$ is.  Therefore $W_\Omega$ is good. By a characterisation of virtually abelian Coxeter groups, see \cite[Theorem~2.1]{MollerVarghese2023} we have a decomposition $W_\Gamma\cong W_{\Gamma_{\rm sph}}\times W_{\Gamma_{{\rm aff}}}$. Hence, by \Cref{thmx:products} follows that $W_\Omega\cong W_{\Omega_{\rm sph}}\times W_{\Omega_{{\rm aff}}}$ and $W_{\Gamma_{\rm sph}}\cong W_{\Omega_{\rm sph}}$ and $\widehat{W_{\Gamma_{\rm aff}}}\cong \widehat{W_{\Omega_{\rm aff}}}$. Since products of irreducible affine Coxeter groups are profinitely rigid by \cite{CorsonHughesMollerVarghese2024}, it follows that $W_{\Gamma_{\rm aff}}\cong W_{\Omega_{\rm aff}}$ and therefore $W_\Gamma\cong W_\Omega$.
\end{proof}

\begin{corollary}
\label{centreless}
Let $W_\Gamma$ and $W_\Omega$ be Coxeter groups. If $\widehat{W_\Gamma}\cong\widehat{W_\Omega}$, then $Z(W_\Gamma)\cong Z(W_\Omega)$. In particular, to be centreless is a $\mathcal{W}$-profinite invariant.
\end{corollary}
\begin{proof}
By \cite[Section 6.3]{Humphreys1990} the centre of a Coxeter group $W$ is contained in the spherical part of $W$, $Z(W)\subseteq Z(W_{{\rm sph}})$ and is trivial or $Z(W)\cong\Z_2 ^n$ for $n\geq 1$. Using \Cref{thmx:products} we conclude that $Z(W_\Gamma)\cong Z(W_\Omega)$. Hence, to be centreless is a $\mathcal{W}$-profinite invariant.
\end{proof}

\section{Coxeter groups of low cohomological dimension} \label{sec low vcd}
In this section we will show that the properties of being virtually free and virtually surface are detected by the profinite completion of a Coxeter group. 

\subsection{Virtually free Coxeter groups} \label{sec low vcd vir free}
We give a cohomological argument that being virtually free is a profinite invariant amongst Coxeter groups.  The proof is inspired by work of Wilton \cite{Wilton2021}. Note that by \cite[Theorem~34]{MihalikTschantz2009}, a Coxeter group $W_\Gamma$ is virtually free if and only if $\Gamma$ is chordal and $W_\Gamma$ is of FC type. However, we cannot detect chordality of the defining graph profinitely, hence, we need a different characterisation for our purposes.

\begin{thm}\label{vir free W profinite}
    Let $W_\Gamma$ be a Coxeter group.  The following are equivalent:
    \begin{enumerate}
        \item $\vcd W_\Gamma \leq 1$;
        \item $W_\Gamma$ is virtually free;
        \item $W_\Gamma$ does not contain a quasiconvex surface subgroup;
        \item $\widehat{W_\Gamma}$ does not contain the profinite completion of an infinite surface group as a subgroup.
    \end{enumerate}
\end{thm}
\begin{proof}
    The equivalence of (1) and (2) is essentially Stallings' Theorem.  The equivalence of (2) and (3) is a theorem of Gordon, Long, and Reid \cite[Theorem~1.1]{GordonLongReid2004} combined with the observation that their proof constructs a quasiconvex surface subgroup.  We now explain how the first three properties imply (4).  Since $W_\Gamma$ is virtually free we may pass to a finite index free subgroup $H$.  Now, the profinite completion $\widehat H$ of $H$ has $\cd_\mathbf{p}\widehat H=1$ for every prime $p$.  But, the profinite completion of a surface group $\widehat{\pi_1\Sigma_g}$ for $g\geq 1$ has $\cd_\mathbf{p}\widehat{\pi_1\Sigma_g}=2$ for every prime $p$.  Since $p$-cohomological dimension of profinite groups is monotonic with respect to subgroups, we see that $\widehat{\pi_1\Sigma_g}$ is never a subgroup of $\widehat H$.  
    
    Finally, we show (4) implies (3), more precisely we prove the contrapositive.  To this end, suppose $W_\Gamma$ contains a surface subgroup.  We aim to show that $W_\Gamma$ contains another surface subgroup as a virtual retract.  In this case one of the following holds:
    \begin{enumerate}[label=\alph*)]
        \item every proper parabolic subgroup of $W_\Gamma$ is finite;
        \item $\Gamma$ contains an induced $n$-cycle for some $n\geq 4$.
    \end{enumerate}
    
    In case (a) we see that $W_\Gamma$ is either an affine Coxeter group or a hyperbolic triangle group, or $\Gamma$ corresponds to a rank $4$ or $5$ Lann\'er diagram.  If $W_\Gamma$ is an affine Coxeter group, then $W_\Gamma$ is virtually abelian of dimension at least $2$ and so contains a $\Z^2$ as a virtual retract.  If $W_\Gamma$ is a hyperbolic triangle group, then $W_\Gamma$ is virtually a hyperbolic surface group.  If $\Gamma$ corresponds to a rank $4$ or $5$ Lann\'er diagram, then by the argument in \cite[Proof of Theorem~2.3]{GordonLongReid2004}, taking a $2$-cell $F$ of a simplex, the centraliser of the reflection in the hyperbolic plane spanned by $F$ contains a hyperbolic surface group $H$ of finite index.  Such a subgroup $H$ is necessarily quasiconvex and so since hyperbolic Coxeter groups are virtually compact special by Haglund--Wise \cite{HaglundWise2010}, $H$ is a virtual retract by \cite[Theorem~2.7]{LongReid2008}.

    In case (b) $W_\Gamma$ has a proper parabolic subgroup $W_\Omega$ which is virtually an infinite surface group.  Moreover $W_\Gamma$ is hyperbolic and hence affine-free.  It follows that, $W_\Omega$ is a virtual retract by \Cref{lem:affine-free_parabs}.  Restricting this retract to a torsion-free subgroup of $W_\Gamma$, gives a virtual retract onto a hyperbolic surface subgroup of $W_\Gamma$ as required.
\end{proof}

\begin{corollary}\label{profinite.invariant. virtuallyfreeness}
    Let $W_\Gamma$ and $W_\Lambda$ be Coxeter groups with $\widehat {W_\Gamma}\cong \widehat{W_\Lambda}$.  Then, $W_\Gamma$ is a virtually free group if and only if $W_\Lambda$ is a virtually free group.
\end{corollary}

\begin{corollary}\label{vir free almost rigid}
    Let $W_\Gamma$ be a virtually free Coxeter group.  Then $W_\Gamma$ is almost profinitely rigid amongst Coxeter groups.
\end{corollary}
\begin{proof}
    By \Cref{profinite.invariant. virtuallyfreeness}, any other Coxeter group with isomorphic profinite completion is virtually free.  The result now follows from \cite[Theorem~3.3]{GrunewaldZalesskii2011}.
\end{proof}

\begin{prop}\label{amalgam.finitegroups.Outabelian}
Let $W_\Gamma$ be a Coxeter group. Assume that there exist finite special parabolic subgroups $W_{\Gamma_1}$, $W_{\Gamma_2}$, $W_{\Gamma_3}$ such that $W_\Gamma\cong W_{\Gamma_1}*_{W_{\Gamma_3}} W_{\Gamma_2}$. If the outer automorphism group $\Out(W_{\Gamma_3})$ is abelian, then $W_\Gamma$ is $\mathcal{W}$-profinitely rigid.

In particular, if $W_{\Gamma_3}$ is of type $\mathtt{A}_{n(n\geq 1)}$, $\mathtt{B}_{n(n\geq 2)}$, $\mathtt{D}_{n(n\geq 5)}$, $\mathtt{E}_{6}$, $\mathtt{E}_{7}$, $\mathtt{E}_{8}$, $\mathtt{H}_3$, $\mathtt{H}_4$ or $\mathtt{I}_2(m\geq 3)$, then $W_\Gamma$ is $\mathcal{W}$-profinitely rigid.
\end{prop}
\begin{proof}
Let $W_\Omega$ be a Coxeter group such that $\widehat{W_\Gamma}\cong\widehat{W_\Omega}$. The Coxeter group $W_\Gamma$ is by assumption virtually free. Since being virtually free is a $\calw$-profinite invariant by \Cref{profinite.invariant. virtuallyfreeness}, it follows that $W_\Omega$ is virtually free. Now we apply
\cite[Theorem~4.1(3)]{GrunewaldZalesskii2011} to conclude that $W_\Gamma\cong W_\Omega$.

If $W_{\Gamma_3}$ is of type $\mathtt{A}_{n(n\geq 1)}$, $\mathtt{B}_{n(n\geq 2)}$, $\mathtt{D}_{n(n\geq 5)}$, $\mathtt{E}_{6}$, $\mathtt{E}_{7}$, $\mathtt{E}_{8}$, $\mathtt{H}_3$, $\mathtt{H}_4$ or $\mathtt{I}_2(m\geq 3)$, then $\Out(W_{\Gamma_3})$ is abelian, see \cite[Table 1]{Bannai1969}. 
Hence $W_\Gamma$ is $\mathcal{W}$-profinitely rigid.
\end{proof}

\subsection{Virtually surface Coxeter groups} \label{sec low vcd vir surface}

The following argument is inspired by work of Fruchter and Morales \cite{FruchterMorales2022} on residually free groups.

We define the $n$th virtual Betti number of a group $G$ to be
\[\vb_n(G)\coloneqq \sup\{\dim_\QQ H_n(K;\QQ): K\leqslant G,\ |G:K|<\infty\}\in\NN\cup\{\infty\}. \]
It is easy to see that virtual Betti numbers are profinite invariants of good groups.

\begin{lemma}\label{lem vb infty}
    Let $G$ be virtually compact special hyperbolic group.  If $H\leqslant G$ is an infinite index quasiconvex subgroup with $\vb_n(H)>0$, then $\vb_n(G)=\infty$.
\end{lemma}
\begin{proof}
    Abusing notation we replace $H$ by a finite index subgroup satisfying $\mathsf b_n(H)>0$, such a subgroup is still quasiconvex. Now, a theorem of Arzhantseva \cite[Theorem~1]{Arzhantseva2001} implies that there exists a quasiconvex subgroup $K=\Z\ast H\leqslant G$.  By elementary Bass--Serre theory, we see that $K$ has finite index subgroups of the form $K_\ell\coloneqq F_{k_\ell}\ast H_1\ast\dots\ast H_\ell$ for arbitrarily large $\ell$ and some $k_\ell\in \NN$.  A simple Mayer--Vietoris sequence argument shows that $\mathsf{b}_\ell(K_\ell)\geq n$.  Hence, $\vb_n(K)=\infty$.  But, $K$ is a virtual retract of $G$, hence so is every $K_\ell$.  Thus, we have a sequence of finite index subgroups $G_\ell\leqslant G$ such that $K_\ell$ is retract of $G_\ell$.  Since retracts induce an inclusion of cohomology groups we see that $\mathsf{b}_n(G_\ell)\geq \ell$.  Hence, $\vb_n(G)=\infty$ as required.
\end{proof}

By a \emph{surface group} we mean the fundamental group of a closed orientable surface of genus at least $1$. By combining \cite[Theorem~A.2]{DaniThomas2017} with  \cite{GordonLongReid2004} we obtain a visual characterisation of virtually surface Coxeter groups.
\begin{lemma}\label{virtuallySurfaceGroups.Graphcharacterization}
Let $W_\Gamma$ be a Coxeter group. Then $W_\Gamma$ is virtually a surface group if and only if there exist special parabolic subgroups $W_\Omega$ and $W_\Delta$ such that $W_\Omega$ is finite, $W_\Delta$ is infinite and $\Delta$ is an $n$-cycle for some $n\geq 3$ and $W_\Gamma\cong W_\Omega\times W_\Delta$.
\end{lemma}

\begin{prop}\label{surface vb1}
    Let $W_\Gamma$ be a Coxeter group.  The following are equivalent:
    \begin{enumerate}
        \item $W_\Gamma$ is virtually a non-abelian surface group;
        \item $\widehat{W_\Gamma}$ is virtually the profinite completion of a non-abelian surface group;
        \item $W_\Gamma$ is Gromov-hyperbolic and $\vb_2(W_\Gamma)=1$.
    \end{enumerate}
\end{prop}
\begin{proof}
    It is easy to see (1) implies (2).  That (2) implies (3) is proved as follows: Hyperbolicity is detected by the profinite completion of a Coxeter group \Cref{profiniteinvariantHyperbolicity}.   Now, virtual Betti numbers are profinite invariants of good groups, so \Cref{CoxGpGood} completes the claim. 
    
    We now show (3) implies (1).  More precisely we show a strong contrapositive, namely, if $W_\Gamma$ is hyperbolic, satisfies $\vcd W_\Gamma\geq 2$, and is not virtually a non-abelian surface group, then $\vb_2(W_\Gamma)=\infty$.  Since $W_\Gamma$ has cohomological dimension at least $2$, by \Cref{vir free W profinite} we see that $W_\Gamma$ contains a quasiconvex surface group $H$.  Now, $H$ is necessarily infinite index because $W_\Gamma$ is not virtually surface group, and $H$ is necessarily non-abelian because $W_\Gamma$ is Gromov-hyperbolic.  The result follows from \Cref{lem vb infty} since Gromov-hyperbolic Coxeter groups are virtually compact special by \cite[Theorem~8.1]{HaglundWise2010} and \cite[Corollary 1.5]{CapraceMuhlherr2005}.
\end{proof}

\begin{corollary}\label{vir surface profinite}
    Let $W_\Gamma$ and $W_\Lambda$ be Coxeter groups with $\widehat {W_\Gamma}\cong \widehat{W_\Lambda}$.  Then, $W_\Gamma$ is a virtually surface group if and only if $W_\Lambda$ is a virtually surface group.
\end{corollary}
\begin{proof}
    The virtually abelian case is covered by \cite{CorsonHughesMollerVarghese2024}.  Now, by \Cref{profiniteinvariantHyperbolicity}, $W_\Gamma$ is Gromov-hyperbolic if and only if $W_\Lambda$ is.  Hence, if one of them is a virtually non-abelian surface group then the other is Gromov-hyperbolic.  
    
    Since virtual Betti numbers are profinite invariants of Coxeter groups, we see that if one of the groups is virtually surface, then the other group has $2$nd virtual Betti number equal to one. Thus, the conditions of \Cref{surface vb1}(3) are satisfied and both groups are virtually non-abelian surface groups.
\end{proof}

\begin{lemma}\label{Virtually.surface.groups.CF}
    Let $W$ be a Gromov-hyperbolic virtually surface Coxeter group with no non-trivial finite normal subgroup.  Then, the natural map $W\to\widehat{W}$ induces an isomorphism of posets $\CF(W)\to\CF(\widehat W)$.
\end{lemma}
\begin{proof}
    Every finite subgroup of $W$ is soluble.  Since $W$ is Gromov-hyperbolic and virtually compact special we may apply \Cref{prop:CF}\eqref{CF:sol} to obtain a poset isomorphism $\CF(W)\to\CFsol(\widehat W)$.  Now, an argument exactly as in \cite[Theorem~5.1]{BridsonConderReid2016} shows that $\CFsol(\widehat W)=\CF(\widehat W)$.
\end{proof}

\begin{thm}\label{triangleAndmore}
    Let $W_\Gamma$ be a hyperbolic virtually surface Coxeter group.  The following conclusions hold:
    \begin{enumerate}
        \item $|\calg_\calw(W_\Gamma)|$ is finite;
        \item if $\Gamma$ is a triangle, then $|\calg_\calw(W_\Gamma)|=1$; 
        \item if $\Gamma$ is an $n$-gon such that at most one edge label in $\Gamma$ differs from the others, then $|\calg_\calw(W_\Gamma)|=1$;
        \item if $m(E(\Gamma))\subseteq 2\NN$, then $|\calg_\calw(W_\Gamma)|=1$.
    \end{enumerate}
\end{thm}
Note that the case where $m(E(\Gamma))=\{n\}$ corresponds to a reflection group in the sides of a regular hyperbolic polygon.
\begin{proof}
Before we prove the numbered statements we establish some notation and facts.  Let $W_\Gamma$ be a Gromov-hyperbolic virtually surface group. It follows from \Cref{virtuallySurfaceGroups.Graphcharacterization} that there exist special parabolic subgroups $W_\Omega$ and $W_\Delta$ such that $W_\Omega$ is finite, $W_\Delta$ is infinite and Gromov-hyperbolic and $\Delta$ is an $n$-gon with edge labels $e_1,\dots, e_n$ for some $n\geq 3$ and $W_\Gamma\cong W_\Omega\times W_\Delta$. 

The poset $\mathcal{CF}(W_\Delta)$ has a maximal element for each edge in the $n$-gon since maximal finite subgroups are never conjugate, see \Cref{FactsCoxeter}. Thus, we can tell the number of edges from $\CF(W_\Delta)$. Moreover, we know what the edge labels are. In particular, we have a bijection between $\CFmax(W_\Delta)$ and the multi-set $\left\{e_1,\ldots, e_n\right\}$. Note that a multi-set is a collection in which the order of the elements do not matter but the multiplicities do.

Let $W_{\Gamma'}$ be a Coxeter group such that $\widehat{W_\Gamma}\cong\widehat{W_{\Gamma'}}$. Since being a Gromov-hyperbolic virtually surface group is a $\mathcal{W}$-profinite invariant (\Cref{vir surface profinite}, \Cref{profiniteinvariantHyperbolicity}), it follows that $W_{\Gamma'}$ is a Gromov-hyperbolic virtually surface group. Hence there exist special parabolic subgroups $W_{\Omega'}$ and $W_{\Delta'}$ such that $W_{\Omega'}$ is finite, $W_{\Delta'}$ is infinite and Gromov-hyperbolic and $\Delta'$ is an $m$-gon with edge labels $f_1,\dots, f_m$ for some $m\geq 3$ and $W_{\Gamma'}\cong W_{\Omega'}\times W_{\Delta'}$. 

Let $$\Phi\colon W_\Omega\times \widehat{W_\Delta}\to W_{\Omega'}\times\widehat{W_{\Delta'}}$$
be an isomorphism. Then by \Cref{thmx:products} we have
$W_\Omega\cong W_{\Omega'}$ and $\widehat{W_\Delta}\cong\widehat{W_{\Delta'}}$. By applying \Cref{Virtually.surface.groups.CF} we obtain
$$\CF(W_\Delta)=\CF(\widehat{W_\Delta})=\CF(\widehat{W_{\Delta'}})=\CF(W_{\Delta'}).$$
In particular, we have a bijection between the multi-sets $\left\{e_1,\ldots, e_n\right\}\to\left\{f_1,\ldots, f_m\right\}$ and therefore $n=m$.  With this in hand we now prove the numbered statements in the theorem.

We now prove (1). Since there are finitely many $n$-gons with multi-set $\left\{e_1,\ldots, e_n\right\}$ the set $|\calg_\calw(W_\Gamma)|$ is finite.

We now prove (2). If $\Gamma$ is a triangle, then $n=m=3$ and therefore $\Omega$ is also a triangle. Two triangles with the same multi-sets of edge-labels are isomorphic. Hence $W_\Gamma\cong W_\Omega$.

We now prove (3). If at most one edge label in $\Gamma$ differs from the others, then $\Gamma\cong\Omega$ and therefore $W_\Gamma\cong W_\Omega$.

Finally, we prove (4). If all edge labels in $\Gamma$ are even, then the edges with edge-labels $e_i$ and $e_j$ intersect in a vertex if and only if the corresponding conjugacy classes in $\CF(W_\Gamma)$, $[D_{e_i}]$ and $[D_{e_j}]$ have a non-trivial common lower bound. Since $\CF(W_\Delta)=\CF(W_{\Delta'})$ we obtain $\Delta\cong\Delta'$.
\end{proof}

\bigskip
\section{Actions on trees}\label{sec trees}
In this section, we prove that Serre's property $\mathsf{FA}$ is a profinite invariant amongst Coxeter groups. Moreover, we show that the number of ends is a profinite invariant amongst Coxeter groups, that Coxeter groups lie in class $\mathcal{A}$ of \cite{AndersonBessaZalesskii2023}, and prove \Cref{thmx:Reduction1ended}.

For background on groups acting on profinite trees see \cite{Ribes2017}.

\subsection{Property \texorpdfstring{$\FA$}{FA} and invariance of FC type} \label{sec trees FA FC}
  
\begin{defn}
    A group $G$ is said to have Serre's property $\FA$ if every action of $G$ on a simplicial tree $T$ without edge inversion has a fixed point.
\end{defn}

Following \cite{Caprace2006}, a group $G$ is called $2$-\emph{spherical} if it possesses a finite generating set $S$ such that any pair of elements of $S$ generates a finite subgroup. By \cite{Serre2003}, a $2$-spherical group has property $\FA$.

\begin{lemma}\label{CoxetergroupFA}
    A Coxeter group $W_\Gamma$ has property $\FA$ if and only if $\Gamma$ is complete.
\end{lemma}
\begin{proof}
    This is standard.  It follows immediately from the fact that a complete Coxeter group is $2$-spherical coupled with the observation that a non-complete $\Gamma$ always gives rise to a non-trivial amalgamated product decomposition, see \Cref{amalgam.decompostion.Coxeter}.
\end{proof}

\begin{prop}~\cite[Proposition~2.4.9]{Ribes2017}\label{HellyThm}
    Let $T$ be a profinite tree and let $T_1,\ldots,T_n$ be profinite subtrees for some $n\geq 1$. If $T_i\cap T_j\neq \emptyset$ for all $i,j\in \{1,\ldots,n\}$, then $\bigcap_{i=1}^{n}T_i\neq \emptyset$.
\end{prop}

\begin{prop}\label{2-spherical}
Let $G$ be a residually finite $2$-spherical group. Let $\mathbf G_1\coprod_\mathbf{H} \mathbf G_2$ be a profinite amalgamated product.
If $\widehat{G}$ is a subgroup of $\mathbf G_1\coprod_\mathbf{H} \mathbf G_2$, then $\widehat{G}$ is a subgroup of a conjugate of $\mathbf G_1$ or $\mathbf G_2$.
\end{prop}
\begin{proof}
The result is trivial if the profinite amalgamated product is trivial.  So suppose the profinite amalgamated product is non-trivial.  The non-trivial profinite amalgamated product $\mathbf G_1\coprod_\mathbf{H} \mathbf G_2$ induces an action of $\widehat{G}$ on the infinite profinite tree $T_{\mathbf G_1\coprod_\mathbf{H} \mathbf G_2}$ associated to this profinite amalgam. Since $G$ is residually finite, it embeds in its profinite completion and thus acts on the profinite tree $T_{\mathbf G_1\coprod_\mathbf{H} \mathbf G_2}$.

Let $S$ be a $2$-spherical generating set of $G$. Each generator $s_i\in S$ generates a finite subgroup, hence each $s_i$ has a fixed point by \cite[Theorem~3.10]{MelnikovZalesskii1989}. Moreover, by \cite[Theorem~2.8]{MelnikovZalesskii1989}, for every $s_i$ we have that ${\rm Fix}(s_i)$ is a profinite subtree. Now, since $S$ is $2$-spherical, every pair of generators $s_i, s_j$ generates a finite group, which also has a fixed point. Therefore, we conclude ${\rm Fix}(s_i)\cap {\rm Fix}(s_j)\neq \emptyset$. By \Cref{HellyThm}, we have $\bigcap_{i=1}^{n}{\rm Fix}(s_i)\neq \emptyset$. Hence, $G$ has a global fixed point. Since the action on the profinite tree is continuous and $G$ is dense in $\widehat{G}$, we conclude that $\widehat{G}$ has a global fixed point as well, hence $\widehat{G}$ is a subgroup of a conjugate of $\mathbf G_1$ or $\mathbf G_2$.
\end{proof}

In general, Serre's property $\FA$ is not a profinite invariant \cite{CheethamWestLubotzkyReidSpitler2022}. But for Coxeter groups it turns out to be the case.

\begin{thm}
\label{profiniteinvariantFA}
Serre's fixed point property $\FA$ is a $\mathcal{W}$-profinite invariant. 
\end{thm}
\begin{proof}
Let $W_\Gamma$ be a complete Coxeter group and $W_\Omega$ be a Coxeter group such that $\widehat{W_\Gamma}\cong\widehat{W_\Omega}$. 

Assume that $\Omega$ is not complete, then there exist vertices $v, w\in V(\Omega)$, $v\neq w$ such that $\left\{v,w\right\}\notin E(\Omega)$. Thus $W_\Omega=A*_C B$ where $A=W_{V(\Omega)-\left\{v\right\}}$, $B=W_{st(v)}$ and $C=W_{lk(v)}$. By \Cref{parabolic is separable}, special parabolic subgroups are closed in the profinite topology of $W_\Omega$.  Hence 
$\widehat{W_\Omega}\cong \overline{A}\coprod_{\overline{C}} \overline{B}$ is a profinite amalgamated product with $\overline C\neq \overline A,\overline B\neq \widehat W_\Omega$.  This splitting induces an action of $\widehat{W_\Omega}$ on the infinite profinite tree $T_{\overline{A}\coprod_{\overline{C}} \overline{B}}$ associated to this profinite amalgam. Applying \Cref{2-spherical} we conclude that $\widehat{W_\Gamma}$ has a global fixed point, contradicting the construction of the tree $T_{\overline{A}\coprod_{\overline{C}} \overline{B}}$.
\end{proof}

By definition, a Coxeter group $W_\Gamma$ is of \emph{FC type} if every complete special parabolic subgroup is finite.
\begin{thm}\label{profinite.invariant.FC}
Let $W_\Gamma$ and $W_\Omega$ be Coxeter groups with $\widehat{W_\Gamma}\cong\widehat{W_\Omega}$. Then, $W_\Gamma$ is of FC type if and only if $W_\Omega$ is of FC type.
\end{thm}
\begin{proof}
Let $W_\Gamma$ be of FC type and let $W_\Delta$ be a complete parabolic subgroup of $W_\Omega$. We note that Coxeter groups of FC type are affine-free, thus every special parabolic subgroup of $W_\Gamma$ has the full profinite topology by \Cref{lem:affine-free_parabs}.

If $\Gamma$ is complete, then $W_\Gamma$ is finite and therefore $W_\Delta$ is finite. Assume that $\Gamma$ is not complete. Then $\widehat{W_\Gamma}$ is a profinite amalgamated product of profinite completions of special parabolic subgroups $\widehat{W_{\Gamma_1}}\coprod_{\widehat{W_{\Gamma_3}}}\widehat{W_{\Gamma_2}}$. Since $W_\Delta$ is $2$-spherical, we conclude that $\widehat{W_\Delta}$ is contained in a conjugate of $\widehat{W_{\Gamma_1}}$ or $\widehat{W_{\Gamma_2}}$ by \Cref{2-spherical}. Continue splitting $W_{\Gamma_1}$ resp. $W_{\Gamma_2}$ until the groups in the amalgamated product are $2$-spherical. Hence, $\widehat{W_\Delta}$ is contained in a conjugate of $\widehat{W_{\Delta'}}$ where $W_{\Delta'}$ is a complete special parabolic subgroup of $W_\Gamma$. By assumption, $W_\Gamma$ is of FC type, therefore $\widehat{W_{\Delta'}}$ and  $\widehat{W_\Delta}$ are both finite.
\end{proof}

\subsection{Number of ends and the class \texorpdfstring{$\cala$}{A}} \label{sec trees ends class A}

We want to use profinite Bass-Serre theory to show that the number of ends is a profinite invariant for Coxeter groups and prove \Cref{thmx:Reduction1ended}. To do so, we first prove the following lemma.

\begin{lemma}\label{1-endedFixedPoint}
    If $W_\Gamma$ is a $1$-ended Coxeter group, then every action of $\widehat{W_\Gamma}$ on a profinite tree with finite edge stabilisers has a global fixed point.
\end{lemma}
\begin{proof}
    We first observe that the defining graph $\Gamma$ is connected, as else $W_\Gamma$ has infinitely many ends.

    Suppose for a contradiction the lemma is false, so $\widehat{W_\Gamma}$ acts on a profinite tree $T$ with finite edge stabilisers and without global fixed point. We write $S=\left\{s_1,\dots,s_n\right\}$ for a set of Coxeter generators of $W_\Gamma$. Note that by the canonical inclusion, $\iota\colon W_\Gamma\hookrightarrow\widehat{W_\Gamma}$ we have that $W_\Gamma$ acts on $T$ as well.

    First, consider $s_1$, since this element has finite order, by \cite[Theorem~3.10]{MelnikovZalesskii1989}, the subgroup generated by $s_1$ fixes a vertex. Now, consider the subgroup $\langle s_1,s_2\rangle$. Either this fixes a vertex or it does not. If it does, consider the subgroup $\langle s_1,s_2,s_3\rangle$ and check whether it fixes a vertex, else consider the subgroup $\langle s_1,s_3\rangle$. This yields a subset of the generators $S_1$ which fixes a vertex. Note, that by assumption $S_1\neq S$, since we do not have a global fixed point.

    We continue with the generator $s_i$, where $i$ is the lowest number such that $s_i\notin S_1$. Now consider the subgroup $\langle s_1,s_i\rangle$ and check if it fixes a vertex. Continue this process to obtain sets $S_1,\dots,S_k$ with $S_1\cup\dots \cup S_k=S$. Note that this will typically not be a disjoint union.

    Now, $\langle S_1\rangle $ fixes a vertex $v$ and $\langle S_2\rangle$ fixes a vertex $w\neq v$. By \cite[Theorem~3.12]{MelnikovZalesskii1989}, we obtain that $\langle S_1\cap S_2\rangle$ is contained in an edge stabiliser, hence, finite. 

    Furthermore, given $s_j\in S_j$ and $s_j\notin S_\ell$ with $S_j\cap S_\ell\neq \emptyset$, then there is no relation between $s_j$ and $s_\ell$ for every $s_\ell\in S_j-S_\ell$, since otherwise due to construction, the fixed point sets would form a cycle. More precisely, $\Fix(\langle s_j, s_\ell\rangle) \cap \Fix(S_j\cap S_\ell)=\emptyset$, but there is a path connecting $\Fix(S_j)$ and $\Fix(S_\ell)$ through $\Fix(\langle s_j,s_\ell\rangle)$ and through $\Fix(S_\ell\cap S_j)$. 

    Now, we construct a tree $T'$  as follows: vertices are all cosets of $W_{S_i}$ for $i=1,\dots,k$. We draw an edge between two vertices if and only if their intersection is non-empty. Note that, since $\Gamma$ is connected, $T'$ is connected as well. Because $T$ is a profinite tree, we conclude that $T'$ is a tree.
    
    We see that $W$ acts on $T'$ by left-multiplication, without edge inversions, since the action on $T$ had this property. Edge stabilisers are finite, since these are conjugates of finite groups and we do not have a global fixed point by construction.  But now Bass-Serre theory implies that $W$ is infinitely ended, a contradiction.    
\end{proof}

\begin{remark}
    Note that the previous lemma can be easily adapted to detect $\{$virtually $\Z\}$-splittings of one-ended hyperbolic Coxeter groups, although we do not pursue this here.
\end{remark}

Recently, in \cite{AndersonBessaZalesskii2023}, a class of accessible groups has been studied with regard to profinite genus of graphs of groups. We want to point out the following remark.
\begin{remark}\label{ClassA}
    Coxeter groups are in the class $\mathcal{A}$ of \cite{AndersonBessaZalesskii2023}, meaning that they are accessible (this follows by invoking finite presentability and Dunwoody's theorem) and that every vertex group in the JSJ-decomposition has a fixed point for every action on a profinite tree. The latter part is precisely the above lemma.
\end{remark}

Let $W_\Gamma$ be a Coxeter group and $\pi_1(\mathcal{G}, \Delta)$ be the fundamental group of a graph of groups such that $W_\Gamma\cong \pi_1(\mathcal{G}, \Delta)$. We call $\pi_1(\mathcal{G}, \Delta)$
 a \emph{visual graph of groups decomposition}  if each edge group and each vertex group is a special parabolic subgroup of $W_\Gamma$.

In \cite{Dunwoody1985} Dunwoody showed that any finitely presented group has a finite graph of groups decomposition where each edge group is finite and each vertex group is finite or $1$-ended. One might ask how visual this decomposition is for Coxeter groups. In \cite[Corollary 18]{MihalikTschantz2009}, Mihalik and Tschantz proved that every Coxeter group has a visual graph of groups decomposition where each edge group is finite and each vertex group is finite or $1$-ended. Since Coxeter groups are in the class $\mathcal{A}$, we have the following result which is a direct consequence of \cite[Theorem 1.1]{BessaPortoZalesskii2024}.

\begin{thm}\label{profinite.JSJ.Coxeter}
Let $W_\Gamma$ and $W_\Omega$ be Coxeter groups and $\pi_1(\mathcal{G}_1, \Delta_1)$ and $\pi_1(\mathcal{G}_2, \Delta_2)$ their visual graph of groups decompositions where the edge groups are finite and the vertex groups are finite or $1$-ended. If $\widehat{W_\Gamma}\cong\widehat{W_\Omega}$, then there exist bijections $f\colon E(\Delta_1)\to E(\Delta_2)$ and $g\colon V(\Delta_1)\to V(\Delta_2)$ such that $\mathcal{G}_1(e)\cong\mathcal{G}_2(f(e))$ and $\widehat{\mathcal{G}_1(v)}\cong\widehat{\mathcal{G}_2(g(v))}$ for all $e\in E(\Delta_1), v\in V(\Delta_1)$.
\end{thm}

The following example shows how useful \Cref{profinite.JSJ.Coxeter} is for showing that a Coxeter group is $\calw$-profinitely rigid.

\begin{example}
Let $W_\Gamma$ be a Coxeter group where $\Gamma$ is isomorphic to the graph in \Cref{fig:twotriangles} and $m(E(\Gamma))=\left\{3\right\}$.

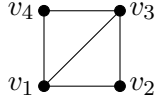
\begin{figure}[htb]
	\begin{center}
	\captionsetup{justification=centering}
		\begin{tikzpicture}

\draw[fill=black]  (3,-2) circle (2pt);
            \draw[fill=black]  (3,-3) circle (2pt);
            \draw[fill=black]  (4,-2) circle (2pt);
            \draw[fill=black]  (4,-3) circle (2pt);
            \draw (3,-3)--(4,-3);
            \draw (3,-3)--(4,-2);
            \draw (3,-3)--(3,-2);
            \draw (3,-2)--(4,-2);
            \draw (4,-3)--(4,-2);
            \node at (2.7, -3) {$v_1$};
            \node at (4.3, -3) {$v_2$};
            \node at (2.7, -2) {$v_4$};
            \node at (4.3, -2) {$v_3$};
            \end{tikzpicture}
        \caption{Two triangles glued together along an edge.}
        \label{fig:twotriangles}
    \end{center}
\end{figure}

Then $W_\Gamma\cong\langle v_1, v_2, v_3\rangle *_{\langle v_1, v_3\rangle}\langle v_1, v_3, v_4\rangle$ is a visual decomposition where the edge group is finite and the vertex groups are $1$-ended. 

Let $W_\Omega$ be a Coxeter group such that $\widehat{W_\Gamma}\cong\widehat{W_\Omega}$. Applying \Cref{profinite.JSJ.Coxeter} we know that there exist special parabolic subgroups $A, B, C\subseteq W_\Omega$ such that $\widehat{W_\Omega}\cong \widehat{A}\coprod_C\widehat{B}$ and $\widehat{A}\cong\widehat{\langle v_1, v_2, v_3\rangle}$, $\widehat{B}\cong \widehat{\langle v_1, v_3, v_4\rangle}$ and $C\cong \langle v_1, v_3\rangle$. Since $\langle v_1, v_2, v_3\rangle, \langle v_1, v_3, v_4\rangle$ are profinitely rigid by \cite{CorsonHughesMollerVarghese2024}, we have $A\cong\langle v_1, v_2, v_3\rangle$, $B\cong\langle v_1, v_3, v_4\rangle$. Further, the Coxeter groups $\langle v_1, v_2, v_3\rangle, \langle v_1, v_3, v_4\rangle$ and $\langle v_1, v_3\rangle$ are additionally graph rigid (\Cref{graph.rigid.Coxetergroups}(3)), so $\Omega$ can be covered by two triangles where each edge label is $3$ and the intersection of the triangles is an edge. Hence there is only one possibility to do that, we conclude that $\Gamma\cong\Omega$ and $W_\Gamma\cong W_\Omega$. 
\end{example}
A direct consequence from \Cref{profinite.JSJ.Coxeter} is:
\begin{corollary}\label{amalgam.profinitelyrigid}
Let $W_\Gamma$ be a Coxeter group. 
Assume $W_\Gamma\cong W_{\Gamma_1}*_{W_{\Gamma_3}} W_{\Gamma_2}$ where
\begin{enumerate}
\item $W_{\Gamma_1}$, $W_{\Gamma_2}$, $W_{\Gamma_3}$ are special parabolic subgroups,
\item $W_{\Gamma_3}$ is finite and graph rigid,
\item $W_{\Gamma_1}$, $W_{\Gamma_2}$ are finite or $1$-ended groups and graph rigid.
\end{enumerate}
Assume additionally that every graph $\Delta$ which can be covered by the full subgraphs $\Gamma_1$ and $\Gamma_2$, i. e. $\Delta=\Gamma_1\cup\Gamma_2$ such that  $\Gamma_1\cap\Gamma_2=\Gamma_3$ is isomorphic to $\Gamma$.

If $W_{\Gamma_1}$ and $W_{\Gamma_2}$ are $\mathcal{W}$-profinitely rigid, then $W_\Gamma$ is $\mathcal{W}$-profinitely rigid.
\end{corollary}

\begin{example}

Let $W_\Gamma$ be a Coxeter group. If $\Gamma$ is isomorphic to the graph in \Cref{fig:TinyHouse} and $m\neq 4k+2$, $k\geq 1$, then $W_\Gamma$ is $\mathcal{W}$-profinitely rigid by \Cref{amalgam.profinitelyrigid}. 

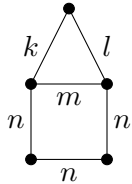
\begin{figure}[htb]
	\begin{center}
	\captionsetup{justification=centering}
		\begin{tikzpicture}

            \draw[fill=black]  (3,-2) circle (2pt);
            \draw[fill=black]  (3,-3) circle (2pt);
            \draw[fill=black]  (4,-2) circle (2pt);
            \draw[fill=black]  (4,-3) circle (2pt);
            \draw[fill=black]  (3.5,-1) circle (2pt);
            \draw (3,-3)--(4,-3);
            \draw (3,-3)--(3,-2);
            \draw (3,-2)--(4,-2);
            \draw (4,-3)--(4,-2);
            \draw (3.5, -1)--(3,-2);
            \draw (3.5, -1)--(4,-2);
            
            \node at (3.5, -3.2) {$n$};
            \node at (3.5, -2.2) {$m$};
            \node at (2.8, -2.5) {$n$};
            \node at (4.2, -2.5) {$n$};

            \node at (3, -1.5) {$k$};
            \node at (4, -1.5) {$l$};
                     
            \end{tikzpicture}
        \caption{Tiny house.}
        \label{fig:TinyHouse}
    \end{center}
\end{figure}
\end{example}

\begin{thm}
\label{profiniteinvariantNumberofends}
    Let $W_\Gamma$ be a Coxeter group. Then the number of ends is a profinite invariant amongst Coxeter groups
\end{thm}
\begin{proof}
    First note that by Hopf's Theorem on ends \cite{Hopf44}, $W_\Gamma$ either has $0,1,2$ or infinitely many ends.

    The case $e(W_\Gamma)=0$ corresponds to $W_\Gamma$ being finite, so this case is trivial.  If $e(W_\Gamma)=2$, then, by \cite[I.8.32]{BridsonHaefliger1999}, $W_\Gamma$ is virtually $\Z$. Hence, $W_\Gamma$ is profinitely rigid by \Cref{virtually.abelian.Coxetergroups.profinitely.rigid}.  It remains to show that we can profinitely distinguish $1$-ended Coxeter groups from infinitely ended Coxeter groups.

    So towards a contradiction, we assume that there exist Coxeter groups $W_\Gamma$, $W_\Omega$ such that $W_\Gamma$ is $1$-ended and $W_\Omega$ has infinitely many ends and $\widehat{W_\Gamma}\cong\widehat{W_\Omega}$. Since $W_\Omega$ has infinitely many ends, by \cite[I.8.32]{BridsonHaefliger1999}, we see that $W_\Omega$ splits as a non-trivial amalgamated free product $A\ast_C B$ where $A,B,C$ are special subgroups and $C$ is finite. Thus we write $A=W_X$, $B=W_Y$ and $C=W_Z$.

     We may write $\widehat{W_\Omega}\cong \overline{W_X}\coprod_{W_Z} \overline{W_Y}$ as a profinite amalgamated product. Note that $\overline{W_Z}=W_Z$ and that $W_X,W_Y,W_Z$ are all separable within $W_\Omega$ by \Cref{parabolic is separable} (but may not have the full profinite topology). Hence, we deduce that $\widehat{W_\Omega}\cong \overline{W_X}\coprod_{W_Z} \overline{W_Y}$ with $W_Z \neq \overline W_X,\overline W_Y\neq \widehat W_\Omega$. This profinite splitting induces an action on an infinite profinite Bass-Serre tree $T$ without global fixed point, without inversion of edges, and with finite edge stabilizers.

     Since $\widehat{W_\Gamma}\cong\widehat{W_\Omega}$, we see that $\widehat{W_\Gamma}$ also acts on $T$ without a global fixed point, without edge inversion, and with finite edge stabilizers. But $W_\Gamma$ is one-ended, so its profinite completion cannot admit such an action by \Cref{1-endedFixedPoint}, a contradiction.
\end{proof}

These arguments allow us to reduce almost $\mathcal{W}$-profinite rigidity to the case of $1$-ended Coxeter groups.
\setcounter{thmx}{3}
\begin{thmx}\label{thmx:Reduction1ended}
    Coxeter groups are almost $\mathcal{W}$-profinitely rigid if and only if $1$-ended Coxeter groups are almost $\mathcal{W}$-profinitely rigid.
\end{thmx}
\begin{proof}
    Suppose that all $1$-ended Coxeter groups are almost $\mathcal{W}$-profinitely rigid. Let $W_\Gamma$ be a Coxeter group.  We show that the $\mathcal{W}$-profinite genus is finite. 
    
 Let $W_\Omega$ be a Coxeter group such that $\widehat{W_\Gamma}\cong \widehat{W_\Omega}$. 
 
 We decompose both groups into graphs of groups with one-ended (or finite) vertex groups and finite edge groups such that all edge and vertex groups are special parabolic subgroups (see \cite[Corollary 18]{MihalikTschantz2009}), 
 $$W_\Gamma\cong \pi_1(\mathcal{G}_1, \Delta_1)\text{ and }W_\Omega \cong \pi_1(\mathcal{G}_2,\Delta_2) .$$ 
 In particular,  $\pi_1(\mathcal{G}_1, \Delta_1)$ and $\pi_1(\mathcal{G}_2,\Delta_2)$ are visual graph of groups decompositions. 
 We note that there exist defining graphs of the vertex groups in $\mathcal{G}_1$ resp. $\mathcal{G}_2$ such that the graph $\Gamma$ resp. $\Omega$ are unions of these graphs. 
 
 By assumption $\widehat{W_\Gamma}\cong \widehat{W_\Omega}$, hence by \Cref{profinite.JSJ.Coxeter} there exist bijections $f\colon E(\Delta_1)\to E(\Delta_2)$ and $g\colon V(\Delta_1)\to V(\Delta_2)$ such that $\mathcal{G}_1(e)\cong\mathcal{G}_2(f(e))$ and $\widehat{\mathcal{G}_1(v)}\cong\widehat{\mathcal{G}_2(g(v))}$ for all $e\in E(\Delta_1), v\in V(\Delta_1)$.

 Let $V(\Delta_1)=\left\{v_1, \ldots, v_k\right\}$ and $\mathcal{G}_1(v_1),\ldots,\mathcal{G}_1(v_k)$ be the vertex groups in $\pi_1(\mathcal{G}_1, \Delta_1)$. By assumption, there exist only finitely many non-isomorphic Coxeter groups $W$ such that $\widehat{W}\cong\widehat{\mathcal{G}_1(v_i)}$. Thus,
 for  $\widehat{\mathcal{G}_1(v_i)}$ we define 
 $$c_i:=\max\left\{\pseudorank(W)\mid \widehat{W}\cong \widehat{\mathcal{G}_1(v_i)}\right\}.$$
 Since $\widehat{\mathcal{G}_1(v)}\cong\widehat{\mathcal{G}_2(g(v))}$ for all $ v\in V(\Delta_1)$ we have 
 $$\rank(W_\Omega)\leq \sum_{i=1}^k c_i.$$ 
By \Cref{cor:almostRigidity}, $W_\Gamma$ is almost $\calw$-profinitely rigid.
\end{proof}

\section{Cocompact simplicial reflection groups}\label{sec lanner}

Let $\mathtt{X}$ be a Coxeter-Dynkin diagram and $W_\mathtt{X}$  be the associated Coxeter group. The Coxeter group $W_\mathtt{X}$ is a  \emph{cocompact simplicial reflection group} if  $W_\mathtt{X}$ acts properly cocompactly on a Riemannian manifold of constant sectional curvature.   The goal of this section is to prove these Coxeter groups are $\calw$-profinitely rigid. 

\begin{thm}\label{thm.Lanner}
    If $W$ is a cocompact Coxeter simplicial reflection group, then $W$ is profinitely rigid amongst Coxeter groups.
\end{thm}
\begin{proof}
    There are three cases: the trivial case when $W$ is finite, the case when $W$ is affine which is covered by \cite{CorsonHughesMollerVarghese2024}, and the hyperbolic case.  In the hyperbolic case $W$ has rank at most $5$ by the classification of Lann\'er \cite{Lanner1950} and Vinberg \cite[Propositions~3.2 and 4.2]{Vinberg1985}.   The rank $2$ case corresponds to $D_\infty$ and so is covered by the affine case.  In rank $3$ we have that $W$ is a hyperbolic triangle group and so is covered by \Cref{{triangleAndmore}}.  In ranks $4$ and $5$ the Coxeter--Dynkin diagram for $W$ is a Lann\'er diagram (c.f. Vinberg above).  That Coxeter groups with Lann\'er diagrams are $\calw$-profinitely rigid is given by \Cref{Lanner rk 4 - Wrigid}.
\end{proof}

\begin{table}[]
\begin{tabular}{|c|c|c|}
\hline
 Label    & $\mathtt{X}$ & $\CFmax(W_\mathtt{X})$  \\ \hline \hline

$\call_1$ & $\vcenter{\hbox{
\begin{tikzpicture}
    \node[circle, fill=black, scale=0.5] (1) at (0,0) {};
    \node[circle, fill=black, scale=0.5] (2) at (1,0) {};
    \node[circle, fill=black, scale=0.5] (3) at (2,0) {};
    \node[circle, fill=black, scale=0.5] (4) at (3,0) {};
    \draw (1) --  (2) -- node[above,midway] {5} (3) -- (4);
\end{tikzpicture}}}$
& $\mathtt{A}_1\times \mathtt{A}_2$,  $\mathtt{A}_1\times \mathtt{A}_2$, $\mathtt{H}_3$, $\mathtt{H}_3$ \\ \hline

$\call_2$ & $\vcenter{\hbox{
\begin{tikzpicture}
    \node[circle, fill=black, scale=0.5] (1) at (0,0) {};
    \node[circle, fill=black, scale=0.5] (2) at (1,0) {};
    \node[circle, fill=black, scale=0.5] (3) at (2,0) {};
    \node[circle, fill=black, scale=0.5] (4) at (3,0) {};
    \draw (1) -- node[above,midway] {5} (2) --  (3) -- node[above,midway] {4} (4);
\end{tikzpicture}}}$
& $\mathtt{A}_1\times \mathtt{B}_2$, $\mathtt{A}_1\times \mathtt{H}_2$, $\mathtt{B}_3$, $\mathtt{H}_3$  \\ \hline

$\call_3$ & $\vcenter{\hbox{
\begin{tikzpicture}
    \node[circle, fill=black, scale=0.5] (1) at (0,0) {};
    \node[circle, fill=black, scale=0.5] (2) at (1,0) {};
    \node[circle, fill=black, scale=0.5] (3) at (2,0) {};
    \node[circle, fill=black, scale=0.5] (4) at (3,0) {};
    \draw (1) -- node[above,midway] {5} (2) --  (3) -- node[above,midway] {5} (4);
\end{tikzpicture}}}$
& $\mathtt{A}_1\times \mathtt{H}_2$, $\mathtt{A}_1\times \mathtt{H}_2$, $\mathtt{H}_3$, $\mathtt{H}_3$  \\ \hline

$\call_4$ & $\vcenter{\hbox{
\begin{tikzpicture}
    \node[circle, fill=black, scale=0.5] (1) at (0,0) {};
    \node[circle, fill=black, scale=0.5] (2) at (1,0) {};
    \node[circle, fill=black, scale=0.5] (3) at (2,0.5) {};
    \node[circle, fill=black, scale=0.5] (4) at (2,-0.5) {};
    \draw (1) -- node[above,midway] {5} (2);
    \draw (3) -- (2) -- (4);
\end{tikzpicture}}}$
& $\mathtt{A}_1^3$, $\mathtt{A}_3$, $\mathtt{H}_3$, $\mathtt{H}_3$  \\ \hline

$\call_5$ & $\vcenter{\hbox{ 
\begin{tikzpicture}
    \node[circle, fill=black, scale=0.5] (1) at (0,1) {};
    \node[circle, fill=black, scale=0.5] (2) at (1,1) {};
    \node[circle, fill=black, scale=0.5] (3) at (1,0) {};
    \node[circle, fill=black, scale=0.5] (4) at (0,0) {};
    \draw (1) -- node[above,midway] {4} (2) -- (3) -- (4) -- (1);
\end{tikzpicture}}}$
& $\mathtt{A}_3$, $\mathtt{A}_3$, $\mathtt{B}_3$, $\mathtt{B}_3$   \\ \hline

$\call_6$ & $\vcenter{\hbox{
\begin{tikzpicture}
    \node[circle, fill=black, scale=0.5] (1) at (0,1) {};
    \node[circle, fill=black, scale=0.5] (2) at (1,1) {};
    \node[circle, fill=black, scale=0.5] (3) at (1,0) {};
    \node[circle, fill=black, scale=0.5] (4) at (0,0) {};
    \draw (1) -- node[above,midway] {4} (2) -- (3) -- node[below,midway] {4} (4) -- (1);
\end{tikzpicture}}}$
& $\mathtt{B}_3$, $\mathtt{B}_3$, $\mathtt{B}_3$, $\mathtt{B}_3$  \\ \hline

$\call_7$ &
$\vcenter{\hbox{ 
\begin{tikzpicture}
    \node[circle, fill=black, scale=0.5] (1) at (0,1) {};
    \node[circle, fill=black, scale=0.5] (2) at (1,1) {};
    \node[circle, fill=black, scale=0.5] (3) at (1,0) {};
    \node[circle, fill=black, scale=0.5] (4) at (0,0) {};
    \draw (1) -- node[above,midway] {5} (2) -- (3) -- (4) -- (1);
\end{tikzpicture}}}$
& $\mathtt{A}_3$, $\mathtt{A}_3$, $\mathtt{H}_3$, $\mathtt{H}_3$  \\ \hline

$\call_8$ & $\vcenter{\hbox{
\begin{tikzpicture}
    \node[circle, fill=black, scale=0.5] (1) at (0,1) {};
    \node[circle, fill=black, scale=0.5] (2) at (1,1) {};
    \node[circle, fill=black, scale=0.5] (3) at (1,0) {};
    \node[circle, fill=black, scale=0.5] (4) at (0,0) {};
    \draw (1) -- node[above,midway] {5} (2) -- (3) -- node[below,midway] {4} (4) -- (1);
\end{tikzpicture}}}$ & $\mathtt{B}_3$, $\mathtt{B}_3$, $\mathtt{H}_3$, $\mathtt{H}_3$  \\ \hline

$\call_9$ & $\vcenter{\hbox{
\begin{tikzpicture}
    \node[circle, fill=black, scale=0.5] (1) at (0,1) {};
    \node[circle, fill=black, scale=0.5] (2) at (1,1) {};
    \node[circle, fill=black, scale=0.5] (3) at (1,0) {};
    \node[circle, fill=black, scale=0.5] (4) at (0,0) {};
    \draw (1) -- node[above,midway] {5} (2) -- (3) -- node[below,midway] {5} (4) -- (1);
\end{tikzpicture}}}$
& $\mathtt{H}_3$, $\mathtt{H}_3$, $\mathtt{H}_3$, $\mathtt{H}_3$\\ \hline              
\end{tabular}
\caption{The rank $4$ Lann\'er diagrams and maximal simplices.}
\label{rank 3 lanner}
\end{table}

\begin{table}[]
\begin{tabular}{|c|c|c|}
\hline
 Label    & $\mathtt{X}$ & $\CFmax(W_\mathtt{X})$  \\ \hline \hline

$\call_1$ & $\vcenter{\hbox{
\begin{tikzpicture}
    \node[circle, fill=black, scale=0.5] (1) at (0,0) {};
    \node[circle, fill=black, scale=0.5] (2) at (1,0) {};
    \node[circle, fill=black, scale=0.5] (3) at (2,0) {};
    \node[circle, fill=black, scale=0.5] (4) at (3,0) {};
    \node[circle, fill=black, scale=0.5] (5) at (4,0) {};
    \draw (1) -- node[above,midway] {5} (2) --  (3) -- (4) -- (5);
\end{tikzpicture}}}$
& $\mathtt{A}_1\times\mathtt{A}_3$,  $\mathtt{A}_1\times\mathtt{H}_3$, $\mathtt{A}_2\times\mathtt{H}_2$, $\mathtt{H}_4$, $\mathtt{A}_4$ \\ \hline

$\call_2$ & $\vcenter{\hbox{
\begin{tikzpicture}
    \node[circle, fill=black, scale=0.5] (1) at (0,0) {};
    \node[circle, fill=black, scale=0.5] (2) at (1,0) {};
    \node[circle, fill=black, scale=0.5] (3) at (2,0) {};
    \node[circle, fill=black, scale=0.5] (4) at (3,0) {};
    \node[circle, fill=black, scale=0.5] (5) at (4,0) {};
    \draw (1) -- node[above,midway] {5} (2) --  (3) -- (4) -- node[above,midway] {4} (5);
\end{tikzpicture}}}$
& $\mathtt A_1\times \mathtt H_3$, $\mathtt A_1\times \mathtt H_3$, $\mathtt H_2^2$, $\mathtt H_4$, $\mathtt H_4$  \\ \hline

$\call_3$ & $\vcenter{\hbox{
\begin{tikzpicture}
    \node[circle, fill=black, scale=0.5] (1) at (0,0) {};
    \node[circle, fill=black, scale=0.5] (2) at (1,0) {};
    \node[circle, fill=black, scale=0.5] (3) at (2,0) {};
    \node[circle, fill=black, scale=0.5] (4) at (3,0) {};
    \node[circle, fill=black, scale=0.5] (5) at (4,0) {};
    \draw (1) -- node[above,midway] {5} (2) --  (3) -- (4) -- node[above,midway] {5} (5);
\end{tikzpicture}}}$
& $\mathtt A_1\times \mathtt B_3$, $\mathtt A_1\times \mathtt H_3$, $\mathtt B_2\times \mathtt H_2$, $\mathtt B_4$, $\mathtt H_4$  \\ \hline

$\call_4$ & $\vcenter{\hbox{
\begin{tikzpicture}
    \node[circle, fill=black, scale=0.5] (1) at (0,0) {};
    \node[circle, fill=black, scale=0.5] (2) at (1,0) {};
    \node[circle, fill=black, scale=0.5] (3) at (2,0) {};
    \node[circle, fill=black, scale=0.5] (4) at (3,0) {};
    \node[circle, fill=black, scale=0.5] (5) at (1,-1) {};
    \draw (2) --  (5);
    \draw (1) -- (2) -- (3) -- node[above,midway] {5} (4);
\end{tikzpicture}}}$
& $\mathtt A_1\times \mathtt A_3$, $\mathtt A_1^2\times \mathtt H_2$, $\mathtt D_4$, $\mathtt H_4$, $\mathtt H_4$  \\ \hline

$\call_5$ & $\vcenter{\hbox{ 
\begin{tikzpicture}
    \foreach \a in {0,1,2,3,4} 
        \node[circle, fill=black, scale=0.5] (\a) at (72*\a-18:1cm) {};
    \draw (0) -- (1) -- node[above,midway] {4} (2) --  (3) -- (4) --  (0);
\end{tikzpicture}}}$
& $\mathtt A_4$, $\mathtt A_4$, $\mathtt B_4$, $\mathtt B_4$, $\mathtt F_4$   \\ \hline
    \end{tabular}
    \caption{The rank 5 Lann\'er diagrams and maximal simplices.}
    \label{tab rank 5 Lanner}
\end{table}

\subsection{Pseudo-rank of a Coxeter--Dynkin diagram}
\label{sec.pseudo-rank}

We recall, that the \emph{pseudo-rank} of a (finite) Coxeter group $W_\mathtt{X}$ is defined to be the maximum
\[\max\left\{\rank(W_\mathtt{Y}) \mid \mathtt{Y} \text{ a Coxeter--Dynkin diagram with } W_\mathtt{Y}\cong W_\mathtt{X}\right\}. \]

\begin{lemma}
Let $\mathtt{X}_n$ be a Coxeter-Dynkin diagram from \Cref{fig:FiniteCoxeterGroups}. 
\begin{enumerate}
\item If $\mathtt{X}_n$ is not isomorphic to $\mathtt{B}_{2k+1}$ or $\mathtt{I}_2(4k+2)$ for $k\geq 1$, then $\rank(W_{\mathtt{X}_n})=\pseudorank(W_{\mathtt{X}_n})=n$.
\item If $\mathtt{X}_n=\mathtt{B}_{2k+1}$, then $\pseudorank(W_{\mathtt{X}_n})=2k+2$.
\item If $\mathtt{X}_n=\mathtt{I}_2(2(2k+1))$ for $k\geq 1$, then $\pseudorank(W_{\mathtt{X}_n})=3$.
\end{enumerate}
\end{lemma}
\begin{proof}
This follows from \cite[Theorem 2.17 and 3.3]{Nuida2006} (see also \cite{Paris2007}).  More precisely, $W_{\mathtt{X}_n}$ is directly decomposable in products of Coxeter groups if and only if $\mathtt{X}_n$ is of type $\mathtt{B}_{2k+1}$ or $\mathtt{I}_2(2(2k+1))$, $k\geq 1$ and the direct decomposition of these groups is as follows: $W_{\mathtt{B}_{2k+1}}\cong \Z_2\times W_{\mathtt{D}_{2k+1}}$, $W_{\mathtt{I}_2(2(2k+1))}\cong \Z_2 \times W_{\mathtt{I}_2(2k+1)}$.
\end{proof}

\begin{remark}
    The proofs in this section use extensive results about the subgroup structure of finite Coxeter groups.  Our data source is the abstract groups database on the $L$-functions and modular forms database \cite{lmfdb}.  In the database, the relevant Coxeter groups have the following id's: $\mathtt{H}_3$ is \texttt{120.35}, $\mathtt{A}_4$ is \texttt{120.34}, $\mathtt{B}_4$ is \texttt{384.5602}, $\mathtt{D}_4$ is \texttt{192.1493}, $\mathtt{F}_4$ is \texttt{1152.157478}, $\mathtt{A}_5$ is \texttt{720.763}, $\mathtt{B}_5$ is \texttt{3840.ch}, $\mathtt{D}_5$ is \texttt{1920.240996}, and $\mathtt{E}_6$ is \texttt{51840.b}.  For $\mathtt H_4$ we use direct computations in MAGMA.
\end{remark}

The following lemma is well known but can easily be deduced from the above data source.

\begin{lemma}\label{solvable.finite.Coxeter.Groups}
Let $\mathtt{X}$ be a connected Coxeter-Dynkin diagram such that $W_\mathtt{X}$ is finite. Then $W_\mathtt{X}$ is solvable if and only if $\mathtt{X}$ is isomorphic to one of the following Coxeter-Dynkin diagrams: $\mathtt{A}_1$, $\mathtt{A}_2$, $\mathtt{A}_3$, $\mathtt{B}_2$, $\mathtt{B}_3$, $\mathtt{B}_4$, $\mathtt{D}_4$, $\mathtt{F}_4$, $\mathtt{G}_2$, $\mathtt{I}_2(m)$. 
\end{lemma}

\begin{table}[]
\makebox[\linewidth][c]{%
\begin{tabular}{|c|c|}
\hline
Diagram     & Soluble subgroup                      \\ \hline
$\mathtt{E}_{6,7,8}$ & contains $\Z_2 \cdot (\sym(3)\wr \Z_2 )$ \\ \hline
$\mathtt{A}_5$       & maximal $\sym(3)\wr \Z_2 $              \\ \hline
$\mathtt{B}_5$       & maximal $\Z_2\wr F_5$                  \\ \hline
$\mathtt{D}_5$       & maximal $\Z_2^4 \rtimes F_5$           \\ \hline
$\mathtt{A}_4$       & maximal $F_5$                         \\ \hline
$\mathtt{H}_4$       & maximal $D_5^2\rtimes \Z_4$     \\ \hline
$\mathtt{H}_3$       & maximal $\Z_2\times \alt(4)$ and $D_{10}$         \\ \hline
\end{tabular}

}
\caption{Soluble subgroups which distinguish irreducible finite (non-soluble) Coxeter groups up to rank $4$.}
\label{tab maximal sol subs}
\end{table}

\begin{lemma}
    Let $\mathtt{X}$ be a Coxeter--Dynkin diagram such that $G=W_{\mathtt{X}}$ is finite and $\mathtt{X}$ has rank at most $4$.  Then, there is a list of maximal soluble subgroups of $G$ that distinguishes $G$ from all other finite Coxeter groups.
\end{lemma}
\begin{proof}
    We claim the information in \Cref{tab maximal sol subs} is sufficient to determine a finite Coxeter group of Coxeter--Dynkin rank at most $4$ up to isomorphism.  To do this we are using two reductions, firstly a soluble group is clearly determined by itself.  Secondly, we may reduce to the irreducible case since maximal soluble subgroups of direct products of Coxeter groups are exactly direct products of the maximal soluble subgroups of the factors.  Now, one can clearly determine the property of $\mathtt{X}$ being rank at most $4$ using the data in \Cref{tab maximal sol subs}.  Thus, it remains to distinguish the groups that are rank at most $4$.  In this case they are all distinguished, either by being soluble, or from the data in \Cref{tab maximal sol subs}. 
\end{proof}

\begin{lemma}\label{Lanner rk 4 - CFsol to CFmax}
    Let $n=4,5$, $\mathtt{\Gamma}$, $\mathtt{\Omega}$ be Coxeter-Dynkin diagrams and  $W_\mathtt{\Gamma}$ and $W_\mathtt{\Omega}$ be the associated Coxeter groups. Suppose $\mathtt{\Gamma}$ has no rank $n$ or pseudo-rank $n$ finite subdiagrams.  If there is a poset isomorphism $\CFsol(W_\mathtt{\Omega})\to\CFsol(W_{\mathtt{\Gamma}})$, then $\mathtt{\Omega}$ has no pseudo-rank $n$ spherical subdiagrams and the multiset of pseudo-rank $n-1$ spherical subdiagrams is equal to that of $\mathtt{\Gamma}$.
\end{lemma}
\begin{proof}
    Suppose not, then $W_\mathtt{\Omega}$ has a conjugacy class $[X]$ of maximal finite subgroups not contained in $W_\mathtt{\Gamma}$.  We may assume that $X$ is non-soluble, since otherwise there would be a maximal soluble subgroup in $\CFsol(W_\mathtt{\Omega})$ which is not contained in $\CFsol(W_\mathtt{\Gamma})$.  Suppose $n=5$.  We see that either $X$ has rank $4$ and is equal to $W_{\mathtt{A}_4}$ or $W_{\mathtt{H}_4}$, or has rank at least $5$.
    
    In the rank $5$ or higher case, by \Cref{tab maximal sol subs} there exists a subgroup in $\CFsol(W_\mathtt{\Omega})$ which is not isomorphic to a subgroup of $\CFsol(W_{\mathtt{\Gamma}})$ so we may rule out this possibility.  
    
    It remains to treat the rank $4$ case armed with the knowledge there are no rank $5$ finite subdiagrams in $\Omega$.  Again \Cref{tab maximal sol subs} allows us to distinguish the various groups involved so we just have to check the multisets of pseudo-rank $4$ subdiagrams are equal.   Essentially the only situation that can go wrong now is that we have several conjugacy classes of $X$, but their maximal soluble subgroups are all conjugate and so $\CFsol(W_\mathtt{\Omega})$ cannot distinguish them. 

    Suppose $X\cong W_{\mathtt{A}_4}$.  In this case we would have two conjugacy classes of maximal $\sym(5)$ subgroups $[S_1]$ and $[S_2]$ of $W_\mathtt{\Omega}$ and maximal-soluble $5$-Frobenius subgroups $F_1\leqslant S_1$ and $F_2\leqslant S_2$ such that $F_1^g=F_2$.  But now, the parabolic subgroup $S_1^g \cap S_2$ contains $F_2$ and any parabolic subgroup which contains a copy of $5$-Frobenius group has rank at least $4$.  Whence, $S_1^g=S_2$.

    The case where $X\cong W_{\mathtt{H}_4}$ is entirely analogous, instead using the maximal-soluble subgroup $D_5^2\rtimes \Z_4$ (small group ID $\langle400,129\rangle$).  

    Similarly, the case where $X\cong W_{\mathtt{A}_1\times \mathtt{H}_3}$ is analogous, instead using the maximal finite soluble subgroups $\Z_2^2 \times \alt(4)$ and $\Z_2^2\times D_{10}$.
    
    We are now in the case $n=4$ and armed with the knowledge that there are no pseudo-rank $5$ subdiagrams.  \Cref{tab maximal sol subs} and the poset isomorphism allow us to rule out the existence of any pseudo-rank $4$ groups. Since all pseudo-rank $3$ groups are soluble except $W_{\mathtt{H}_3}$, by the same arguments as above, it remains to rule out the case where we undercount conjugacy classes of $X\cong W_{\mathtt{H}_3}$.  But here the same argument as in the previous paragraphs applies to the subgroup $\Z_2 \times \alt(4)$.  Whence the lemma.
\end{proof}

\subsection{Detecting the virtual cohomological dimension}
By a result of Davis \cite[Corollary 8.5.5]{Davis2008},
the virtual cohomological dimension of a Coxeter system $(W,S)$ can be calculated using the nerve $\mathcal{N}(W,S)$. The \emph{nerve} of a Coxeter system $(W,S)$ is the simplicial complex with vertices $S$ and $n$-simplices given by all spherical subsets of $S$ of size $n$.

\begin{thm}[Davis]\label{thm vcd Cox}
    Let $(W,S)$ be a Coxeter system.  Then,
    \[\vcd(W)=  \max \left\{n \mid \widetilde H^n\left(\mathcal{N}(W_{S-T},S-T);\Z\right)\neq 0  \right\}\]
    where $T$ ranges over all subsets of $S$ such that $W_T$ is finite.
\end{thm}

\begin{lemma}\label{Lanner rk 4 - maximal cells}
    Let $W_\mathtt{\Gamma}$ be a Coxeter group such that $\mathtt{\Gamma}$ is a rank $n$ Lann\'er diagram, $n=4,5$.  If $W_\mathtt{\Omega}$ is another Coxeter group with $\widehat{W_\mathtt{\Omega}}\cong \widehat{W_\mathtt{\Gamma}}$, then $V(\mathtt{\Omega})$ is $2$-spherical, $W_\mathtt{\Omega}$ is a word-hyperbolic Coxeter group, and any Dynkin diagram for $W_\mathtt{\Omega}$ has no pseudo-rank $m$ spherical subdiagrams, for any $m \geq n$.  Moreover, the multiset of pseudo rank $n-1$ spherical subdiagrams of $\mathtt{\Omega}$ is equal to that of $\mathtt{\Gamma}$.
\end{lemma}

\begin{proof}
    That $V(\mathtt{\Omega})$ is $2$-spherical follows from the profinite invariance of property $\FA$ for Coxeter groups, that is \Cref{profiniteinvariantFA}.  That $W_\mathtt{\Omega}$ is a Gromov-hyperbolic group is given by \Cref{profiniteinvariantHyperbolicity}.

    Now, both $W_\mathtt{\Omega}$ and $W_\mathtt{\Gamma}$ are Gromov-hyperbolic and virtually compact special.  Hence, by \Cref{prop:CF}\eqref{CF:sol} we have that a poset isomorphism $\CFsol(W_\mathtt{\Gamma})\cong\CFsol(W_\mathtt{\Omega})$ exists.  Applying \Cref{Lanner rk 4 - CFsol to CFmax} we see that any Coxeter--Dynkin diagram for $W_\mathtt{\Omega}$ has no pseudo-rank $m$ spherical subdiagrams, for any $m \geq n$. 
\end{proof}

\begin{lemma}\label{Lanner rk 4 vcd}
    Let $W_\mathtt{\Gamma}$ be a Coxeter group such that $\mathtt{\Gamma}$ is a rank $n$ Lann\'er diagram $n=4,5$.  If $W_\mathtt{\Omega}$ is another Coxeter group with $\widehat{W_\mathtt{\Omega}}\cong \widehat{W_\mathtt{\Gamma}}$, then $\vcd W_\mathtt{\Omega}= n-1$ and $\vb_{n-1}(W_\mathtt{\Omega})=1$
\end{lemma}
\begin{proof}
    Each Coxeter group on a rank $n$ Lann\'er diagram corresponds to a closed hyperbolic $(n-1)$-orbifold which admits a closed orientable hyperbolic $(n-1)$-manifold as a finite cover.  In particular, there exists a finite index subgroup $H\leqslant W_\mathtt{\Gamma}$ with $H^{n-1}(H;\FF_p)= \FF_p$ for every prime $p$.  
    
    Note that since $W_\mathtt{\Gamma}$ is good (\Cref{CoxGpGood}), so is $H$. Now, $H$ is torsion-free so it follows from goodness that its profinite completion $\widehat H$ is also torsion-free and by the computation in the previous paragraph we have
    \[H^{n-1}(H;\FF_p)\cong H^{n-1}(\widehat H;\FF_p)\cong \FF_p.\]
    Thus, the $p$-cohomological dimension $\cd_p \widehat H$ is at least $n-1$ (in fact it is exactly $n-1$ but we do not need this). 
    
    Let $H_\mathtt{\Omega}$ be the subgroup corresponding to $H$ in $W_\mathtt{\Omega}$ under the isomorphism of profinite completions $\widehat{W_\mathtt{\Gamma}}\cong \widehat{W_\mathtt{\Omega}}$. By \Cref{CoxGpGood}, $W_\mathtt{\Omega}$ and $H_\mathtt{\Omega}$ are good.  Moreover, $\widehat H\cong \widehat{H_\mathtt{\Omega}}$, so by goodness we see that 
    \[H^{n-1}(\widehat H;\FF_p)\cong H^{n-1}(H_\Omega;\FF_p)=\FF_p. \]
    Thus, we have $\vcd W_\mathtt{\Omega} \geq n-1$.  Now, by \Cref{Lanner rk 4 - maximal cells} the nerve $N\mathtt{\Omega}$ contains no cells of dimension at least $n$.  So we conclude $\vcd W_\mathtt{\Omega}=n-1$. Finally, the claim $\vb_{n-1}(W_\mathtt{\Omega})=1$ follows from goodness and the fact that $\vb_{n-1}(W_\mathtt{\Gamma})=1$ because it is virtually the fundamental group of a closed aspherical $(n-1)$-manifold.
\end{proof}

\subsection{Finishing the proof}

\begin{prop}\label{Lanner rk 4 - Wrigid}
    Let $W_\mathtt{\Gamma}$ be a Coxeter group such that $\mathtt{\Gamma}$ is a rank $n=4,5$ Lann\'er diagram.  If $W_\mathtt{\Omega}$ is another Coxeter group with $\widehat{W_\mathtt{\Omega}}\cong \widehat{W_\mathtt{\Gamma}}$, then $W_\mathtt{\Omega}\cong W_\mathtt{\Gamma}$.
\end{prop}
\begin{proof}
    By \Cref{Lanner rk 4 - maximal cells} we have that $W_\mathtt{\Omega}$ is a $2$-spherical Coxeter group admitting a diagram $\mathtt{\Omega}$, such that $N\mathtt{\Omega}$ has no $n$-cells and at most $n$ many $(n-1)$-cells.  By \Cref{Lanner rk 4 vcd} and Davis' Theorem on the virtual cohomological dimension of Coxeter groups (\Cref{thm vcd Cox}), there exists a possibly empty spherical subset $T\subseteq\mathtt{\Omega}$ such that $H^{n-1}(N(\mathtt{\Omega}\setminus T);\Z)\neq 0$.  It follows immediately that the $n$ many $(n-1)$-cells in $N\mathtt{\Omega}$ form the boundary of an $(n-1)$-simplex $\partial\Delta^{n-1}$.  Let $S'=\{v_1,\dots,v_n\}$ denote the vertices spanning this subcomplex.

    We now suppose for contradiction that $N\mathtt{\Omega}$ contains more than $n$ vertices.  In this case $\partial\Delta^n$ is a full subcomplex of $N\mathtt{\Omega}$.  In particular, by \Cref{lem:affine-free_parabs}, $W_{S'}$ is a quasiconvex subgroup of the virtually compact special hyperbolic group $W_\mathtt{\Omega}$.  But $\vb_{n-1}(W_{S'})=1$, so \Cref{lem vb infty} applies and we conclude $\vb_{n-1}(W_\mathtt{\Omega})=\infty$.  This contradicts \Cref{Lanner rk 4 vcd}. So we must have $N\mathtt{\Omega}=\partial\Delta^{n-1}$.  
    
    In this case, the Coxeter diagram is completely determined by the set of spherical $(n-1)$-cells.  Indeed, the only other complete Coxeter groups with this nerve are affine Coxeter groups.  But these groups are not hyperbolic.  That the Lann\'er diagrams are distinguished from each other follows readily from \Cref{rank 3 lanner} and \Cref{tab rank 5 Lanner}.
\end{proof}

\section{Rigidity of special types of Coxeter groups}\label{sec Special types}
In this last section we prove (almost) rigidity results for special types of Coxeter groups. 
\subsection{Gromov-hyperbolic FC type}

\begin{prop}\label{prop:CFhyperbolicFC}
Let $\Gamma$ be a Coxeter graph and $W_\Gamma$ be the associated Coxeter group. If $W_\Gamma$ is of $FC$ type, then the canonical inclusion $\iota\colon W_\Gamma\to\widehat{W_\Gamma}$ induces an epimorphism $\varphi\colon\CF(W_\Gamma)\twoheadrightarrow\CF(\widehat{W_\Gamma})$.
In particular, if $W_\Gamma$ is Gromov-hyperbolic and of FC type, then the induced map $\varphi\colon\CF(W_\Gamma)\to\CF(\widehat{W_\Gamma})$ is an order isomorphism.

\end{prop}
\begin{proof}
If $\Gamma$ is a clique, then $W_\Gamma$ is finite and the conclusion of the proposition follows immediately. Thus, we may assume that $\Gamma$ is not a clique. By \Cref{amalgam.decompostion.Coxeter} the Coxeter group $W_\Gamma$ is an amalgamated product $W_\Gamma\cong A*_C B$ of special parabolic subgroups $A, B, C$. Special parabolic subgroups of $W_\Gamma$ are virtual retracts of $W_\Gamma$, therefore we can apply \Cref{AmalgamProfCompletion} to obtain $\widehat{W_\Gamma}\cong\widehat{A}\coprod_{\widehat{C}} \widehat{B}$.

Let $[F]\in\mathcal{CF}(\widehat{W_\Gamma})$. Then $F$ is contained in a conjugate of $\widehat{A}$ or $\widehat{B}$ by \cite[Theorem 4.18]{Ribes2017} using profinite Bass-Serre theory. If $A$ is not a clique, then we decompose $A$ again into an amalgamated product. Repeating this process finitely many times we obtain that $F$ is contained in a conjugate of $\widehat{A'}$ where $A'$ is a clique and hence a finite subgroup. Thus, $\widehat{A'}=A'$ and $F\subseteq gA'g^{-1}$. In particular, there exists a finite subgroup $A^{''}\subseteq W_\Gamma$ such that $F=gA^{''}g^{-1}$. Hence, $\varphi([A^{''}])=[F]$, which shows the surjectivity of $\varphi$. Clearly, by construction $\varphi$ is order preserving.

Assume now that $W_\Gamma$ is additionally Gromov-hyperbolic. Since Gromov-hyperbolic Coxeter groups are virtually compact special and virtually toral relatively hyperbolic,
the induced map $\varphi\colon\CF(W)\to\CF(\widehat{W})$ is also injective by \Cref{{prop:CF}}(4).  
Hence $\varphi\colon\CF(W_\Gamma)\to\CF(\widehat{W_\Gamma})$ is an order isomorphism.
\end{proof}

\begin{corollary}\label{almost hyp FC type}
Let $W_\Gamma$ be a Gromov-hyperbolic Coxeter group. If $W_\Gamma$ is of FC type, then $W_\Gamma$ is almost $\mathcal{W}$-profinitely rigid.
\end{corollary}
\begin{proof}
Let $W_\Omega$ be a Coxeter group such that $\widehat{W_\Gamma}\cong\widehat{W_\Omega}$.  
Similar to the proof of \Cref{thmx:Reduction1ended}, to prove the corollary it suffices to find a bound on $|V(\Omega)|$ in terms of the combinatorics of $V(\Gamma)$.  

By \Cref{prop:CFhyperbolicFC} we have
$\CF(W_\Gamma)=\CF(\widehat{W_\Gamma}).$ Since being Gromov-hyperbolic is a $\mathcal{W}$-profinite invariant by \Cref{profiniteinvariantHyperbolicity}, $W_\Omega$ is Gromov-hyperbolic and we have
$\CF(W_\Omega)\hookrightarrow \CF(\widehat{W_\Omega})$ by \Cref{prop:CF}(4).
Hence
$$\CF(W_\Gamma)=\CF(\widehat{W_\Gamma})=\CF(\widehat{W_\Omega})\hookleftarrow \CF(W_\Omega)$$
The order monomorphism $\CF(W_\Omega)\hookrightarrow\CF(W_\Gamma)$ shows that the graph $\Omega$ can be covered by defining graphs of some representatives $[A]\in \CF(W_\Gamma)$ where $A$ is a Coxeter group. Let $\CF(W_\Gamma)=\left\{[A_1],\ldots,[A_n]\right\}$ and define $\mathcal{C}:=\left\{A_i\mid A_i\text{ is a Coxeter group}\right\}$. Then 
$$|V(\Omega)|\leq \sum_{C\in\mathcal{C}}\pseudorank(C).$$
By \Cref{cor:almostRigidity} follows that $W_\Gamma$ is almost $\calw$-profinitely rigid.
\end{proof}

\subsection{Extra large type}\label{sec XL}
A Coxeter group $W_\Gamma$ is of \emph{extra large type} if every edge label in $\Gamma$ is at least $4$.  Note that such a Coxeter group is hyperbolic by \Cref{NoZsquare}.

\begin{lemma}\label{extralargetype.CF}
Let $W_\Gamma$ and $W_\Omega$ be connected Coxeter groups such that $\widehat{W_\Gamma}\cong\widehat{W_\Omega}$. If $W_\Gamma$ is of extra large type, then 
$$\CF(W_\Gamma)=\CF(\widehat{W_\Gamma})=\CF(\widehat{W_\Omega})=\CF(W_\Omega).$$ 
Assume additionally that $m(e)\neq 4k+2$, $k\geq 1$ for every $e\in E(\Gamma)$. Then $W_\Omega$ is also of extra large type.
\end{lemma}
\begin{proof}
Since $W_\Gamma$ is of extra large type it is hyperbolic.  Thus, $W_\Omega$ is also hyperbolic by \Cref{profiniteinvariantHyperbolicity} and hence is good.  Since $m(E(\Gamma))\subseteq\mathbb{N}_{\geq 4}$ it follows from the classification of finite irreducible Coxeter groups, see \Cref{fig:FiniteCoxeterGroups}, that maximal finite subgroups of $W_\Gamma$ are special parabolic subgroups $W_\Delta$ where $\Delta$ is an edge in $\Gamma$. In particular, $\CF(W_\Gamma)=\CFsol(W_\Gamma)$. 

Let $H\subseteq W_\Gamma$ be a finite subgroup. Then $H$ is contained in a conjugate of $W_\Delta$ where $\Delta=(\left\{v,w\right\}, \left\{\left\{v,w\right\}\right\})$ and $m(\left\{v,w\right\})=l$, so $W_\Delta$ is isomorphic to the Dihedral group $D_l$.
The non-trivial subgroups of $D_l$ are of two types: cyclic or Dihedral. 
Hence every representative of an element in $\CF(W_\Gamma)$ is trivial or cyclic or Dihedral.

Since $m(E(\Gamma))\subseteq\mathbb{N}_{\geq 4}$ it follows from \Cref{NoZsquare} that $W_\Gamma$ is Gromov-hyperbolic and therefore $W_\Omega$ is also Gromov-hyperbolic and by \Cref{prop:CF} we have 
$$\CFsol(W_\Gamma)=\CFsol(W_\Omega).$$
Our goal now is to show that $\CFsol(W_\Omega)=\CF(W_\Omega)$. Assume towards a contradiction that $\CFsol(W_\Omega)\neq \CF(W_\Omega)$. Then $W_\Omega$ has a special parabolic subgroup isomorphic to  $\mathtt{A}_n$, $n\geq 4$ or $\mathtt{B}_n$, $n\geq 5$ or $\mathtt{D}_n$, $n\geq 5$, $\mathtt{E}_6$, $\mathtt{E}_7$, $\mathtt{E}_8$, $\mathtt{H}_3$ or $\mathtt{H}_4$ by \Cref{solvable.finite.Coxeter.Groups}. All these non-solvable finite Coxeter groups have non-cyclic and non-Dihedral solvable subgroups by \Cref{tab maximal sol subs}. Note that $\mathtt{A}_4\subseteq \mathtt{A}_n$ for $n\geq 4$, $\mathtt{B}_5\subseteq \mathtt{B}_n$ for $n\geq 5$ and $\mathtt{D}_5\subseteq \mathtt{D}_n$ for $n\geq 5$. Since $\CFsol(W_\Omega)$ has only cyclic and Dihedral subgroups, it follows that $W_\Omega$ can not have a finite parabolic non-solvable subgroup. Hence $\CF(W_\Omega)=\CFsol(W_\Omega)$.
Thus we obtain

$$\CF(W_\Gamma)=\CF(\widehat{W_\Gamma})=\CF(\widehat{W_\Omega})=\CF(W_\Omega).$$

In particular $\CFmax(W_\Gamma)=\CFmax(W_\Omega)$.
If $m_{\Gamma}(e)\neq 4k+2$ for $e\in E(\Gamma)$, then the defining graphs for the representatives in $\CFmax(W_\Omega)$ are edges with labels $\geq 4$ by \Cref{DihedralGroups.Graphrigid}. Therefore, $W_\Omega$ is of extra large type.
\end{proof}

\begin{prop}\label{XL almost}
Coxeter groups of extra large type are almost $\mathcal{W}$-profinitely rigid.
\end{prop}
\begin{proof} 
Let $W_\Gamma$ be a Coxeter group of extra large type and $W_\Omega$ be a Coxeter group such that $\widehat{W_\Gamma}\cong\widehat{W_\Omega}$.  Since $W_\Gamma$ is hyperbolic we see that $W_\Omega$ is also hyperbolic by \Cref{profiniteinvariantHyperbolicity}.   Without loss of generality we can assume that $\Gamma$ and $\Omega$ are connected graphs. By \Cref{extralargetype.CF} we have
$$\CF(W_\Gamma)=\CF(W_\Omega).$$
Note that the Coxeter graph $\Omega$ can be covered by defining graphs of representatives in $\CFmax(W_\Omega)$. Each representative in $\CFmax(W_\Omega)$ is  isomorphic to a Dihedral group and Dihedral groups have pseudo-rank $\leq 3$.
Hence $|V(\Omega)|\leq 3\cdot |\CFmax(W_\Gamma)|$. By \Cref{cor:almostRigidity} follows that $W_\Gamma$ is almost $\calw$-profinitely rigid. 
\end{proof}

\subsection{Strongly even Coxeter groups}\label{sec even}

A Coxeter group $W_\Gamma$ is called \emph{strongly even} if $E(\Gamma)=\emptyset$ or $m(E(\Gamma))\subseteq\left\{2\right\}\cup 4\mathbb{N}$. 

\begin{prop}\label{strongly.even.Coxeter.CF}
Let $W_\Gamma$ and $W_\Omega$ be strongly even Coxeter groups. The following statements are equivalent:
\begin{enumerate}
\item $W_\Gamma\cong W_\Omega$
\item $\Gamma\cong \Omega$
\item $\CF(W_\Gamma)=\CF(W_\Omega)$
\end{enumerate}
\end{prop}
\begin{proof}
For right-angled Coxeter groups the result was proven in \cite[Theorem 2.5]{CorsonHughesMollerVarghese2023}. Furthermore, the equivalence of (1) and (2) for strongly even Coxeter groups was shown in \cite[Theorem 4.11]{Radcliffe2001}. Following the proof strategy of \cite[Theorem 2.5]{CorsonHughesMollerVarghese2023} and using \cite[Theorem 4.11]{Radcliffe2001} we obtain the equivalence of (2) and (3).
\end{proof}

\begin{thm}\label{thm:even.profinite}
Let $W_\Gamma$ be a Coxeter group. If $m(E(\Gamma))\subseteq 4\mathbb{N}$, then $W_\Gamma$ is $\mathcal{W}$-profinitely rigid. 
\end{thm}
\begin{proof}
Let $W_\Omega$ be a Coxeter group such that $\widehat{W_\Gamma}\cong\widehat{W_\Omega}$. Without loss of generality we can assume that $\Gamma$ and $\Omega$ are connected graphs, see \cite[Theorem 3.10]{CorsonHughesMollerVarghese2023}.  Then $\CF(W_\Gamma)=\CF(W_\Omega)$ by \Cref{extralargetype.CF}. In particular, maximal finite subgroups in $W_\Omega$ are Dihedral groups or order $2\cdot 4k$. These Dihedral groups are graph rigid by \Cref{DihedralGroups.Graphrigid}, therefore every edge label in $\Omega$ is divisible by $4$. 
Now we apply \Cref{strongly.even.Coxeter.CF} and we obtain $W_\Gamma\cong W_\Omega$.
\end{proof}

\subsection{Odd Coxeter groups}\label{Section:OddCoxeterGroups}

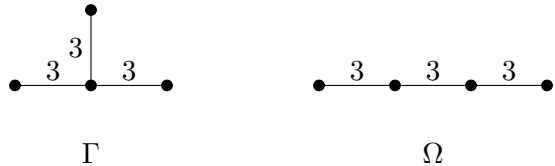
\begin{figure}[htb]
	\begin{center}
	\captionsetup{justification=centering}
		\begin{tikzpicture}
			\draw[fill=black]  (0,0) circle (2pt);
            \draw[fill=black]  (1,0) circle (2pt);
            \draw[fill=black]  (2,0) circle (2pt);
            \draw[fill=black]  (1,1) circle (2pt);
            \draw (0,0)--(2,0);
            \draw (1,0)--(1,1);
            \node at (0.5, 0.2) {$3$};
            \node at (1.5, 0.2) {$3$};
            \node at (0.8, 0.5) {$3$};
            \draw[fill=black]  (4,0) circle (2pt);
            \draw[fill=black]  (5,0) circle (2pt);
            \draw[fill=black]  (6,0) circle (2pt);
            \draw[fill=black]  (7,0) circle (2pt);
            \draw (4,0)--(7,0);
            \node at (4.5, 0.2) {$3$};
            \node at (5.5, 0.2) {$3$};
            \node at (6.5, 0.2) {$3$};
           \node at (1, -0.9) {$\Gamma$};
            \node at (5.5, -0.9) {$\Omega$};           
        \end{tikzpicture}
        \caption{Two non-isomorphic Coxeter diagrams with isomorphic Coxeter groups, first observed by M\"uhlherr \cite{Muhlherr2000}.}
	      \label{fig:noNonRigidOddCoxeterGroups}
    \end{center}
\end{figure}

Let $\Gamma$ be a Coxeter graph. The associated Coxeter group $W_\Gamma$ is called \emph{odd} if all edge labels in $\Gamma$ are odd.  Note that odd Coxeter groups in general are not determined by their Coxeter graph, see \Cref{fig:noNonRigidOddCoxeterGroups}.  In this section, we prove two theorems: the first states that Coxeter groups with Coxeter graph an odd labelled forest are profinitely rigid amongst Coxeter groups (\Cref{oddForest_rigid})); the second states that odd Coxeter groups are almost profinitely rigid amongst Coxeter groups.  
We denote by $\cc_2(W)$ the number of conjugacy classes of involutions in $W$.

\begin{prop}\label{props odd}
    Let $W_\Gamma$ and $W_\Omega$ be Coxeter groups.  The following conclusions hold:
    \begin{enumerate}
        \item \cite[Corollary~1.3(1)]{MichaelSantosRegosSchwerVarghese2024} \label{CC2=1} A Coxeter group $W_\Gamma$ with $\Gamma$ connected is odd if and only if $\cc_2(W_\Gamma)=1$.
    \end{enumerate}
    Suppose additionally $W_\Gamma$ is odd.
    \begin{enumerate}[resume]
    \item If $W_\Gamma\cong W_\Omega$, then $W_\Omega$ is odd. \label{OddCoxeterGroupsCC2}
    \item If $\widehat{W_\Gamma}\cong \widehat{W_\Omega}$, then $W_\Omega$ is odd.\label{OddCoxeterProf}
    \item $W_\Gamma$ is Gromov-hyperbolic if and only if $\Gamma$ contains no $\Delta(3,3,3)$ triangle subgraphs.\label{OddHyp}
    \item $W_\Gamma$ is virtually free if and only if $\Gamma$ is a tree.\label{OddVirFree}
    \end{enumerate}
    Suppose additionally $\Gamma$ is an odd labelled tree.
    \begin{enumerate}[resume]
    \item \cite[Example 5.1]{BradyMcCammondMühlherrNeumann2002} Then, $W_\Gamma\cong W_\Omega$ if and only if $\Omega$ is a tree and the multiset of edge labels for $\Gamma$ is the same as the multiset of edge labels for $\Omega$.\label{OddTrees_twisting}
    \item There is an order isomorphism $\CF(W_\Gamma)\to\CFsol(\widehat {W_\Gamma})$.\label{oddTrees_CF}
    \item If $\widehat{W_\Gamma}\cong \widehat{W_\Omega}$, then $\Omega$ is an odd-tree and $\CF(W_\Gamma)=\CF(W_\Omega)$. \label{odd_trees_CF_prof}
    \end{enumerate}
\end{prop}
\begin{proof}
    The first claim \eqref{CC2=1} is as cited so there is nothing to prove.  To see \eqref{OddCoxeterGroupsCC2} note that Coxeter group $W_\Gamma$ is odd if and only if every connected component of $\Gamma$ is odd. The number of connected components is an isomorphism invariant since a Coxeter group splits as a free product if and only if the Coxeter graph is disconnected. Now, apply \eqref{CC2=1} to every connected component.

    We now prove \eqref{OddCoxeterProf}.  By \cite[Theorem~3.10]{CorsonHughesMollerVarghese2023}, we may assume without loss of generality that $\Gamma$ is connected.  By \Cref{prop:CF}\eqref{CF:p}, we have order isomorphisms \[\CF_2(W_\Gamma)\to \CF_2(\widehat{W_\Gamma}) \to \CF_2(\widehat{W_\Omega})\xleftarrow{} \CF_2(W_\Omega).\]  Hence, $\cc_2(W_\Omega)=1$ and so the claim follows from \eqref{OddCoxeterGroupsCC2}.

    Item \eqref{OddHyp} follows from Moussong's Theorem \cite{Moussong88}. Next, \eqref{OddVirFree} follows from the fact that a Coxeter group is virtually free if and only if the Coxeter graph is chordal and every complete parabolic subgroup is finite.  

    We explain \eqref{OddTrees_twisting} since we are claiming a slightly stronger result than in \cite[Example 5.1]{BradyMcCammondMühlherrNeumann2002}.  The key point is that loc. cit. assumes that $\Omega$ is also an odd-labelled tree.  However, the oddness assumption is superfluous by \eqref{OddCoxeterGroupsCC2}.
    
    We next prove \eqref{oddTrees_CF}. Since $W_\Gamma$ is virtually free, we see that $W_\Gamma$ is a hyperbolic group.   By \Cref{prop:CF}\eqref{CF:sol}, we have an order isomorphism $\CFsol(W_\Gamma)\to \CFsol(\widehat{W_\Gamma})$.  But, every finite subgroup of $W_\Gamma$ is soluble.  Hence, the claim.   Finally, \eqref{odd_trees_CF_prof} follows from \eqref{OddCoxeterProf}, \Cref{vir free W profinite}, \eqref{OddVirFree}, and then applying \eqref{oddTrees_CF} twice.
\end{proof}

Let $W_\Gamma$ be a connected odd Coxeter group. It follows from \cite{Brink1996} that the fundamental group of $\Gamma$ is an isomorphism invariant. In the next proposition we show that, if the fundamental group of $\Gamma$ is abelian, then it is also a $\mathcal{W}$-profinite invariant.

\begin{prop}
\label{oddFundamentalGroup}
Let $W_\Gamma$ and  $W_\Omega$ be odd Coxeter groups with $\Gamma,\Omega$ connected. We denote by $\pi_1(\Gamma)$ resp. $\pi_1(\Omega)$ the fundamental group of $\Gamma$ resp. $\Omega$. Assume that $\pi_1(\Gamma)$ is trivial or $\mathbb{Z}$. If $\widehat{W_\Gamma}\cong\widehat{W_\Omega}$, then
$\pi_1(\Gamma)\cong\pi_1(\Omega)$.
\end{prop}
\begin{proof}
By \Cref{props odd} $\cc_2(W_\Gamma)=\cc_2(W_\Omega)=1$. Let $v\in V(\Gamma)$ and $w\in W_\Gamma$ be an involution. Since $v$ and $w$ are conjugate, the centraliser of $v$ is conjugate to the centraliser of $w$. The structure of centralisers of standard generators of a Coxeter group was calculated in \cite{Brink1996}, hence $C_{W_\Gamma}(v)\cong A\rtimes \pi_1(\Gamma)$ where $A\subseteq W_\Gamma$ is a Coxeter group. Let $a_i$ be a Coxeter generator of $A$. Assume that $a_i\neq v$. Then the subgroup $\langle a_i, v\rangle$ is isomorphic to $\Z_2^2$ and therefore is contained in a conjugate of a special parabolic subgroup which is isomorphic to a Dihedral group of order $2m$, where $m$ is odd. Hence $4$ divides $2m$ which is impossible, thus $A=\langle v\rangle$ and $C_{W_\Gamma}(v)\cong \Z_2 \rtimes\pi_1(\Gamma)$.

Let $f\colon\widehat{W_\Gamma}\to\widehat{W_\Omega}$ be an isomorphism. By \Cref{prop:CF} we have
$C_{\widehat{W_\Gamma}}(v)=\overline{C_{W_\Gamma}(v)}$. Hence, if $\pi_1(\Gamma)$ is trivial or isomorphic to $\mathbb{Z}$, then $C_{\widehat{W_\Gamma}}(v)\cong\overline{\Z_2 \rtimes\pi_1(\Gamma)}\cong \Z_2 \rtimes\widehat{\pi_1(\Gamma)}$  by \cite[Lemma 4.5]{CapraceMinasyan2013}.

Now, $f(C_{\widehat{W_\Gamma}}(v))=C_{\widehat{W_\Omega}}(f(v))$. Hence, if $\pi_1(\Gamma)$ is trivial, then $C_{W_\Gamma}(v)\cong \Z_2 \cong C_{\widehat{W_\Gamma}}(v)$ and therefore  $C_{\widehat{W_\Omega}}(f(v))\cong \Z_2 $. Since $\CF_2(W_\Omega)=\CF_2(\widehat{W_\Omega})$ we conclude that the centraliser of every involution in $W_\Omega$ is finite and therefore $\pi_1(\Omega)$ is trivial.

If $\pi_1(\Gamma)$ is isomorphic to $\mathbb{Z}$, then $C_{\widehat{W_\Omega}}(f(v))\cong \Z_2 \rtimes\widehat{\mathbb{Z}}$. Thus, by the same argument as above the centraliser of every involution in $W_\Omega$ is isomorphic to $\Z_2 \rtimes\mathbb{Z}$. Hence $\pi_1(\Omega)\cong\mathbb{Z}$.
\end{proof}

We call a Coxeter group $W_\Gamma$ an \emph{odd forest Coxeter group} if each connected component of $\Gamma$ is an odd-labelled tree.

\begin{thm}\label{oddForest_rigid}
    If $W_\Gamma$ is an odd forest Coxeter group, then $W_\Gamma$ is profinitely rigid amongst Coxeter groups
\end{thm}
\begin{proof}
    By \cite[Theorem~3.10]{CorsonHughesMollerVarghese2023} we may assume that $\Gamma$ is a tree.   Since $W_\Gamma$ is hyperbolic by \Cref{profiniteinvariantHyperbolicity} we see that any Coxeter group in its genus is also hyperbolic (and hence good).  The result now follows immediately from items \eqref{oddTrees_CF} and \eqref{OddTrees_twisting} of \Cref{props odd}, noting that the multiset of edge labels of $\Gamma$ is completely determined by the maximal elements of the poset $\CF(W_\Gamma)$.
\end{proof}

Let $W_\Gamma$ be a Coxeter group. We define $A_2:=m^{-1}(\left\{2\right\})$, where $m\colon E(\Gamma)\to\mathbb{N}_{\geq 2}$ is the edge-labelling of $\Gamma$. For $\left\{v,w_1\right\}, \left\{v, w_2\right\}\in A_2$ we write $\left\{v,w_1\right\}\approx \left\{v, w_2\right\}$ if $m(\left\{w_1, w_2\right\})$ is odd. Let $\sim$ be the equivalence relation on $A_2$ generated by $\approx$.

We denote by $\nu(W_\Gamma)$ the number of equivalence classes of $\sim$ on $A_2$.  We denote by $\mu(W_\Gamma)$ the number of edges in $\Gamma$ with edge-label $\geq 3$.  The abelianization of $W_\Gamma$ is $\left(\Z_2 \right)^n$ for some $n$.  We define $\xi(W_\Gamma)=n$.  

\begin{thm}~\cite[Theorem A]{Howlett1988}\label{thm.Howlett}
Let $W_\Gamma$ be a Coxeter group. Then
$$M(W_\Gamma)=\Z_2 ^{\nu(W_\Gamma)+\mu(W_\Gamma)+\xi(W_\Gamma)-|V(\Gamma)|}$$
\end{thm}

\begin{corollary}
\label{Odd_MultProfinite}
Let $W_\Gamma$ and $W_\Omega$ be Coxeter groups. If $\widehat{W_\Gamma}\cong\widehat{W_\Omega}$, then 
$$\nu(W_\Gamma)+\mu(W_\Gamma)-|V(\Gamma)|= \nu(W_\Omega)+\mu(W_\Omega)-|V(\Omega)|.$$
If additionally $\Gamma$ is connected and odd, then 
\[|E(\Gamma)|-|V(\Gamma)|=|E(\Omega)|-|V(\Omega)|.\]
\end{corollary}
\begin{proof}
    The first claim follows from applying \Cref{M(W) prof} and \Cref{thm.Howlett}.  We now prove the second claim.  First, observe that by \cite[Theorem~3.10]{CorsonHughesMollerVarghese2023}, $\Omega$ is connected and by \Cref{props odd}\eqref{OddCoxeterProf}, $W_\Omega$ is odd.  For an odd Coxeter group $W_\Lambda$ note that $\nu(W_\Lambda)=0$ and that $\mu(W_\Lambda)=|E(W_\Lambda)|$.  The result now follows from the first claim.
\end{proof}

\begin{thm}
\label{OddCoxeterAlmost}
Odd Coxeter groups are almost $\mathcal{W}$-profinitely rigid.
\end{thm}
\begin{proof}
Let $W_\Gamma$ be an odd Coxeter group and $W_\Omega$ be a Coxeter group such that $\widehat{W_\Gamma}\cong\widehat{W_\Omega}$.  By \cite[Theorem~3.10]{CorsonHughesMollerVarghese2024}, without loss of generality we may assume that $\Gamma$ is connected and hence so is $\Omega$.  It follows from \Cref{props odd}\eqref{OddCoxeterProf}, that $W_\Omega$ is odd.   

By \Cref{Odd_MultProfinite} we have
\begin{equation}
    |E(\Omega)|=|E(\Gamma)|-|V(\Gamma)|+|V(\Omega)| \label{eqn.Odd.VE}
\end{equation}
By \cite{Chiswell1992} and profinite invariance of the Euler characteristic of Coxeter groups (\Cref{profInv_X_H*}) we have
\begin{align}
    \label{X_odd_W_Gamma} \chi(W_\Gamma)&=1-\frac{|V(\Gamma)|}{2}+\sum_{e\in E(\Gamma)}\frac{1}{2\cdot m_\Gamma(e)}\\
    \nonumber &=1-\frac{|V(\Omega)|}{2}+\sum_{e\in E(\Omega)}\frac{1}{2\cdot m_\Omega(e)}=\chi(W_\Omega),
\end{align}
where $m_\Gamma$ and $m_\Omega $ are the respective edge-labellings.  Now, $2\leq m_\Gamma(e)\leq d$ for all $e\in E(\Gamma)$ and similarly $2\leq m_\Omega(e)\leq d$ for all $e\in E(\Omega)$.  We now combine the extremes of these two inequalities, that is the case where $m_\Gamma(e)=d$ for all $e\in E(\Gamma)$, and the case where $m_\Omega(e)=2$ for all $e\in E(\Omega)$, with equation \eqref{X_odd_W_Gamma} to obtain
\begin{align}
\nonumber 1-\frac{|V(\Gamma)|}{2}+\frac{|E(\Gamma)|}{2d} & \leq 1-\frac{|V(\Gamma)|}{2}+\sum_{e\in E(\Gamma)}\frac{1}{2\cdot m(e)}\\
\nonumber &=1-\frac{|V(\Omega)|}{2}+\sum_{e\in E(\Omega)}\frac{1}{2\cdot m(e)}\\
\nonumber &\leq 1-\frac{|V(\Omega)|}{2}+|E(\Omega)|\cdot\frac{1}{2\cdot 2}.
\intertext{Thus, we obtain}
\nonumber 1-\frac{|V(\Gamma)|}{2}+\frac{|E(\Gamma)|}{2d} &\leq 1-\frac{|V(\Omega)|}{2}+\frac{|E(\Omega)|}{4}.
\end{align}
Rearranging gives
\begin{align}
    |E(\Omega)|&\geq 2\left(|V(\Omega)|-|V(\Gamma)|+\frac{1}{d}|E(\Gamma)|\right).\label{eqn.odd.VE.2}
\intertext{Substituting \eqref{eqn.odd.VE.2} into \eqref{eqn.Odd.VE} we obtain}
    |E(\Gamma)|-|V(\Gamma)|+|V(\Omega)| &\geq 2\left(|V(\Omega)|-|V(\Gamma)|+\frac{1}{d}|E(\Gamma)|\right). \nonumber\\
\intertext{Hence,}
    |V(\Omega)| &\geq 2|V(\Omega)|-|V(\Gamma)|+\left(\frac{2}{d}-1\right)|E(\Gamma)|,\nonumber
\intertext{then subtracting $2|V(\Omega)|$ and multiplying by $-1$ yields}
    |V(\Omega)|&\leq |V(\Gamma)| +\left(1-\frac{2}{d}\right)|E(\Gamma)|.
\end{align}
which is the desired inequality. By \Cref{cor:almostRigidity} follows that $W_\Gamma$ is almost $\calw$-profinitely rigid. 
\end{proof}

\subsection{Complete Coxeter groups}\label{sec complete}
\begin{lemma}
\label{complete.odd.vertices}
Let $W_\Gamma$ and $W_\Omega$ be Coxeter groups. Assume that $\widehat{W_\Gamma}\cong\widehat{W_\Omega}$. If $\Gamma$ is odd and complete, then $W_\Omega$ is odd, complete and $|V(\Gamma)|=|V(\Omega)|$.
\end{lemma}
\begin{proof}
    By \Cref{props odd} and \Cref{profiniteinvariantFA} follows that $\Omega$ is an odd complete graph.
    Let $|V(\Gamma)|=n$ and $|V(\Omega)|=m$. Wlog, we can assume that $n,m\geq 2$.
    We note, that a complete graph with $n$ vertices has $\frac{n\cdot(n-1)}{2}$ edges. By \Cref{Odd_MultProfinite} we have
    $$\frac{n\cdot(n-1)}{2}-n=\frac{m\cdot(m-1)}{2}-m$$
    where $|V(\Gamma)|=n$ and $|V(\Omega)|=m$. Since, for $x\geq 2$ the function $f(x)=x^2-3x$ is strictly increasing, we obtain $n=m$.
\end{proof}

\begin{prop}
\label{OddComplete.profinite}
Let $W_\Gamma$ be an odd  Coxeter group. If $\Gamma$ is complete such that $m(E(\Gamma))=\left\{n\right\}$, then $W_\Gamma$ is $\mathcal{W}$-profinitely rigid.
\end{prop}
\begin{proof}
Let $W_\Omega$ be a Coxeter group such that $\widehat{W_\Gamma}\cong\widehat{W_\Omega}$. 
By \Cref{complete.odd.vertices} the graph $\Omega$ is a complete odd-graph with $|V(\Gamma)|=|V(\Omega)|$. We note that maximal finite subgroups in $W_\Gamma$ and in $W_\Omega$ corresponds to labeled edges  in the graphs $\Gamma$ and $\Omega$. 

\begin{enumerate}
\item Let $n=3$.  It remains to show that $m(E(\Omega))=\left\{3\right\}$. 
Assume that there exists an edge $e\in E(\Omega)$ with edge label $l\neq 3$. Then there exist a prime number $p$ and $k\in\mathbb{N}$, such that $p^k|l$ but $p^k\nmid 3$, ($p\neq 2$, since $\Omega$ is odd). Thus $W_\Omega$ has an element or order $p^k$. By \Cref{prop:CF} $\CF_p$ is a $\mathcal{W}$-profinite invariant, hence $\CF_p(W_\Gamma)=\CF_p(W_\Omega)$. In particular, $W_\Gamma$ has an element of order $p^k$, which is impossible, since every torsion element is conjugate to an element of the Dihedral group $D_3$ that has order $6$ and $p^k\neq 2, 3$. Hence, every edge-label in $\Omega$ is equal to $3$ and $W_\Gamma\cong W_\Omega$.

\item If $n\geq 5$, then $W_\Gamma$ is of extra large type and by \Cref{extralargetype.CF} we have
$$\CF(W_\Gamma)=\CF(W_\Omega)$$

The conjugacy classes of maximal finite subgroups correspond to labeled edges in the graphs $\Gamma$ and $\Omega$. Thus  $m(E(\Omega))=\left\{n\right\}$ and therefore $\Gamma\cong \Omega$ and $W_\Gamma=W_\Omega$.\qedhere
\end{enumerate}
\end{proof}

\begin{prop}
\label{evenComplete.profinite}
Let $W_\Gamma$ be a Coxeter group. If $\Gamma$ is complete such that $m(E(\Gamma))=\left\{n\right\}$, $n$ even, $n\neq 4k+2$ for $k\geq 1$, then $W_\Gamma$ is $\mathcal{W}$-profinitely rigid.
\end{prop}
\begin{proof}
Let $W_\Omega$ be a Coxeter group such that $\widehat{W_\Gamma}\cong\widehat{W_\Omega}$. Since being a complete graph is a $\mathcal{W}$-profinite invariant by \Cref{profiniteinvariantFA}, it follows that $\Omega$ is a complete graph. If $n=2$, then $W_\Gamma$ is finite and therefore $W_\Gamma\cong W_\Omega$. Let $n\geq 4$. Then $W_\Gamma$ is of extra large type and by \Cref{extralargetype.CF} we have
$$\CF(W_\Gamma)=\CFsol(W_\Gamma)=\CFsol(W_\Omega)=\CF(W_\Omega).$$

Let $[A],[B]\in\CF_{max}(W_\Gamma)$. Then there are only two possibilities: 1) the greatest common lower bound of $[A],[B]$ is trivial or 2)  the greatest common lower bound of $[A],[B]$ is $[C]$ where $C\cong \Z_2 $.

Since $n\neq 4k+2$, the defining graphs of the maximal finite subgroups in $W_\Omega$ are edges. Hence $\Gamma\cong\Omega$.
\end{proof}

\begin{prop}\label{complete n neq 3 4k+2}
   Let $W_\Gamma$ be a Coxeter group be such that all edge-labels are equal to $n\neq 4k+2$ for $k\geq 1$.  If $\Gamma$ is complete, then $W_\Gamma$ is $\calw$-profinitely rigid.
\end{prop}
\begin{proof}
    This is the content of \Cref{OddComplete.profinite} and \Cref{evenComplete.profinite}.
\end{proof}

\subsection{Coxeter groups of pseudo-rank at most four}\label{sec rank 4}

\begin{prop}\label{rank at most 3 rigid}
Let $W_\Gamma$ be a Coxeter group. If $|V(\Gamma)|\leq 3$, then $W_\Gamma$ is $\calw$-profinitely rigid.
\end{prop}
\begin{proof}
Let $W_\Omega$ be a Coxeter group such that $\widehat{W_\Gamma}\cong\widehat{W_\Omega}$. Without loss of generality we can assume that $\Gamma$ is connected, see \cite[Theorem 3.10]{CorsonHughesMollerVarghese2023}.
If $|V(\Gamma)|\leq 2$, then $W_\Gamma$ is finite and therefore $W_\Gamma\cong W_\Omega$.

If $|V(\Gamma)|=3$, then we have to consider two cases:
\begin{enumerate}
\item Let $|E(\Gamma)|=2$. Then $W_\Gamma$ can be decomposed as a visual amalgamated product $A*_C B$ where $A, B$ are finite Dihedral groups and $C$ is of type $\mathtt{A}_1$. Hence by \Cref{amalgam.finitegroups.Outabelian} follows that $W_\Gamma\cong W_\Omega$. 
\item Let $|E(\Gamma)|=3$. Then $W_\Gamma$ is a triangle Coxeter group. Hence $W_\Gamma\cong W_\Omega$ by \Cref{triangleAndmore}(2).\qedhere
\end{enumerate}
\end{proof}

\begin{thm}
\label{CoxeterGroups4Vertices}
Let $W_\Gamma$ be a connected Coxeter group.  If $|V(\Gamma)|=4$ and $m(E(\Gamma))=\left\{n\right\}$, $n\neq 4k+2$ for $k\geq 1$,
then $W_\Gamma$ is profinitely rigid amongst Coxeter groups.
\end{thm}
\begin{proof}
It is known that a connected simplicial graph with $4$ vertices is isomorphic to one of the graphs in \Cref{fig:connectedgraphs4vertices}. 

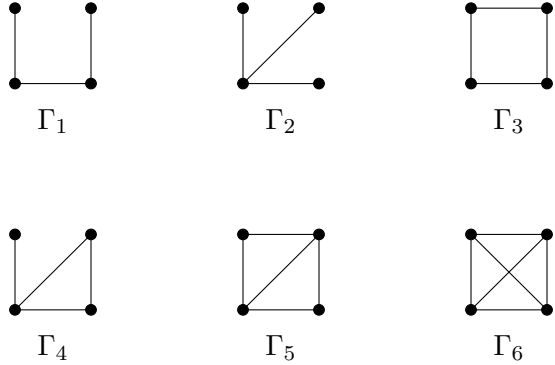
\begin{figure}[htb]
	\begin{center}
	\captionsetup{justification=centering}
		\begin{tikzpicture}
			\draw[fill=black]  (0,0) circle (2pt);
            \draw[fill=black]  (0,1) circle (2pt);
            \draw[fill=black]  (1,0) circle (2pt);
            \draw[fill=black]  (1,1) circle (2pt);
            \draw (0,0)--(1,0);
            \draw (0,1)--(0,0);
            \draw (1,0)--(1,1);
            \node at (0.5, -0.5) {$\Gamma_1$};

            \draw[fill=black]  (3,0) circle (2pt);
            \draw[fill=black]  (3,1) circle (2pt);
            \draw[fill=black]  (4,0) circle (2pt);
            \draw[fill=black]  (4,1) circle (2pt);
            \draw (3,0)--(4,0);
            \draw (3,0)--(4,1);
            \draw (3,0)--(3,1);
            \node at (3.5, -0.5) {$\Gamma_2$};

            \draw[fill=black]  (6,0) circle (2pt);
            \draw[fill=black]  (6,1) circle (2pt);
            \draw[fill=black]  (7,0) circle (2pt);
            \draw[fill=black]  (7,1) circle (2pt);
            \draw (6,0)--(6,1);
            \draw (6,0)--(7,0);
            \draw (7,0)--(7,1);
            \draw (6,1)--(7,1);
            \node at (6.5, -0.5) {$\Gamma_3$};          

            \draw[fill=black]  (0,-2) circle (2pt);
            \draw[fill=black]  (0,-3) circle (2pt);
            \draw[fill=black]  (1,-2) circle (2pt);
            \draw[fill=black]  (1,-3) circle (2pt);
            \draw (0,-3)--(1,-3);
            \draw (0,-3)--(1,-2);
            \draw (0,-3)--(0,-2);
            \draw (1,-3)--(1,-2);
            \node at (0.5, -3.5) {$\Gamma_4$};

            \draw[fill=black]  (3,-2) circle (2pt);
            \draw[fill=black]  (3,-3) circle (2pt);
            \draw[fill=black]  (4,-2) circle (2pt);
            \draw[fill=black]  (4,-3) circle (2pt);
            \draw (3,-3)--(4,-3);
            \draw (3,-3)--(4,-2);
            \draw (3,-3)--(3,-2);
            \draw (3,-2)--(4,-2);
            \draw (4,-3)--(4,-2);
            \node at (3.5, -3.5) {$\Gamma_5$};

            \draw[fill=black]  (6,-2) circle (2pt);
            \draw[fill=black]  (6,-3) circle (2pt);
            \draw[fill=black]  (7,-2) circle (2pt);
            \draw[fill=black]  (7,-3) circle (2pt);
            \draw (6,-3)--(7,-3);
            \draw (6,-3)--(7,-2);
            \draw (6,-3)--(6,-2);
            \draw (7,-3)--(7,-2);
            \draw (7,-3)--(6,-2);
            \draw (6,-2)--(7,-2);
            \node at (6.5, -3.5) {$\Gamma_6$};
                      
        \end{tikzpicture}
        \caption{Connected graphs with four vertices.}
	      \label{fig:connectedgraphs4vertices}
    \end{center}
\end{figure}

Let $W_\Omega$ be a Coxeter group such that $\widehat{W_\Gamma}\cong \widehat{W_\Omega}$. 

\begin{enumerate}
\item If $\Gamma\cong \Gamma_1$ or $\Gamma_2$, then by a characterisation of virtually free Coxeter groups via combinatorial properties of the defining graphs \cite[Theorem 34]{MihalikTschantz2009} $W_\Gamma$ is virtually free. Since being virtually free in a $\mathcal{W}$-profinite invariant, see \Cref{{profinite.invariant. virtuallyfreeness}} we know that $W_\Omega$ is also virtually free  and therefore we have order isomorphisms by \Cref{prop:CFhyperbolicFC}.
$$\CF(W_\Gamma)=\CF(\widehat{W_\Gamma})=\CF(\widehat{W_\Omega})=\CF(W_\Omega).$$
\begin{enumerate}
\item If $n$ is odd, then $W_\Gamma\cong W_\Omega$ by \Cref{oddForest_rigid}.
\item If $n$ is even, then we consider two cases: (i) if $n=2$, then $W_\Gamma$ is a right-angled Coxeter group and so it is $\mathcal{W}$-profinitely rigid by \cite[Theorem 3.8]{CorsonHughesMollerVarghese2023}.

(ii) if $n\geq 4$, then it follows by \Cref{thm:even.profinite} that $W_\Gamma\cong W_\Omega$. 
\end{enumerate}

\item If $\Gamma\cong \Gamma_3$, then $W_\Gamma\cong W_\Omega$ by \Cref{triangleAndmore}. 

\item If $\Gamma\cong \Gamma_4$ or $\Gamma_5$, then by \Cref{amalgam.profinitelyrigid} it follows that $W_\Gamma\cong W_\Omega$, since triangle groups and Dihedral groups $D_n$, where $n\neq 4k+2$ are graph rigid and $\mathcal{W}$-profinitely rigid by \Cref{triangleAndmore}.

\item Let $\Gamma\cong \Gamma_6$. If $n$ is odd, then $W_\Gamma\cong W_\Omega$ by \Cref{OddComplete.profinite}. 
Let $n$ be even, $n\neq 4k+2$ for $k\geq 1$.  Then $W_\Gamma\cong W_\Omega$ by \Cref{evenComplete.profinite}.
\end{enumerate}
This completes the proof.
\end{proof}

\bigskip\bigskip
\appendix
\section{\texorpdfstring{$\ell^2$}{ℓ²}-invariants of Coxeter groups}\label{sec appendix}
\smallskip
\begin{center}by \textsc{Sam P.~Fisher} and \textsc{Sam Hughes}\end{center}
\medskip

\subsection{Background on \texorpdfstring{$\ell^2$}{ℓ²}-invariants and the Atiyah conjecture}

We give a very brief review on $\ell^2$-invariants. The reader is invited to consult \cite{Luck02,KammeyerBook2019} for more details. Let $G$ be a countable discrete group. It acts by left multiplication on $\ell^2(G)$, the Hilbert space of square summable complex series supported in $G$. The algebra of bounded operators acting on $\ell^2(G)$ on the right and commuting with the left $G$-action is denoted by $\mathcal N(G)$ and called the \emph{von Neumann algebra} of $G$. Note that the complex group algebra $\CC G$ is naturally a subring of $\mathcal N(G)$. The non-zero divisors of $\mathcal N(G)$ form an Ore set, and the localisation at the set is denoted by $\mathcal U(G)$ and called the \emph{algebra of operators affiliated to $G$}. There is a dimension function $\dim_{\mathcal U(G)}$ on $\mathcal U(G)$-modules; we do not define it here.

\begin{defn}
    Let $G$ be a countable group. The \emph{$\ell^2$-Betti numbers} of $G$ are given by
    \[
        b_i^{(2)}(G) = \dim_{\mathcal U(G)} \operatorname{Tor}_i^{\CC G}(\mathcal U(G), \CC) = \dim_{\mathcal U(G)} \operatorname{H}_i (G;\mathcal U(G)).
    \]
\end{defn}

The Atiyah conjecture for a group $G$ concerns the possible values of the dimension function $\dim_{\mathcal U(G)}$.

\begin{conjecture}[The strong Atiyah conjecture]
    Let $G$ be a countable group such that there is a bound on the orders of finite subgroups. Let $\mathrm{lcm}(G)$ be the least common multiple of the orders of the finite subgroups of $G$. If $k \subseteq \CC$ is a subfield, we say that $G$ satisfies the \emph{strong Atiyah conjecture over $k$} if
    \[
        \dim_{\mathcal U(G)} (\mathcal U(G) \otimes_{kG} M) \in \frac{1}{\mathrm{lcm}(G)} \Z
    \]
    for every finitely presented left $kG$-module $M$.
\end{conjecture}

\subsection{The Atiyah conjecture for Coxeter groups}

Following Schreve \cite{Schreve2014}, we say that a group $G$ has the \emph{factorisation property} if every epimorphism of the form $G \rightarrow Q$ with $Q$ finite factors through a map $G \rightarrow E$ where $E$ is torsion-free and elementary amenable. By using the work of Schreve \cite{Schreve2014} and Genevois \cite{Genevois2024}, we can prove that Coxeter groups virtually have the factorisation property by arguing similarly as in \cref{CoxGpGood}.

\begin{lemma}\label{lem:Coxeter_factorisation}
    If $W_\Gamma$ is a Coxeter group, then there is a finite-index subgroup $H \leqslant W_\Gamma$ with the factorisation property.
\end{lemma}
\begin{proof}
    If $W_\Gamma$ is a right-angled Coxeter group, then it is virtually compact special since it acts cocompactly on its Davis complex, and thus $W_\Gamma$ virtually has the factorisation property by \cite[Corollary~4.3]{Schreve2014}. Now suppose that $W_\Gamma$ is a general Coxeter group. By \cref{cor:Cox_vir_retract_RACG}, there is a finite-index subgroup $H_0 \leqslant W_\Gamma$ that is a retract of a finite-index subgroup $G_0$ of some larger right-angled Coxeter group $W_\Omega$. Let $G_1 \leqslant W_\Omega$ be a finite-index subgroup with the factorisation property. Then $G := G_0 \cap G_1$ has the factorisation property \cite[Lemma~2.1]{Schreve2014} and it retracts onto its image $H$, which is of finite-index in $H_0$. Since the factorisation property passes to retracts \cite[Lemma~2.2]{Schreve2014}, $H$ has the factorisation property. \qedhere
\end{proof}

We can now establish the strong Atiyah conjecture for general Coxeter groups, extending the results of \cite{LinnellOkunSchick_atiyahCoxeter} and \cite[Corollary~4.5]{Schreve2014} where it is proven in the right-angled case.

\begin{corollary}
    Coxeter groups satisfy the strong Atiyah conjecture over $\CC$. More generally, if $1 \rightarrow W_\Gamma \rightarrow G \rightarrow E \rightarrow 1$ is a short exact sequence where $W_\Gamma$ is a Coxeter group and $E$ is elementary amenable, then $G$ satisfies the strong Atiyah conjecture over $\CC$.
\end{corollary}
\begin{proof}
    By \cref{CoxGpGood} and \cref{lem:Coxeter_factorisation}, we can choose a characteristic finite-index subgroup $H \leqslant W_\Gamma$ that is good and has the factorisation property. Then $H$ is normal in $G$ and the quotient $G/H$ is still elementary amenable, being an extension of a finite group by $E$. It then follows from \cite[Theorem 1.1]{Schreve2014} that $G$ satisfies the strong Atiyah conjecture over $\CC$. \qedhere
\end{proof}

\subsection{Profinite invariance of \texorpdfstring{$\ell^2$}{ℓ²}-Betti numbers}

In this section, we will focus on the class of residually (locally indicable amenable) groups, since this is essentially the most general class where our arguments work. Note that Agol's RFRS groups \cite{AgolCritVirtFib} are residually (locally indicable and Abelian), so they fall within this class. Thus, typical examples of groups to which the following results apply are (virtually) special groups.

We begin by briefly introducing positive characteristic analogues of $\ell^2$-Betti numbers. Linnell \cite{LinnellDivRings93} showed that if $G$ is torsion-free, then $G$ satisfies the strong Atiyah conjecture over $k \subseteq \CC$ if and only if the division closure of $kG$ in $\mathcal U(G)$ is a division ring. If $G$ is locally indicable, then $G$ satisfies the strong Atiyah conjecture over $\CC$, and therefore the division closure of $\CC G$ in $\mathcal U(G)$ is the unique \emph{Hughes-free division ring} containing $\CC G$ and is denoted by $\mathcal D_{\CC G}$ (see \cite{JaikinZapirain2020Div} for a definition of the Hughes-free property and more background). The $\ell^2$-Betti numbers of $G$ can be calculated using
\[
    b_i^{(2)}(G) = \dim_{\mathcal D_{\CC G}} \operatorname{Tor}_i^{\CC G}(\mathcal D_{\CC G}, \CC),
\]
where $\dim_{\mathcal D_{\CC G}}$ returns the usual rank of a (necessarily free) $\mathcal D_{\CC G}$-module.

If $G$ is a locally indicable group and $\FF$ is a field, then in some (conjecturally all) cases there exists a Hughes-free division ring $\mathcal D_{\FF G}$ containing $\FF G$, and if it exists it is unique by \cite{HughesDivRings1970}. The $\ell^2$-Betti numbers of $G$ over $\FF$ are defined by
\[
    b_i^{(2)}(G;\FF) = \dim_{\mathcal D_{\FF G}} \operatorname{Tor}_i^{\FF G}(\mathcal D_{\FF G}, \FF),
\]
in complete analogy with the classical $\ell^2$-Betti numbers. Jaikin-Zapirain showed that if $G$ is a residually (locally indicable amenable) group, then $\mathcal D_{\FF G}$ always exists, and satisfies an additional property called \emph{universality} \cite[Corollary~1.3]{JaikinZapirain2020Div}.

The following result gives a weak analogue of L\"uck approximation in positive characteristic. An essentially equivalent statement has already appeared in \cite[Theorem 3.6]{AvramidiOkunSchreve_7manifold} in the case of residually (torsion-free nilpotent) groups, and the proof is almost identical. For an $\FF G$-module $M$, we write $\dim_{\mathcal D_{\FF G}}(M)$ for $\dim_{\mathcal D_{\FF G}}(\mathcal D_{\FF G} \otimes_{\FF G} M)$.

\begin{lemma}\label{lem:approx_fp_module}
    Let $G$ be a residually (locally indicable and amenable) group. If $M$ is a finitely presented $\FF G$-module, then
    \[
        \dim_{\mathcal D_{\FF G}}(M) = \inf_{G_0 \leqslant_{\mathrm{f.i.}} G} \left\{ \frac{\dim_\FF(\FF \otimes_{\FF G_0} M)}{[G:G_0]} \right\},
    \]
    where the infimum is taken over all finite-index subgroups of $G$.
\end{lemma}
\begin{proof}
    Let $G_0 \leqslant G$ be an arbitrary subgroup of finite-index. Then
    \[
        \dim_{\mathcal D_{\FF G}}(M) = \frac{\dim_{\mathcal D_{\FF G_0}}(M)}{[G:G_0]} \leqslant \frac{\dim_\FF(\FF \otimes_{\FF G_0} M)}{[G:G_0]}
    \]
    where the equality follows from Hughes-freeness of $\mathcal D_{\FF G}$ and the inequality follows from the universality of $\mathcal D_{\FF G_0}$ (see \cite[Corollary~1.3]{JaikinZapirain2020Div}).

    To prove the reverse inequality, let $\FF G^m \xrightarrow{A} \FF G^n \to M \to 0$ be a finite presentation of $M$, where $A$ is a matrix with entries in $\FF G$. By \cite[Theorem 1.2]{JaikinZapirain2020Div}, there is a normal subgroup $N \trianglelefteqslant G$ such that $Q = G/N$ is amenable and locally indicable and $\operatorname{rk}_{\mathcal D_{\FF G}}(A) = \operatorname{rk}_{\mathcal D_{\FF Q}}(A)$. By \cite[Theorem 0.2]{LinnellLuckSauer_amenableApprox},
    \begin{align*}
        \dim_{\mathcal D_{\FF G}}(M) &= \dim_{\mathcal D_{\FF Q}}  (\FF Q \otimes_{\FF G} M) \\
        &= \lim_{i \to \infty} \frac{\dim_\FF(\FF \otimes_{\FF Q_i} (\FF Q \otimes_{\FF G} M))}{[Q:Q_i]} \\
        &= \lim_{i \to \infty} \frac{\dim_\FF(\FF \otimes_{\FF G_i} M)}{[Q:Q_i]},
    \end{align*}
    where $Q \geqslant Q_1 \geqslant Q_2 \geqslant \dots$ is an arbitrary residual normal chain of finite index subgroups of $Q$ and $G_i$ is the preimage of $Q_i$ in $G$ (we have used the fact that $\FF Q \otimes_{\FF G} M$ is a finitely presented $\FF Q$-module). Choosing a sufficiently deep finite-index subgroup $G_i \leqslant G$ gives the reverse inequality. \qedhere
\end{proof}

Let $R$ be a commutative ring with unity.  A group $G$ is \emph{type $\mathsf{FP}_n(R)$} if there exists a projective resolution $P_\bullet\to R$ of the trivial $RG$-module $R$ such that $P_i$ for $i\leq n$ is finitely generated.

\begin{corollary}\label{cor:betti_approx}
    Let $G$ be a residually (locally indicable and amenable) group of type $\mathsf{FP}_{n+1}(\FF)$ for some field $\FF$. Then 
    \[
        b_i^{(2)}(G;\mathbb F) = \inf_{G_0 \leqslant_{\mathrm{f.i.}} G} \left\{ \frac{b_i(G_0;\FF)}{[G:G_0]} \right\}
    \]
    for all $i \leqslant n$.
\end{corollary}
\begin{proof}
    Let $\cdots \rightarrow C_1 \rightarrow C_0 \rightarrow \FF \rightarrow 0$ be a free resolution of the trivial $\FF G$-module $\FF$ such that $C_i$ is finitely generated for all $i \leqslant n+1$. Denote the boundary maps by $\partial_i \colon C_i \rightarrow C_{i-1}$. Fix some $i \leqslant n$ and let $r$ be the rank of $C_i$. Then
    \begin{align*}
        b_i^{(2)}(G;\mathbb F) &= \dim_{\mathcal D_{\FF G}} \operatorname{H}_i(\mathcal D_{\FF G} \otimes_{\FF G} C_\bullet) \\
        &= \dim_{\mathcal D_{\FF G}}(\operatorname{coker}(\partial_i)) + \dim_{\mathcal D_{\FF G}}(\operatorname{coker}(\partial_{i+1})) - r \\
        &= \inf_{G_0 \leqslant_{\mathrm{f.i.}} G} \left\{ \frac{\dim_\FF(\operatorname{coker}(\partial_i)) + \dim_\FF(\operatorname{coker}(\partial_{i+1})) - [G:G_0]r}{[G:G_0]} \right\} \\
        &= \inf_{G_0 \leqslant_{\mathrm{f.i.}} G} \left\{ \frac{b_i(G_0;\FF)}{[G:G_0]} \right\},
    \end{align*}
    where $\dim_\FF(M) = \dim_\FF(\FF \otimes_{\FF G} M)$. Note that the third equality does not formally follow from the statement of \cref{lem:approx_fp_module}, but from the proof. More precisely, we use that there is a single locally indicable amenable quotient $Q$ of $G$ such that $\dim_{\mathcal D_{\FF G}}(\operatorname{coker}(\partial_j)) = \dim_{\mathcal D_{\FF Q}}(\operatorname{coker}(\partial_j))$ for $j = i, i+1$ and the fact that L\"uck approximation holds for locally indicable amenable groups \cite[Theorem 0.2]{LinnellLuckSauer_amenableApprox}. \qedhere
\end{proof}

Taking $\FF = \QQ$ in the below establishes the profinite invariance of the classical $\ell^2$-Betti numbers of residually (locally indicable amenable) groups; note that this contrasts sharply with \cite{KammeyerKionkeRaimbaultSauer2020}. Note that this vastly generalises \cite[Theorem~5.11]{HughesKielak2025}.

\begin{thm}\label{thm:L2_Betti_profinite}
    Let $G$ and $H$ be $n$-good virtually residually (locally indicable amenable) groups of type $\mathsf{FP}_{n+1}$. If $\widehat G \cong \widehat H$, then, for each field $\FF$,
    \[
        b_i^{(2)}(G;\FF) = b_i^{(2)}(H;\FF)
    \]
    for all $i \leqslant n$.
\end{thm}
\begin{proof}
    We may assume that $G$ and $H$ are both residually (locally indicable and amenable). Indeed, if they are not pass to residually (locally indicable amenable) finite-index subgroups of $G$ and $H$ that are profinitely isomorphic and of the same index in $G$ and $H$; this preserves goodness and if the theorem holds for finite-index subgroups of the same index then it also holds for $G$ and $H$ by \cite[Lemma~6.3]{Fisher_Improved}.
    
    By \cref{cor:betti_approx}, the quantity $b_i^{(2)}(G;\FF)$ depends only on the characteristic of $\FF$. Thus, it is enough to prove the theorem for prime fields. First assume that $\FF$ is a finite field. Then
    \begin{multline*}
        b_i^{(2)}(G;\mathbb F) = \inf_{G_0 \leqslant_{\mathrm{f.i.}} G} \left\{ \frac{b_i(G_0;\FF)}{[G:G_0]} \right\} = \inf_{G_0 \leqslant_{\mathrm{f.i.}} G} \left\{ \frac{b_i(\widehat{G_0};\FF)}{[G:G_0]} \right\} \\
        = \inf_{H_0 \leqslant_{\mathrm{f.i.}} H} \left\{ \frac{b_i(\widehat{H_0};\FF)}{[H:H_0]} \right\} = \inf_{H_0 \leqslant_{\mathrm{f.i.}} H} \left\{ \frac{b_i(H_0;\FF)}{[H:H_0]} \right\} = b_i^{(2)}(H;\mathbb F)
    \end{multline*}
    for all $i \leqslant n$, by goodness and two applications of \cref{cor:betti_approx}. 

    The case $\FF = \QQ$ follows immediately since, for a sufficiently large prime $p$, 
    \[
        b_i^{(2)}(G;\mathbb Q) = b_i^{(2)}(G;\mathbb F_p) = b_i^{(2)}(H;\mathbb F_p) = b_i^{(2)}(H;\mathbb Q).
    \]
    This is explained in \cite[Corollary~4.2]{AvramidiOkunSchreve_7manifold} for virtually residually (torsion-free nilpotent) fundamental groups of compact CW complexes, but the argument applies in the generality of this result as well. \qedhere
\end{proof}

A rich source of examples of groups satisfying the hypotheses of \cref{thm:L2_Betti_profinite} are virtually \emph{compact} special groups, since they are virtually residually (locally indicable amenable) and good in the sense of Serre (see, e.g., \cite[Corollary~4.3]{Schreve2014}). 

We record some applications to good (virtually) RFRS groups.  Let $R$ be a commutative ring with unity.  A group $G$ is \emph{$\mathsf{FP}_n(R)$-fibred} if there exists an epimorphism $\phi\colon G\onto \Z$ such that $\ker \phi$ is type $\mathsf{FP}_n(R)$.

\begin{corollary}\label{cor:RFRS:fibring:prof}
    Let $G$ and $H$ be good virtually RFRS groups of type $\mathsf{FP}_{n+1}(\FF)$ for some field $\FF$ and suppose that $\widehat{G} \cong \widehat{H}$. Then $G$ is virtually $\mathsf{FP}_n(\FF)$-fibred if and only if $H$ is.
\end{corollary}
\begin{proof}
    This is an immediate consequence of \cite[Theorem B]{Fisher_Improved} and \cref{thm:L2_Betti_profinite}. \qedhere
\end{proof}

\begin{corollary}
    The properties of being virtually free-by-cyclic and virtually (finitely generated free)-by-cyclic are both profinite invariants in the class of good RFRS groups of cohomological dimension at most $2$.
\end{corollary}
\begin{proof}
    The profinite invariance of being virtually free-by-cyclic is an immediate consequence of \cite[Theorem A]{Fisher_freebyZ} (see \cite{KielakLinton2024} for the case virtually compact special hyperbolic groups) and \cref{thm:L2_Betti_profinite}.

    Now suppose that $G$ is a good RFRS group of cohomological dimension at most $2$, and suppose that it is profinitely isomorphic to a RFRS virtually $F_n$-by-$\Z$ group. Then $b_1^{(2)}(G) = b_2^{(2)}(G) = 0$, so $G$ is virtually $\mathrm{FP}$-fibred by \cite[Theorem 5.4]{KielakRFRS}. By \cite[Theorem 2.4]{Feldman71}, this implies that the kernel of the virtual fibration must in fact be a finitely generated free group. \qedhere
\end{proof}

Finally, we conclude with the applications to Coxeter groups.

\begin{thm}\label{app:cor.Cox}
    If $W_\Gamma$ and $W_\Lambda$ are Coxeter groups such that $\widehat{W_\Gamma} \cong \widehat{W_\Lambda}$, then
    \begin{enumerate}
        \item $b_i^{(2)}(W_\Gamma;\FF) = b_i^{(2)}(W_\Lambda;\FF)$ for all $i \geqslant 0$;
        \item and $W_\Gamma$ is virtually $\mathsf{FP}_n(\FF)$-fibred if and only if $W_\Lambda$ is.
    \end{enumerate}
\end{thm}
\begin{proof}
    To see (1), observe that  Coxeter groups satisfy the hypothesis of \Cref{thm:L2_Betti_profinite} by \cite{HaglundWise2010}.  To see (2), we note that Coxeter groups are virtually RFRS, either by Haglund--Wise \cite{HaglundWise2010}, or by Genevois \cite{Genevois2024} and it is classical that they are type $\mathsf{F}_\infty$.  Thus, we may apply \cref{cor:RFRS:fibring:prof}.
\end{proof}

\bibliographystyle{halpha}
\bibliography{refs.bib}
\bigskip
\end{document}